# Mathematical Analysis of the Problems faced by the People With Disabilities (PWDs)

## (With Specific Reference to Tamil Nadu in India)

**W. B. Vasantha Kandasamy**
**Florentin Smarandache**
**A. Praveen Prakash**

**2012**

# Mathematical Analysis of the Problems faced by the People With Disabilities (PWDs)

(With Specific Reference to Tamil Nadu in India)

**W. B. Vasantha Kandasamy**
**Florentin Smarandache**
**A.  Praveen Prakash**

# CONTENTS

# PREFACE

The authors in this book have analyzed the socio-economic and psychological problems faced by People with Disabilities (PWDs) and their families. The study was made by collecting data using both fuzzy linguistic questionnaire / by interviews in case they are not literates from 2,15,811 lakhs people. This data was collected using the five Non Government Organizations (NGOs) from northern Tamil Nadu.

Now any reader would be interested to know whether the Tamils (natives of Tamil Nadu) had ever spoken about people with disability. Even before 2000 years tamils had heroic poetry Purananuru (28[th] poem) about the war fare methods. In that poem mention was made about 8 types of disability.

The 8 types of disability (by birth not due to disease or accidents or due to war) mentioned in order are as follows: (1) Blindness, (2) Aborted embryo (3) Hump back (lying with the body disposed in a crooked position) (4) Dwarf (5) Dumbness (speech impaired) (6) Deafness, (hearing impaired) (7) No proper development of bones (8) Bewilderment of mind (Not in a position to standup these persons mind never grows) (reference Purananuru poem 28).

Further the tamils had the practice; if the child was born dead (aborted embryo) or dies after birth they had the custom of cutting it by a sword and burying it. This was only to show they had the custom of cutting it by a sword and burying it. This was only to show that they had the valour and the war spirit in their everyday life. (ref. Purananuru 74).

Thus the tamils in ancient days had mentioned about the disabilities and also had ordered them.

Now we wish to talk of disability in the integrated India. India has two great epics Ramayana and Mahabaratha. Also India has the Laws of Manu which speaks of the code and conduct for the kings and the people in general.

The Laws of Manu claims that the lame or single eyed persons while doing prayer for the dead should not present even if they are their close relatives (Laws of Manu chapter VI, 242 sloka) [5].

When the king takes important decisions with ministers it was strictly followed that persons who are blind, dumb, mentally disabled or deaf or any PWD should be driven out of vicinity (Manu chapter VII, sloka 149).

Blind or dumb or deaf or lame or mentally retarded cannot inherit property (chapter IX in Laws of Manu, 209 sloka) [5].

In keeping with this even in the great old epic Mahabaratha the elder brother Dirdhirashtra who was blind was denied the kingship and his younger brother Pandu was made the king.

Likewise in the great Hindu epic Ramayana the PWD's were portrayed to be very cruel and never were good at heart or mind (Kuni the servant in the palace of Rama is an example of it).

Both the kings and Brahmins can take over the property of the sinners and those kings and Brahmins will never give birth to PWDs (Chapter IX Laws of Manu, sloka 247) [5]. So it was believed that those who enjoyed the properties of sinners will give birth to children with disability. If in a Brahmin there was a mentally retarded one he was considered to be better than other 4 castes and dalits (Chapter IX Laws of Manu, sloka 317) [5].

Laws of Manu which is supposed to be the oldest which gave laws and conduct was written by Brahmins and the first four chapters are only about them and it was they who gave the four types of castes and the notion of untouchability.

The Laws of Manu were the ancient laws but followed even today by certain upper castes (Brahmins). Further the four caste-system mentioned in the Laws of Manu are in vogue even today. This has made a lot of social and psychological problems among the natives of India. Further the Laws of Manu had been written by the Brahmins to see they enjoyed the maximum benefit of the land where they came into as nomads. Further the Brahmins not only divided the land into four castes they created lower castes please refer chapter X, slokas 1 to 67 of Laws of Manu [5].

Here it is clearly brought out the only the Laws of Manu created the lowest castes and these caste people were made so due to intercaste marriages and these people must dwell outside the village and their wealth shall be dogs and donkeys (Chapter X, sloka 51). They should wear the garments of the dead and eat food from broken dishes and at night they shall not walk about in villages and in towns (Chapter X, sloka 52 to 54). Their work is to carry out corpses (of persons) who have no relatives (Chapter X, sloka 55). Now we have just mentioned a word or two about the disability in a disabled caste system given by the Laws of Manu [5].

Now we just mention about how the social and psychological problems of the PWDs were studied and analysed using mathematical and fuzzy linguistic models. At the outset the authors who are not experts about the technicality of naming the disability may mention the type of disability in a common way say as 'blind' instead of "impaired vision" or so, likewise 'deaf' instead of hearing impaired and so on. This is by no means to offend or make one feel bad, it is completely the inability of the authors to use proper technical terms.

Secondly the PWDs most of them physical is analysed using statistical methods and fuzzy models by the authors.

The mental retardation or mental disability is analysed using fuzzy linguistic models. We have categorized the mental disability into five types mainly for our fuzzy linguistic working and they are not in any way medically or technically assigned. The data was collected from interviews from the caretakers of the mentally disabled children / adults from the homes run by NGO's / charity organization in / near Chennai. Then we used experts, the experts were caretakers or parents or persons closely associated with the mentally disabled people.

The suggestion / conclusions are from the experts and we are not very sure whether it is true or false for it is what they feel or they have experienced while taking care of these children / adults this only a firsthand experience who have been with these mentally retarded children for over a decade or so. Most of the information gathered is from their practical knowledge. Also when we enquired a few experts in the field some of them said the mental disability cannot be easily detected by birth and lots of research ought to be done in this direction and they also said it was incurable. However our study is mainly the social and the psychological problems faced by the parents / relatives of the mentally retarded / PWD, as even today the birth of such children is only thought of as a curse so a curse on the family, so public shun them and their family and they are socially stigmatized. "How to overcome this situation?" is our study.

This book has eight chapters. Chapter one just gives references to the materials used in this book. Analysis of the collected date is carried out in chapter two. Chapter three introduces a New Fuzzy Cognitive Relational Maps (FCRMs) bimodel to study the socio-economic problems of PWDs.

Economic problems faced by PWDs is analysed in chapter four using Special Fuzzy Cognitive Maps (SFCMs) model introduced by the authors. Chapter five carries out the study of interrelated problems faced by the PWDs and the caretakers. In chapter six we use the new FCRM bimodel to analyse the problems faced by the rural PWDs.

Since the study of mental retardation or mental disability cannot be studied using the opinion of the disabled we have to approach the problem by a new methods. For most caretakers of these mentally disabled people / children were not literate or had no time to sit and fill questionnaire so their answers to questions was mainly linguistic in nature, so the authors opted to used fuzzy linguistic models to analyse and study the problem. Further this data was collected from charity homes run by NGOs and others in and near Chennai is carried out in chapter seven. Conclusions based on this study are given the end of this chapter.

The final chapter gives the suggestions and conclusions from our study in chapter two to six and also from the interviews and discussions with them.

We also thank NGOs; ROSE, LAMP, SERVITES, VDS, MAKAM, WCDS, SHARE and CARE etc. We also like to thank Mr.A.John Aruldoss for helping in data collection. Finally we thank Dr. K.Kandasamy for various unostentatious help and suggestions while collecting the data of the mentally disabled people and for proof reading, without whose cooperation and support this book would have not seen this form.

W.B.VASANTHA KANDASAMY

FLORENTIN SMARANDACHE

A. PRAVEEN PRAKASH

Chapter One

# INTRODUCTION TO BASIC CONCEPTS

In this chapter we just indicate the references that can supply the necessary material for the researcher / reader to follow the fuzzy mathematical models described in this book in an unterrupted way. Here we have collected the data using linguistic questionnaire and in the second chapter they are only analysed.

In chapter two we use a new fuzzy bimodel called the, New Fuzzy Cognitive Relational Maps (FCRMs) bimodel. So the reader should be familiar with the concept of Fuzzy Cognitive Maps (FCMs). For this concept please refer [19-21, 44].

FCM models have been studied by several authors [19-21, 44]. Further we also use the Fuzzy Relational Maps (FRM) model to build the new FCRMs bimodel. Fuzzy relational maps can be realized as the generalization / simplification / modification of the fuzzy cognitive maps for these function with a domain space and a range space and in case of FRMs the hidden pattern of the dynamical system for any given state vector is pair of row vectors, which may be a fixed point or a limit cycle. Further the matrix associated with the dynamical system is a rectangular matrix. This model is applicable when the concepts / attributes can be divided into two disjoint sets and there exists no relation among the attributes in the range space (or among the attributes in the domain space). For more about this concept please refer [19-21, 44]. The new FCRMs bimodel is built using both the FCM model and FRM model. The advantage of using this and the functioning of it will be described in chapter three of this book.

Further this book for the first time uses the Fuzzy Linguistic Cognitive Maps (FLCM) model and Fuzzy Linguistic Relational Maps (FLRM) model to study the problems faced by the mentally retarded children. For basic notions about this model please refer [50]. This new model functions just like the fuzzy model instead of taking the membership values from the unit interval [0, 1] we take the membership values from the fuzzy linguistic set L where L consists of only fuzzy linguistic terms like {fast, very fast, slow, very slow, medium fast, medium slow, ..., 0}.

It can be measuring the speed of a car / any vehicle on the road. Likewise one can measure the performance of students in a class; say the set L can be as follows: L = {good, 0, best, average, very good, fair, very fair, bad, very bad, worst}. This set L is replaced by [0, 1] and now the graphs associated with these fuzzy linguistic models will be known as the fuzzy linguistic graphs.

They are also directed graphs the vertices are the attributes associated with the problem and the edges are the linguistic terms from L. Now the associated matrix of this fuzzy linguistic graphs are matrices whose entries are from the set L so instead of putting values from [0, 1] or from {−1, 0, 1} we use the entries from L.

These matrices are defined as fuzzy linguistic matrices. Now we will be using these matrices as the dynamical system to arrive at the fuzzy linguistic hidden pattern which will be a fixed point or a limit cycle.

For more about these new models please refer [50]. These models are used in the chapter seven of this book.

Chapter Two

# ANALYSIS OF THE COLLECTED DATA OF THE PWDS FROM MELMALAYANUR AND KURINJI PADI BLOCKS IN TAMIL NADU

We have collected data using linguistic questionnaire / by interviews, in case they are not literates, from 2,15,811 lakhs people. The data was collected using the five Non Government Organizations (NGOs) listed in the following table. The blocks and the number of village panchayats from which the data has been collected are also listed. The data is tabulated in the order, which is mainly in keeping with the priority of the PWDs feelings expressed from the interviews. The conclusions based on these interviews are given in the end of this chapter. However the physical interpretation of the tables are given then and there.

The main objective of this was to understand the basic psychological and social problems faced by the PWDs and their caretakers. To interpret the collected data, tabulations and analysis were made. Efforts were made to highlight the inferences in each table and in few places comments were made in order to use it for fuzzy analysis. Following is the list of organizations and the total population from 5 blocks is listed below.

**NGO wise distribution of the PWDs in 93 Panchayats in Melmalayanur and Kurinjipadi Blocks**

| Sl. No. | Name of Implementing Organisation | Target Block | No. of Village Panchayats | Total Population | | |
|---|---|---|---|---|---|---|
| | | | | M | F | Total |
| 1 | ROSE | Kurinjipadi | 22 | 20,650 | 20,312 | 40,962 |
| 2 | LAMP | Kurinjipadi | 16 | 22,368 | 21,190 | 43,558 |
| 3 | SERVITES | Mel-malayanur | 15 | 18,280 | 18,205 | 36,485 |
| 4 | VDS | Mel-malayanur | 20 | 23,051 | 22,382 | 45,433 |
| 5 | MAKAM | Mel-malayanur | 20 | 24,994 | 24,379 | 49,373 |
| Total Population | | | 93 | 109,343 | 106,468 | 215,811 |
| Total Persons With Disability | | | 93 | 2,115 | 1,424 | 3,539 |
| % of Total PWDs on Total Population | | | | 1.93% | 1.34% | 1.64% |

From the above table it is clear that the study was conducted in 93 Panchayats in a population of 2,15,811. The overall population of the PWDs on the total population is 1.64 per cent, which is lesser than the national average. Here, with regard to the PWDs population, the male population is comparatively higher. The table also depicts the coverage of Panchayats and population by each implementing organisation both in Melmalayanur and Kurinjipadi Blocks of Tamilnadu state in India. Now we proceed onto give the age wise physical interpretation of the data using the following table.

## The age-group wise stratification of PWDs

| Sl. No. | Age Group | Total PWDs | | | Types of PWDs | | |
|---|---|---|---|---|---|---|---|
| | | | | | Ortho | | |
| | | M | F | Total | M | F | Total |
| 1 | Below 3 | 18 | 13 | 31 | 5 | 7 | 12 |
| 2 | 4 to 6 | 76 | 66 | 142 | 23 | 12 | 35 |
| 3 | 7 to 14 | 184 | 157 | 341 | 47 | 52 | 99 |
| 4 | 15 to 18 | 136 | 103 | 239 | 55 | 44 | 99 |
| 5 | 19 to 25 | 315 | 256 | 571 | 203 | 139 | 342 |
| 6 | 26 to 35 | 372 | 248 | 620 | 229 | 133 | 362 |
| 7 | 36 to 45 | 332 | 210 | 542 | 203 | 104 | 307 |
| 8 | 46 to 60 | 479 | 282 | 761 | 278 | 138 | 416 |
| 9 | 61 & above | 203 | 89 | 292 | 87 | 29 | 116 |
| Total | | 2,115 | 1,424 | 3,539 | 1,130 | 658 | 1,788 |
| Percentage | | 60 | 40 | 100 | 32 | 19 | 51 |

| Types of PWDs | | | | | | | | |
|---|---|---|---|---|---|---|---|---|
| CP | | | MR | | | LC | | |
| M | F | Total | M | F | Total | M | F | Total |
| 1 | 2 | 3 | 7 | 1 | 8 | - | - | - |
| 5 | 6 | 11 | 14 | 15 | 29 | - | - | - |
| 7 | 5 | 12 | 36 | 26 | 62 | - | - | - |
| 11 | 1 | 12 | 31 | 17 | 48 | - | - | - |
| 5 | 5 | 10 | 34 | 25 | 59 | - | - | - |
| 7 | 2 | 9 | 26 | 18 | 44 | 5 | 1 | 6 |
| 7 | 2 | 9 | 10 | 10 | 20 | 14 | 12 | 26 |
| 8 | 6 | 14 | 3 | 5 | 8 | 61 | 41 | 102 |
| - | 2 | 2 | 1 | 1 | 2 | 54 | 23 | 77 |
| 51 | 31 | 82 | 162 | 118 | 280 | 134 | 77 | 211 |
| 1 | 1 | 2 | 5 | 3 | 8 | 4 | 2 | 6 |

| Types of PWDs | | | | | |
|---|---|---|---|---|---|
| LV | | | TB | | |
| M | F | Total | M | F | Total |
| - | - | - | - | - | - |
| 4 | - | 4 | - | - | - |
| 14 | 3 | 17 | 4 | 5 | 9 |
| 5 | 3 | 8 | 4 | 4 | 8 |
| 5 | 7 | 12 | 5 | 11 | 16 |
| 14 | 6 | 20 | 10 | 12 | 22 |
| 10 | 5 | 15 | 12 | 12 | 24 |
| 17 | 19 | 36 | 30 | 23 | 53 |
| 5 | 1 | 6 | 14 | 6 | 20 |
| **74** | **44** | **118** | **79** | **73** | **152** |
| **2** | **1** | **3** | **2** | **2** | **4** |

| Types of PWDs | | | | | | | | |
|---|---|---|---|---|---|---|---|---|
| MI | | | MD | | | S&HI | | |
| M | F | Total | M | F | Total | M | F | Total |
| - | - | - | 5 | 1 | 6 | - | 2 | 2 |
| 1 | - | 1 | 14 | 17 | 31 | 15 | 16 | 31 |
| 1 | 1 | 2 | 40 | 22 | 62 | 35 | 43 | 78 |
| 1 | 1 | 2 | 10 | 11 | 21 | 19 | 22 | 41 |
| 5 | 9 | 14 | 18 | 12 | 30 | 40 | 48 | 88 |
| 18 | 18 | 36 | 14 | 23 | 37 | 49 | 35 | 84 |
| 11 | 7 | 18 | 10 | 7 | 17 | 55 | 51 | 106 |
| 12 | 6 | 18 | 20 | 7 | 27 | 50 | 37 | 87 |
| 3 | 3 | 6 | 12 | 6 | 18 | 27 | 18 | 45 |
| **52** | **45** | **97** | **143** | **106** | **249** | **290** | **272** | **562** |
| **1** | **1** | **3** | **4** | **3** | **7** | **8** | **8** | **16** |

Out of the total of 3539 Persons With Disabilities (PWDs), a majority of 2115 are men, that form 60 percent and the rest 40 percent are women. Of the total PWDs 14 percent are children (0-14 yrs) and 40 percent are in the reproductive age group (15-35 yrs). Among the nine categories of disabilities, a majority of 51 percent (31% male and 20% female) are orthopedic, followed by 16 percent (8% of each sex) are speech and Hearing impaired. We also notice the proof for the claim that leprosy has been completely eradicated upto 25 years of age. There is not a single person with leprosy and only six persons are found leprosy cured in the 25-35 age group. As from our interviews most of them complained that the caretakers were least interested in getting them married off in the reproductive age group. We were forced to give a numerical physical analysis of the marital status as we found it to be one of the major psychological problems suffered by the PWDs.

### Marital status of PWDs in the reproductive age group

| Sl. No. | Age Group | Total PWDs | | | Marital Status | | |
|---|---|---|---|---|---|---|---|
| | | | | | Married | | |
| | | M | F | Total | M | F | Total |
| 1 | 19 to 25 | 317 | 256 | 573 | 45 | 62 | 107 |
| 2 | 26 to 35 | 371 | 247 | 618 | 225 | 114 | 339 |
| Total | | 688 | 503 | 1,191 | 270 | 176 | 446 |
| Percentage | | 58 | 42 | 100 | 23 | 15 | 37 |

| Marital Status | | | | | | | | | | | |
|---|---|---|---|---|---|---|---|---|---|---|---|
| Unmarried | | | Widowed | | | Divorced | | | Separated | | |
| M | F | Total | M | F | Total | M | F | Total | M | F | Total |
| 271 | 190 | 461 | - | 1 | 1 | - | 2 | 2 | 1 | 1 | 2 |
| 138 | 116 | 254 | 5 | 11 | 16 | 2 | 4 | 6 | 1 | 2 | 3 |
| 409 | 306 | 715 | 5 | 12 | 17 | 2 | 6 | 8 | 2 | 3 | 5 |
| 34 | 26 | 60 | 0 | 1 | 1 | 0 | 1 | 1 | 0 | 0 | 0 |

Of the total of 1191 PWDs in the marriageable age group (19-35 yrs), a majority of 715 PWDs which forms 60 percent are not married. This is due to the fact that the public in general and the parents/caretakers in particular are not worried about the age specific need of the PWDs. Also fear to give or take in marriage alliance for they fear such disability may recur in their family if not in this generation in future generations. The level of social acceptance of the PWDs is reflected by this fact. There are families where the members of the family live on the earning of the PWD, whose socio psychological and physical need - marriage is not given much importance by those family members who live on his/her income.

**Marital status of PWDs in the reproductive age group**

| Sl. No. | Age Group | Total PWDs | | | Marital Status | | |
|---|---|---|---|---|---|---|---|
| | | | | | Married | | |
| | | M | F | Total | M | F | Total |
| 1 | 36 to 45 | 331 | 206 | 537 | 268 | 128 | 396 |
| 2 | 46 to 60 | 478 | 282 | 760 | 391 | 171 | 562 |
| 3 | 61 & above | 203 | 89 | 292 | 159 | 46 | 205 |
| Total | | 1,012 | 577 | 1,589 | 818 | 345 | 1,163 |
| Percentage | | 64 | 36 | 100 | 51 | 22 | 73 |

| Marital Status | | | | | | | | | | | |
|---|---|---|---|---|---|---|---|---|---|---|---|
| Unmarried | | | Widowed | | | Divorced | | | Separated | | |
| M | F | Total | M | F | Total | M | F | Total | M | F | Total |
| 52 | 39 | 91 | 5 | 31 | 36 | 5 | 5 | 10 | 1 | 3 | 4 |
| 61 | 31 | 92 | 19 | 72 | 91 | 4 | 3 | 7 | 3 | 5 | 8 |
| 23 | 7 | 30 | 20 | 34 | 54 | 1 | - | 1 | - | 2 | 2 |
| 136 | 77 | 213 | 44 | 137 | 181 | 10 | 8 | 18 | 4 | 10 | 14 |
| 9 | 5 | 13 | 3 | 9 | 11 | 1 | 1 | 1 | 0.3 | 1 | 1 |

From this table we understand that there are 1589 PWDs of which a majority 1012 forming 64 percent are men and the rest women. Of these population, a majority of 1163 constituting 73 percent are found to be married, 11 percent are widowers, and negligible 1 percent each are divorced and separated. Of those married, majority are men and hardly one-third is women.

While comparing the above two tables, one can easily come to the conclusion that the younger generation are neither willing nor encouraged / motivated / assisted to get married compared to the older generation. The causes for this trend need to be probed further and some of the observations are given in the last chapter of this book.

Next we proceed on to give the table form of the educational status of the PWDs in the following table.

## The age-group wise stratification of PWDs

| Sl. No. | Age Group | Total PWDs | | |
|---|---|---|---|---|
| | | M | F | Total |
| 1 | 4 to 6 | 76 | 66 | 142 |
| 2 | 7 to 14 | 185 | 156 | 341 |
| 3 | 15 to 18 | 136 | 103 | 239 |
| 4 | 19 to 25 | 315 | 256 | 571 |
| 5 | 26 to 35 | 372 | 248 | 620 |
| 6 | 36 to 45 | 332 | 210 | 542 |
| 7 | 46 to 60 | 479 | 282 | 761 |
| 8 | Above 60 | 203 | 89 | 292 |
| Total | | 2,098 | 1,410 | 3,508 |
| Percentage | | 60 | 40 | 100 |

| Educational Status | | | | | |
|---|---|---|---|---|---|
| None | | | Studying | | |
| M | F | Total | M | F | Total |
| - | - | - | 13 | 13 | 26 |
| 25 | 22 | 47 | 112 | 97 | 209 |
| 39 | 45 | 84 | 47 | 23 | 70 |
| 125 | 119 | 244 | 27 | 11 | 38 |
| 177 | 168 | 345 | 2 | 1 | 3 |
| 209 | 175 | 384 | - | - | - |
| 328 | 246 | 574 | - | - | - |
| 132 | 81 | 213 | - | - | - |
| 1,035 | 856 | 1,891 | 201 | 145 | 346 |
| 30 | 24 | 54 | 5.7 | 4.1 | 9.9 |

| Educational Status | | | | | | | | |
|---|---|---|---|---|---|---|---|---|
| Dropout | | | Not yet enrolled | | | Literate | | |
| M | F | Total | M | F | Total | M | F | Total |
| - | - | - | 63 | 53 | 116 | - | - | - |
| 13 | 7 | 20 | 35 | 30 | 65 | - | - | - |
| 50 | 35 | 85 | - | - | - | - | - | - |
| - | - | - | - | - | - | 163 | 126 | 289 |
| - | - | - | - | - | - | 193 | 79 | 272 |
| - | - | - | - | - | - | 123 | 35 | 158 |
| - | - | - | - | - | - | 151 | 36 | 187 |
| - | - | - | - | - | - | 71 | 8 | 79 |
| 63 | 42 | 105 | 98 | 83 | 181 | 701 | 284 | 985 |
| 1.8 | 1.2 | 3 | 2.8 | 2.4 | 5 | 20 | 8 | 28 |

Table reveals the fact that, of the total of 3508 PWDs, which exclude the children below 3yrs, a majority of 59 percent have not entered the school premises. This once again shows that they are not interested to attend school with others.

### Disability Category-wise Educational Status-School drop out

| Sl. No. | Age Group | Total PWDs | | | Types of PWDs | | |
|---|---|---|---|---|---|---|---|
| | | | | | Ortho | | |
| | | M | F | Total | M | F | Total |
| 1 | 7 to 14 | 13 | 7 | 20 | 3 | 1 | 4 |
| 2 | 15 to 18 | 50 | 35 | 85 | 20 | 21 | 41 |
| Total | | 63 | 42 | 105 | 23 | 22 | 45 |
| Percentage | | 11 | 7 | 18 | 4 | 4 | 8 |

| Types of PWDs | | | | | | | | |
|---|---|---|---|---|---|---|---|---|
| CP | | | MR | | | LC | | |
| M | F | Total | M | F | Total | M | F | Total |
| 1 | - | 1 | 3 | 1 | 4 | - | - | - |
| 5 | - | 5 | 9 | 2 | 11 | - | - | - |
| 6 | - | 6 | 12 | 3 | 15 | - | - | - |
| 1 | - | 1 | 2 | 1 | 3 | - | - | - |

| Types of PWDs | | | | | |
|---|---|---|---|---|---|
| LV | | | TB | | |
| M | F | Total | M | F | Total |
| 1 | - | 1 | - | - | - |
| 2 | 1 | 3 | 1 | 1 | 2 |
| 3 | 1 | 4 | 1 | 1 | 2 |
| 1 | 0 | 1 | 0 | 0 | 0 |

| Types of PWDs | | | | | | | | |
|---|---|---|---|---|---|---|---|---|
| MI | | | MD | | | S&HI | | |
| M | F | Total | M | F | Total | M | F | Total |
| - | - | - | 3 | 1 | 4 | 2 | 4 | 6 |
| 1 | 1 | 2 | 3 | 2 | 5 | 9 | 7 | 16 |
| 1 | 1 | 2 | 6 | 3 | 9 | 11 | 11 | 22 |
| 0 | 0 | 0 | 1 | 1 | 2 | 2 | 2 | 4 |

In the age group of 7-18 yrs, there are 220 children, of which 105 children dropped out, which constitute 18 percent. Of these dropouts a majority of 48 children (8 percent) are from ortho, 22 children (4 percent) are from speech and hearing impaired and 15 children (3 percent) are from mentally retarded.

## Disability Category-wise educational status - not yet enrolled

| Sl. No. | Age Group | Total PWDs | | | Types of PWDs | | |
|---|---|---|---|---|---|---|---|
| | | | | | Ortho | | |
| | | M | F | Total | M | F | Total |
| 1 | 4 to 6 | 63 | 53 | 116 | 17 | 11 | 28 |
| 2 | 7 to 14 | 35 | 30 | 65 | 2 | 2 | 4 |
| Total | | 98 | 83 | 181 | 19 | 13 | 32 |
| Percentage | | 20 | 17 | 37 | 4 | 3 | 7 |

| Types of PWDs | | | | | | | | |
|---|---|---|---|---|---|---|---|---|
| CP | | | MR | | | LC | | |
| M | F | Total | M | F | Total | M | F | Total |
| 5 | 6 | 11 | 11 | 14 | 25 | - | - | - |
| - | 2 | 2 | 9 | 9 | 18 | - | - | - |
| 5 | 8 | 13 | 20 | 23 | 43 | - | - | - |
| 1 | 2 | 3 | 4 | 5 | 9 | - | - | - |

| Types of PWDs | | | | | |
|---|---|---|---|---|---|
| LV | | | TB | | |
| M | F | Total | M | F | Total |
| 3 | - | 3 | - | - | - |
| - | - | - | 2 | - | 2 |
| 3 | - | 3 | 2 | - | 2 |
| 1 | - | 1 | 0 | - | 0 |

| Types of PWDs | | | | | | | | |
|---|---|---|---|---|---|---|---|---|
| MI | | | MD | | | S&HI | | |
| M | F | Total | M | F | Total | M | F | Total |
| 1 | - | 1 | 14 | 14 | 28 | 12 | 8 | 20 |
| - | - | - | 19 | 14 | 33 | 3 | 3 | 6 |
| 1 | - | 1 | 33 | 28 | 61 | 15 | 11 | 26 |
| 0 | - | 0 | 7 | 6 | 13 | 3 | 2 | 5 |

This table gives the fact that 181 children out of 483 that constitute 37 percent are yet to be enrolled in the school. Of these majority of 61 children (13 percent) are from multiple disabilities (MD), 43 children (9 percent) are from mentally retarded (MR), 32 children (7 percent) are from orthopaedically handicapped and 26 children (5 percent) from the speech and hearing impaired.

**Disability Category-wise educational status-uneducated**

| Sl. No. | Age Group | Total PWDs | | | Types of PWDs Ortho | | |
|---|---|---|---|---|---|---|---|
| | | M | F | Total | M | F | Total |
| 1 | 7 to 14 | 25 | 22 | 47 | 1 | 1 | 2 |
| 2 | 15 to 18 | 39 | 45 | 84 | 5 | 9 | 14 |
| 3 | 19 to 25 | 125 | 119 | 244 | 49 | 41 | 90 |
| 4 | 26 to 35 | 177 | 168 | 345 | 90 | 73 | 163 |
| 5 | 36 to 45 | 209 | 175 | 384 | 115 | 87 | 202 |
| 6 | 46 to 60 | 328 | 246 | 574 | 183 | 117 | 300 |
| 7 | 61 & above | 132 | 81 | 213 | 54 | 26 | 80 |
| **Total** | | **1,035** | **856** | **1,891** | **497** | **354** | **851** |
| **Percentage** | | **31** | **25** | **56** | **15** | **11** | **25** |

| Types of PWDs | | | | | | | | |
|---|---|---|---|---|---|---|---|---|
| CP | | | MR | | | LC | | |
| M | F | Total | M | F | Total | M | F | Total |
| 2 | 2 | 4 | 7 | 5 | 12 | - | - | - |
| 2 | 1 | 3 | 18 | 14 | 32 | - | - | - |
| 5 | 3 | 8 | 26 | 15 | 41 | - | - | - |
| 2 | 2 | 4 | 17 | 14 | 31 | 4 | 1 | 5 |
| 6 | 1 | 7 | 8 | 8 | 16 | 12 | 11 | 23 |
| 3 | 5 | 8 | 3 | 5 | 8 | 48 | 36 | 84 |
| - | 2 | 2 | 1 | 1 | 2 | 42 | 21 | 63 |
| **20** | **16** | **36** | **80** | **62** | **142** | **106** | **69** | **175** |
| **1** | **0** | **1** | **2** | **2** | **4** | **3** | **2** | **5** |

| Types of PWDs | | | | | |
|---|---|---|---|---|---|
| **LV** | | | **TB** | | |
| **M** | **F** | **Total** | **M** | **F** | **Total** |
| - | - | - | - | 2 | 2 |
| 2 | 1 | 3 | - | 1 | 1 |
| 1 | 2 | 3 | 2 | 6 | 8 |
| 4 | 5 | 9 | 3 | 10 | 13 |
| 6 | 2 | 8 | 11 | 11 | 22 |
| 10 | 18 | 28 | 27 | 21 | 48 |
| 3 | 1 | 4 | 10 | 6 | 16 |
| **26** | **29** | **55** | **53** | **57** | **110** |
| **1** | **1** | **2** | **2** | **2** | **3** |

| Types of PWDs | | | | | | | | |
|---|---|---|---|---|---|---|---|---|
| **MI** | | | **MD** | | | **S&HI** | | |
| **M** | **F** | **Total** | **M** | **F** | **Total** | **M** | **F** | **Total** |
| - | - | - | 11 | 6 | 17 | 4 | 6 | 10 |
| - | - | - | 6 | 9 | 15 | 6 | 10 | 16 |
| - | 5 | 5 | 12 | 9 | 21 | 30 | 38 | 68 |
| 9 | 12 | 21 | 10 | 20 | 30 | 38 | 31 | 69 |
| 5 | 6 | 11 | 6 | 5 | 11 | 40 | 44 | 84 |
| 7 | 4 | 11 | 15 | 6 | 21 | 32 | 34 | 66 |
| 3 | 3 | 6 | 7 | 5 | 12 | 12 | 16 | 28 |
| **24** | **30** | **54** | **67** | **60** | **127** | **162** | **179** | **341** |
| **1** | **1** | **2** | **2** | **2** | **4** | **5** | **5** | **10** |

Of the total of 3316 PWDs in this age group, it is observed that 1891 which constitute 56 percent are not educated. Note also that a majority are from orthopaedically handicapped (25 percent) and from speech and hearing impaired (10 percent).

Next we proceed on to give the table depicting the employment.

## Employment Status of the PWDs

| Sl. No. | Age Group | Total PWDs | | | Current Status Employed | | |
|---|---|---|---|---|---|---|---|
| | | **M** | **F** | **Total** | **M** | **F** | **Total** |
| 1 | 15 to 18 | 136 | 103 | 239 | 19 | 16 | 35 |
| 2 | 19 to 25 | 315 | 256 | 571 | 102 | 48 | 150 |
| 3 | 26 to 35 | 372 | 248 | 620 | 162 | 65 | 227 |
| 4 | 36 to 45 | 332 | 210 | 542 | 126 | 65 | 191 |
| 5 | 46 to 60 | 479 | 282 | 761 | 145 | 73 | 218 |
| 6 | 61 & above | 203 | 89 | 292 | 49 | 22 | 71 |
| **Total** | | 1,837 | 1,188 | 3,025 | 603 | 289 | 892 |
| **Percentage** | | 61 | 39 | 100 | 20 | 10 | 29 |

| Current status | | | | | |
|---|---|---|---|---|---|
| Unemployed | | | Unable to work | | |
| **M** | **F** | **Total** | **M** | **F** | **Total** |
| 37 | 26 | 63 | 80 | 61 | 141 |
| 96 | 75 | 171 | 117 | 133 | 250 |
| 65 | 59 | 124 | 145 | 124 | 269 |
| 59 | 41 | 100 | 147 | 104 | 251 |
| 101 | 65 | 166 | 233 | 144 | 377 |
| 42 | 22 | 64 | 112 | 45 | 157 |
| **400** | **288** | **688** | **834** | **611** | **1,445** |
| **13** | **10** | **23** | **28** | **20** | **48** |

Table gives age wise distribution of the employed, unemployed, unable to work and inapplicable categories. The children numbering 514 out of the total of 3539 PWDs constitute 14 percent of the total PWDs, 25 percent of the total PWDs are employed and 19 percent are unemployed. A significant 41 percent, which constitute two-thirds of the entire PWDs, are not in a position to work due to their disability. Among those employed, men are twice as many as women.

From our study we could get only the following as the reasons given by the Parents / Caretakers for not assisting the PWD's in their daily activities given in the table 2.6 for quick analysis. Next we began to discuss with them (PWDS) about their awareness on their rights in the society. Majority of them were unaware of their rights. This is briefly described in the table with the interpretations in the following.

Reasons as expressed by parents/ caretakers for not helping the PWDs in daily activities.

**Reasons as expressed by parents/ caretakers …**

| Time Spent | No. of Respondent | % |
|---|---|---|
| Busy with Other Duties | 1710 | 72.9 |
| Not Inclined | 229 | 9.8 |
| Not Keeping Good Health | 87 | 3.7 |
| Ashamed/Stigma | 17 | 0.7 |
| Others | 303 | 12.9 |
| **Total** | **2346** | **100.0** |

Out of 2346 respondents 72.9 percent of the parents said that they were busy with other work and hence could not spend their time with their disabled children in doing their routine daily activities; 9.8 percent parents said that they were not inclined to do so; 3.7 percent parents said that they could not help the children due to their ill health and yet another 0.7 percent parents said that they were ashamed in doing so. So, we notice that around 87 percent of the parents do not help their disabled children in their daily activities due to several reasons.

**Awareness of the PWDs on their rights**

| Awareness on their rights | Yes | 44 | 1.2% |
|---|---|---|---|
| | No | 3,492 | 98.8% |
| | **Total** | **3,536** | **100%** |

The question provoked the response of the PWDs to know how far they are knowledgeable on their rights. Out of the 3536 PWDs responded, a majority (98.8 percent) of PWDs expressed that they do not have knowledge on their rights. This gives the grim picture of PWDs on their 'awareness level' on their own rights. They live in a 'culture of silence and ignorance'. They were not aware of the government policies which are given in table.

### Awareness of the PWDs on Government Policies related to their rehabilitation

| Awareness on policy | Yes | 36 | 1.0% |
|---|---|---|---|
| measures relating to | No | 3,493 | 99.0% |
| their special care | **Total** | **3,529** | **100%** |

This table again views the grim sorry state of PWDs on their awareness on government policies. Out of 3529 PWDs who responded, all most all have expressed that they do not have knowledge on the policy perspectives of the government that ensures their basic rights.

Questions were asked about the future plan of the PWDs. Totally 1404 PWDs responded. A majority (1134 - 80 percent) expressed that they plan to take up 'self employment, 72 respondents (5 percent) preferred 'Govt. Job', 72 respondents (5 percent) preferred to go for 'higher education' and the remaining (10 percent) plan to become teachers, agriculturists, vocational training, doctors, collectors, Panchayat Board President, etc. This is given in the following table.

### Future Plan expressed by PWDs

| | | | |
|---|---|---|---|
| 1 | Self Employment | 1134 | 80.8% |
| 2 | Government Job | 72 | 5.1% |
| 3 | Higher Education | 72 | 5.1% |
| 4 | Teacher | 42 | 3.0% |
| 5 | Better Work | 31 | 2.2% |
| 6 | Agricultural | 17 | 1.2% |
| 7 | Vocational Training | 13 | 0.9% |
| 8 | Doctor | 7 | 0.5% |
| 9 | Collector | 6 | 0.4% |
| 10 | Tailor | 5 | 0.4% |
| 11 | Panchayat President | 2 | 0.1% |
| 12 | Bank Manager | 1 | 0.1% |
| 13 | Engineer | 1 | 0.1% |
| 14 | Professor | 1 | 0.1% |
| | **Total** | **1404** | **100.0%** |

From the explorative study conducted in a population of 2,15,811 living in 93 Panchayats of Melmalayanur and Kurinjipadi Blocks of Tamil Nadu, India 3539 persons with disabilities were identified in the 9 categories of disabilities, following are the main findings of the study.

1. Population: Among the 9 types of disabilities, majority (51 percent) of them were found to be orthopedic, followed by the speech and hearing impaired who form 16 percent and the mentally retarded 8 percent, the multiple disabilities form 7 percent, and the remaining categories i.e., 17 percent form other 5 categories of disabilities.

2. Marriage: Among the 'reproductive age group' of PWDs (18-35 years), it was found that 60 percent of them were unmarried, which reflects the low level of social acceptance and high level of social stigma coupled with lack of self esteem.

3. Caste: People living below poverty line in Tamil Nadu mostly come under the MBC (Most Backward Community) and SC (Schedule Caste). The SCs are called as the 'untouchables' or 'dalits'. The survey findings prove the fact that majority of the PWDs are found to be MBC (55 percent) and SC (25 percent) category. Only the remaining one-fifth of the PWDs are from the other communities.

4. Education: Among the school going children (7 to 18 years) age group (i.e., 580 children with disabilities), 18 percent are found to be school dropouts. Among the children in the 4-14 age group, it was found that 37 percent are not enrolled in school. This may be due to various reasons such as ignorance of parents about the importance of education for their children with disabilities, mobility problem of the children with disabilities, economic problem of their families and so on. The overall literacy level of the persons with disabilities is 46 percent.

5. Family: Most of the PWDs (94.4 percent) are found to be living with their families. This is a positive trend, which need to be encouraged and strengthened carefully without affecting the interpersonal relationship.

6. Psychological problems faced by the PWDs: The feelings expressed by the PWDs is quite disheartening; 66 percent of the PWDs expressed that they are 'depressed', 8 percent 'irritated', 12 percent feel that disability is their 'fate' and 8 percent of PWDs feel 'indifferent'. All these connote the kind of psychological agony and distress being experienced by the PWDs. This issue needs to be professionally addressed.

Thirty five percent of the PWDs feel that they are either treated as 'indifferent' or 'abused'. Family support and friends' support were considered as means of overcoming suicidal tendencies. Thirty two percent of PWDs are found to be having 'low self image' of themselves. This needs to be professionally handled.

A majority (57.7 percent) of the PWDs suffer from 'fear of illness', 13.4 percent of them have the fear of their 'parents death' who are their primary caregivers, 10.7 percent of them face 'rejection' by their family members and 12.9 percent of them have the fear of their 'old age'.

7. Community participation: About less than 30 percent of PWDs have expressed that they were not being introduced to the visitors in the families, and were not allowed to socialise with the neighbours. Similarly about 25-30 percent of them expressed their feelings of shame for their disabilities, hence, reluctant to interact with others in the community.

Forty four percent of the Parents / Caregivers expressed that they have not planned the future of the PWDs who were dependent on them.

With regard to the future plan by the parents for PWDs, it was found that 33.6 percent of them suggested 'custodial care', 24.6 percent of the parents of PWDs expressed the need for 'home-based care', 26.6 percent sought assistant for venturing into self employment of PWDs, 8.4 percent sought assistance for taking insurance policy, 3.7 percent suggested institutional care and so on.

8. Awareness on Rights: It is studied that 98.8 and 99 percent of the PWDs are unaware of their 'rights' and 'government policies' respectively.

**Chapter Three**

# THE NEW FUZZY COGNITIVE RELATIONAL MAPS (FCRMS) BIMODEL

In this chapter we for the first time introduce a new model called the Fuzzy Cognitive Relational Maps bimodel (FCRMs bimodel). The main feature of this model is that it can give a pair of hidden patterns for a different pair of attributes. By this one can understand or study or get to know the underlying effect of a pair of attributes over this dynamical bisystem.

Before we go for the definition and description of this model we just give a need for such a model. When we want to consider a pair of attributes one related with the FCM dynamical system and the other with a FRM dynamical system and wish to see the stage by stage effect of an attribute on the system this model will be well suited. For such comparisons of the effects of the attributes as well as the hidden pattern this model would be extremely useful.

Further as all our problems under investigation cannot be given solution in terms of statistical data or be exhibited by numerical values we feel only this model is capable of analyzing the problem at hand and giving the effects of the attributes on the dynamical bisystem.

With this view we now proceed onto define the new bimodel.

## Definition and description of new fuzzy cognitive relational maps bimodel (FCRM bimodel)

Suppose we have a problem at hand, and the data related with it happens to be an unsupervised one, obtained from survey as well as interviews and discussions. It may so happen that the attributes in one case or under one expert the attributes may not be a disjoint one and under another expert the attributes can be made disjoint. In such a case this new bimodel will be best suited. From the survey data, interviews and discussions and analysis, this new bimodel would be appropriate to study the problem.

We have just given references in chapter one about the notion of Fuzzy Cognitive Maps and Fuzzy Relational Maps and their functioning.

A Fuzzy Cognitive Relational Maps (FCRM) is a directed special bigraph with concepts like policies, events, etc as nodes and causalities as edges. It represents causal relationship between concepts.

We first illustrate this by a simple example.

***Example 3.1:*** Suppose one wishes to study the problem of female infanticide (the practice of killing female children at birth or shortly thereafter) which is prevalent in India from early vedic times as women were (and still are) considered as a liability. They are considered to be appendages. As long as women are treated as a property / object the practice of female infanticide will continue in India.

The attributes given by one expert are as follows

$F_1$ - Social status (a stigma if one has only female children)

$F_2$ - Economic condition (burden / loss if one has more female children)

$F_3$ - Torture by in laws for having given birth to only female children.

$F_4$ - Insecurity to have only girl children (As per Hindu laws they leave their parents and live in a different home after marriage)

$F_5$ - A problem if one has many female children (as dowry in time of marriage, and functions thereafter).

These are the attributes used by an expert. These cannot be further subdivided and they are in the opinion of the expert dependent attributes so he chooses to work with the FCMs [19-21, 44].

The second expert however feels that the 5 attributes $F_1$ ,…, $F_5$ given by the first expert to be independent and feels that, the economy of the parents or the caretakers of these children are important as they influence the female infanticide. So he/her gives another set of attributes say $E_1$, …, $E_7$.

One set of attributes, which can be taken as $F_1$, $F_2$, $F_3$, $F_4$ and $F_5$. These are taken as the attributes related with the domain space. The set of attributes given by second expert for range space is as follows:

$E_1$ - very rich
$E_2$ - rich
$E_3$ - upper middle class
$E_4$ - middle class

$E_5$ - lower middle class
$E_6$ - poor
$E_7$ - very poor.

This expert feels the attributes $E_1$ ,..., $E_7$ are independent and hence wants to study the female infanticide problem using FRM [44,47,48,57].

Thus instead of using two models, we combine them and will use the FCRM bimodel.

The directed special bigraph associated with this problem is as follows.

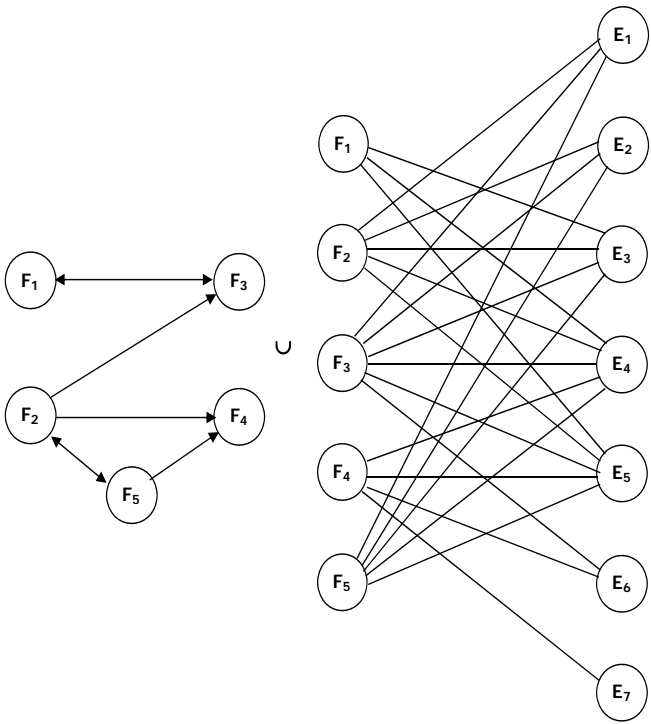

The related special mixed bimatrix M associated with this special bigraph is as follows.

$$M = M_1 \cup M_2 =$$

$$
\begin{array}{c}
\quad F_1\ F_2\ F_3\ F_4\ F_5 \\
\begin{array}{c}
F_1 \\ F_2 \\ F_3 \\ F_4 \\ F_5
\end{array}
\begin{bmatrix}
0 & 0 & 1 & 0 & 0 \\
0 & 0 & 1 & 1 & 1 \\
1 & 0 & 0 & 0 & 0 \\
0 & 0 & 0 & 0 & 1 \\
0 & 1 & 0 & 0 & 0
\end{bmatrix}
\end{array}
\cup
\begin{array}{c}
\quad E_1\ E_2\ E_3\ E_4\ E_5\ E_6\ E_7 \\
\begin{array}{c}
F_1 \\ F_2 \\ F_3 \\ F_4 \\ F_5
\end{array}
\begin{bmatrix}
0 & 0 & 1 & 1 & 1 & 0 & 0 \\
1 & 1 & 1 & 1 & 1 & 0 & 0 \\
1 & 1 & 1 & 1 & 1 & 1 & 0 \\
0 & 0 & 0 & 1 & 1 & 1 & 1 \\
1 & 1 & 1 & 1 & 1 & 0 & 0
\end{bmatrix}
\end{array}
$$

is called the Connection Relational bimatrix or in short CR-bimatrix. We have used the properties of FCM and FRM given [44,47,48,57]. Since the nodes in this model occur as pairs these nodes of FCRM will be known as binodes. If the binodes of the FCRM bimodel are fuzzy bisets then they are called as fuzzy binodes. FCRM with edge weights or causalities from the set $\{-1, 0, 1\}$ are called simple FCRMs. FCRM with edge weights or causalities from the set $\{0, 1\}$ will be called as positive FCRMs. We can have mixed simple-positive or positive simple FCRMs. In such case one component will be simple other positive and vice versa. Consider the binodes / biconcepts $\{(C_1, \ldots, C_n)\ (R_1, \ldots, R_m)\}$ of the FCRM. Suppose the special directed bigraph is drawn using biedge weights $(e_{ij}, \delta_{ij})$; $e_{ij}, \delta_{ij} \in \{0, 1, -1\}$. The bimatrix $E = E_1 \cup E_2$ will be a special mixed bimatrix where $e_{ij}$ is the weight of the directed edge $C_i C_j$ and $\delta_{ij}$ is the weight of the directed edge $C_i R_j$. The components of the bigraph are the directed graph of the FCM and the directed graph of the FRM.

The special mixed bimatrix will be called as the special mixed connection relational bimatrix of FCRM.

Here we assume $C_1, \ldots, C_n$ to be nodes of the FCM of the FCRM and $\{(C_1, \ldots, C_n), (R_1, \ldots, R_m)\}$ of the FRM of the FCRM respectively FCMs of the FCRM will be known as the first component and FRMs of the FCRM bimodel will be known as the second component of the bimodel, FRMs.

$$A = A_1 \cup A_2 = \left(a_1^1,...,a_n^1\right) \cup \left(a_1^2,...,a_m^2\right) \text{ where } a_i \in \{0,1\};$$

A is called the instantaneous state bivector and it denotes the on-off position of the binode at an instant.

$$a_j^i = 0 \text{ if } a_j^i \text{ off and}$$

$$a_j^i = 1 \text{ if } a_j^i \text{ is on;}$$

for i = 1, 2 and $1 \leq j \leq m, n$.

As in case of FCMs and FRMs we in case of FCRMs also assume that the FCRM is a bicycle. We call FCRM which is bicyclic and in which the causal birelations flow in a revolutionary way to be a dynamical bisystem. The biequilibrium state of the dynamical bisystem is called the hidden bipattern.

If the equilibrium bistate of the dynamical bisystem is a unique bistate bivector then it is called a fixed bipoint. If the FCRM settles down with a state bivector repeating in the form then this equilibrium is called a limit bicycle.

As in case of FCMs and FRMs, FCRMs can be combined together to produce the joint effect of all the FCRMs.

If $E_i = E_i^1 \cup E_i^2$ (i = 1, ..., n where $E_i$ is the opinion of the $i^{th}$ expert). We have taken n such experts are the special connection relational mixed bimatrices of the FCRMs with binodes $\{(C_1, ..., C_n), (R_1, ..., R_m)\}$ then the combined FCRMs is got by adding all the mixed bimatrices

$$\begin{aligned} E_1 + ... + E_n \quad &= \left(E_1^1 \cup E_1^2\right) + ... + \left(E_n^1 \cup E_n^2\right) \\ &= \left(E_1^1 + ... + E_n^1\right) \cup \left(E_1^2 + ... + E_n^2\right). \end{aligned}$$

Clearly for the addition to be compatible we assume each of the matrices $E_i^1$ are square matrices of same order $1 \leq i \leq n$ and $E_j^2$ are rectangular matrices of same order for $1 \leq j \leq n$.

We now define the concepts related with the new FCRM bimodel.

**DEFINITION 3.1:** *In a new FCRM we call the pair of associated nodes as binodes. If the order of the bimatrix associated with a FCRM is a n $\times$ n square matrix and a p $\times$ m matrix then the binodes are bivectors of length (n, p) or length (n, m).*

**Example 3.2:** If the new FCRM matrix is given by say

$$M = M_1 \cup M_2 = \begin{pmatrix} 0 & 0 & -1 & 1 \\ 1 & 0 & 1 & -1 \\ -1 & 0 & 0 & 1 \\ 1 & 1 & 0 & 0 \end{pmatrix} \cup \begin{pmatrix} 1 & 0 & -1 \\ 0 & 1 & 0 \\ -1 & 0 & 1 \\ 0 & 1 & 0 \\ 1 & -1 & 0 \\ 1 & 1 & -1 \end{pmatrix}$$

then the binodes can be $(a_1\ a_2\ a_3\ a_4) \cup (b_1\ b_2\ b_3)$ or $(a_1\ a_2\ a_3\ a_4) \cup (c_1, c_2, c_3, c_4, c_5, c_6)$ as the bimatrix associated with it is a $(4 \times 4, 6 \times 3)$ bimatrix hence the binodes can be of length either $(4, 6)$ or $(4, 3)$.

**DEFINITION 3.2:** *Consider the binodes / biconcepts $\{ C_1^1,...,C_n^1 \}$ of the FCM and $\{D_1, ..., D_p\}$ and $\{R_1, ..., R_m\}$ of the FRM of the new FCRM bimodel.*

*Suppose the directed graph is drawn using the edge biweight $e_{ij}^t = \{0, 1, -1\}$; $1 \leq t \leq 2$. The bimatrix $E = E_1 \cup E_2$ is defined by $E = (e_{ij}^1) \cup (e_{ij}^2)$ where $e_{ij}^1$ is the weight of the directed edge $C_i\ C_j$ and $e_{ks}^2$ is the directed edge of $D_k\ R_s$. $E = E_1 \cup E_2$ is called the adjacency bimatrix of the new FCRM* *bimodel, also known as the connecting relational bimatrix of the new FCRM bomodel.*

It is interesting to note that all bimatrices associated with an FCRM are always a square matrix with diagonal entries zero and a rectangular $p \times m$ matrix.

**DEFINITION 3.3:** *The new FCRMs with edge biweight {–1, 0, 1} are called simple FCRMs.*

In this book we work only with simple new FCRMs bimodel.

**DEFINITION 3.4:** *Let $\{C_1, …, C_n\} \cup \{(D_1, …, D_p), (R_1, …, R_m)\}$ be the binodes of an FCRM. $A = A_1 \cup A_2 = (a_1, …, a_n) \cup (b_1, …, b_p)$ (or $(c_1, …, c_m)$) where $a_i, b_j, c_t \in \{0, 1\}$; $1 \le i \le n$, $1 \le j \le p$ and $1 \le t \le m$. A is called the instantaneous state bivector and it denotes the on-off position of the node at an instant.*

$$a_i = 0 \text{ if } a_i \text{ is off and}$$
$$a_i = 1 \text{ if } a_i \text{ is on for } i = 1, 2, …, n.$$

$$b_i = 0 \text{ if } b_i \text{ is off and}$$
$$b_i = 1 \text{ if } b_i \text{ is on for } i = 1, 2, …, p.$$

$$c_i = 0 \text{ if } c_i \text{ is off and}$$
$$c_i = 1 \text{ if } c_i \text{ is on for } i = 1, 2, …, m.$$

**DEFINITION 3.5:** *Let $\{(C_1, …, C_n)\} \cup \{(D_1, …, D_p), (R_1, …, R_m)\}$ be the binodes of the new FCRM. Let $\overline{C_i C_j} \cup D_s R_k$ be the biedges of the FCRMs; $1 \ne j$, $1 \le i, j \le n$, $1 \le s \le p$ and $1 \le k \le m$.*

*Then the biedges form a directed bicycle. An FCRM is said to be bicyclic if it possesses a directed bicycle. An FCRM is said to be abicyclic if it does not possess any directed bicycle.*

**DEFINITION 3.6:** *An FCRM with bicycles is said to have a feedback.*

**DEFINITION 3.7:** *When there is a feed back in an FCRM i.e., when the causal relations flow through a cycle in a revolutionary way, the FCRM is called a dynamical bisystem.*

**DEFINITION 3.8:** *Let $\{ \overline{C_1 C_2}, \overline{C_2, C_3}, ..., \overline{C_{n-1} C_n} \} \cup \{D_i \ R_j \ (or \ R_j \ D_i) \mid 1 \leq i \leq p, \ 1 \leq j \leq m\}$ be a bicycle when $C_i \cup R_j$ (or $D_i$) is switched on and if the causality flows through the edges of the bicycle and if it again causes $C_i \cup R_j$ (or $D_i$) we say that the dynamical bisystem goes round and round. This is true for the binodes $C_i \cup R_j$ (or $D_i$) for $1 \leq i \leq n$ and $1 \leq j \leq m$ (or $1 \leq j \leq p$). The equilibrium bistate for the dynamical bisystem is called the hidden bipattern.*

**DEFINITION 3.9:** *If the equilibrium bistate of the dynamical bisystem is a unique bistate bivector then it is called a fixed bipoint.*

We shall just illustrate by a simple example.

***Example 3.3:*** Consider the FCRMs with $(C_1, ..., C_n) \cup (D_1, ..., D_p)$ (or $R_1, ..., R_m$) as binodes. For instant if we start the dynamical bisystem by switching $C_1 \cup R_1$ (or $D_1$) on. Let us assume the FCRM settles down with $C_1$ and $C_n \cup (R_1$ and $R_m)$ or $(D_1$ and $D_p)$ i.e., the state bivector remains as (1 0 0 0 0 … 0 1) $\cup$ (1 0 0 … 1) in R (or (1 0 0 … 1) in D). This state bivector is called the fixed bipoint.

It is to be noted in the case of FCRM we get a pair of fixed bipoint say $A = A_1 \cup D$ and $A = A_1 \cup R$; D denotes the state vector in the domain space of the FRM component of the FCRM described in the early part of this chapter. Likewise R denotes the state vector of the range space of the FRM component of the FCRM bimodel.

**DEFINITION 3.10:** *If the FCRM settles down with a bistate bivector repeating in the form*

$$A_1 \rightarrow A_2 \rightarrow ... \rightarrow A_1 \cup B_1 \rightarrow B_2 \rightarrow ... \rightarrow B_j \rightarrow B_1 \text{ (or } D_1 \rightarrow D_2 \rightarrow ... \rightarrow D_k \rightarrow D_1\text{)}$$ *then this equilibrium is called a limit bicycle.*

*Notation:* Suppose $A = A^1 \cup A^2$ is a bivector which is passed into a dynamical bisystem $E = E_1 \cup E_2$. Then $AE = A^1E_1 \cup A^2E_2 = (x'_1,...,x'_n) \cup (y'_1,...,y'_p)$ (or $(z'_1,...,z'_m)$ after thresholding and updating the bivector; suppose we get $(x_1, ..., x_n) \cup (y_1, ..., y_p)$ (or $(z_1, ..., z_m)$) we denote that by $(x'_1,...,x'_n) \cup (y'_1,...,y'_p)$ (or $(z'_1,...,z'_m)$) $\hookrightarrow (x_1, ..., x_n) \cup (y_1, ..., y_p)$ (or $(z_1, ..., z_m)$).

Thus the symbol '$\hookrightarrow$' means the resultant bivector has been thresholded and updated.

FCRMs have several advantages as well as some disadvantages:

1)  The main advantage is that this method gives the occurrence of the resultant in pairs by which stage by stage comparison is possible.

2)  When the data happens to be an upsupervised one the FCRMs model becomes handy.

3)  This bimodel also gives the hidden bipattern of the situation.

4)  This bimodel gives equal status to each and every expert who gives his / her opinion.

The only disadvantage is when we form the combined FCRM bimodel and if in that bimodel we have weightages 1 and –1 for some nodes we have the sum adding to zero. Thus at all times the connection bimatrices may not be comfortable for addition.

Combined conflicting opinion tend to cancel out and assisted by the strong law of large numbers. A consensus emerges as the sample opinion approximates the underlying population opinion. This problem will be easily overcome if in the FCRM bimodel the entries are only 0 and 1.

Here we mention some of the properties of the FCRM bimodel

1) FCRM are more applicable when the data in the first place is an unsupervised one. FCRMs work on the opinion of experts. FCRMs bimodel the world as a collection of biclasses and casual relation between biclasses.

2) FCRMs are fuzzy signed directed bigraphs with feedback. The directed biedge $e_{ij} = (e_{ij}^1) \cup (e_{ks}^2)$ from the bicausal concept $C_iC_j \cup (D_kR_s)$ (or $R_sD_k$) measures how much ($C_i$ causes $C_j$) $\cup$ ($D_k$ causes $R_s$ (or $R_s$ causes $D_k$)). The time varying concept bifunction $C_i(t) \cup R_s(t)$ (or $D_s(t)$) measures the non negative occurrence of some fuzzy event, perhaps the strength of a political sentiment, problems faced by persons with disability etc.

3) FCRM model can be used to model several types of problems. In this thesis we only model the problems faced by the persons with disability and the problems faced by the caretakers from the society relative to the persons with disability.

4) The bidedges $e_{ij} = (e_{ij}^1) \cup (e_{ks}^2)$ take values in the fuzzy causal biinterval $[-1, 1] \cup [-1, 1]$.

i)  $e_{ij} = 0 \cup 0$ indicates no causality between the binodes.

ii)  $e_{ij} > 0$ implies both $e_{ij}^1 > 0$ and $e_{ks}^2 > 0$; implies increase in the binodes $C_i \cup D_k$ (or $R_s$) implies increase in the binodes $C_j \cup R_s$ (or $D_s$).

iii)  $e_{ij} < 0$ implies both $e_{ij}^1 < 0$ and $e_{ks}^2 < 0$. Similarly decrease in the binodes $C_i \cup D_k$ (or $R_s$) implies decrease in the binodes $C_j \cup R_s$ (or $D_k$).

However unlike the FCM and FRM model we can have the following possibilities other than that of $e_{ij} = 0 \cup 0$, $e_{ij} > 0$ and $e_{ij} < 0$.

a.  $e_{ij} = e_{ij}^1 \cup e_{ks}^2$ can be such that $e_{ij}^1 = 0$ and $e_{ks}^2 > 0$. No relation in one binode and an increase in other node.

b.  In $e_{ij} = e_{ij}^1 \cup e_{ks}^2$ we can have $e_{ij}^1 = 0$ and $e_{ks}^2 < 0$. No causulality in the FCM node and decreasing relation in the FRM node.

c.  In $e_{ij} = e_{ij}^1 \cup e_{ks}^2$ we can have $e_{ij}^1 < 0$ and $e_{ks}^2 > 0$.

d.  In $e_{ij} = e_{ij}^1 \cup e_{ks}^2$ we can have $e_{ij}^1 < 0$ and $e_{ks}^2 = 0$.

e.  In $e_{ij} = e_{ij}^1 \cup e_{ks}^2$ we can have $e_{ij}^1 > 0$ and $e_{ks}^2 = 0$.

f.  In $e_{ij} = e_{ij}^1 \cup e_{ks}^2$ we can have $e_{ij}^1 > 0$ and $e_{ks}^2 < 0$.

Thus in case of FCRM we can have 9 possibilities where as in FCMs or FRMs we have only 3 possibilities. Thus the extra 6 possibilities can help in making the solution of the problem more sensitive or accurate.

We will illustrate the working for the bimodel described in of this chapter.

Using the bimatrix

$$M = M_1 \cup M_2 =$$

$$\begin{bmatrix} 0 & 0 & 1 & 0 & 0 \\ 0 & 0 & 1 & 1 & 1 \\ 1 & 0 & 0 & 0 & 0 \\ 0 & 0 & 0 & 0 & 1 \\ 0 & 1 & 0 & 0 & 0 \end{bmatrix} \cup \begin{bmatrix} 0 & 0 & 1 & 1 & 1 & 0 & 0 \\ 1 & 0 & 1 & 1 & 1 & 0 & 0 \\ 1 & 1 & 1 & 1 & 1 & 1 & 0 \\ 0 & 0 & 0 & 1 & 1 & 1 & 1 \\ 1 & 1 & 1 & 1 & 1 & 0 & 0 \end{bmatrix}$$

However we mention that this example is only an illustrate model to show how the dynamical bisystem functions.

Suppose $A = A_1 \cup A_2 = (1\ 0\ 0\ 0\ 0) \cup (0\ 0\ 1\ 0\ 0)$ be the bivector under consideration in which the node $F_1$ i.e., social status is in the on state in the FCM component of the FCRM and the node $F_3$ is in the on state in the FRM component of the FCRM.

To find the impact of A on the dynamical bisystem M.

$$\begin{aligned} AM \quad &= \quad (A_1 \cup A_2)\ (M_1 \cup M_2) \\ &= \quad A_1 M_1 \cup A_2 M_2 \\ &= \quad [(1\ 0\ 0\ 0\ 0) \cup (0\ 0\ 1\ 0\ 0)] \times \end{aligned}$$

$$\left\{ \begin{bmatrix} 0 & 0 & 1 & 0 & 0 \\ 0 & 0 & 1 & 1 & 1 \\ 1 & 0 & 0 & 0 & 0 \\ 0 & 0 & 0 & 0 & 1 \\ 0 & 1 & 0 & 0 & 0 \end{bmatrix} \cup \begin{bmatrix} 0 & 0 & 1 & 1 & 1 & 0 & 0 \\ 1 & 0 & 1 & 1 & 1 & 0 & 0 \\ 1 & 1 & 1 & 1 & 1 & 1 & 0 \\ 0 & 0 & 0 & 1 & 1 & 1 & 1 \\ 1 & 1 & 1 & 1 & 1 & 0 & 0 \end{bmatrix} \right\}$$

$$= (1\ 0\ 0\ 0\ 0) \begin{bmatrix} 0 & 0 & 1 & 0 & 0 \\ 0 & 0 & 1 & 1 & 1 \\ 1 & 0 & 0 & 0 & 0 \\ 0 & 0 & 0 & 0 & 1 \\ 0 & 1 & 0 & 0 & 0 \end{bmatrix} \cup$$

$$(0\ 0\ 1\ 0\ 0) \begin{bmatrix} 0 & 0 & 1 & 1 & 1 & 0 & 0 \\ 1 & 0 & 1 & 1 & 1 & 0 & 0 \\ 1 & 1 & 1 & 1 & 1 & 1 & 0 \\ 0 & 0 & 0 & 1 & 1 & 1 & 1 \\ 1 & 1 & 1 & 1 & 1 & 0 & 0 \end{bmatrix}$$

$$= [0\ 0\ 1\ 0\ 0] \cup [1\ 1\ 1\ 1\ 1\ 1\ 0].$$

Now after updating the first component of the bivector and as the bivectors are only 0's and 1's one need not threshold it.

Thus if $[1\ 0\ 1\ 0\ 0] \cup [1\ 1\ 1\ 1\ 1\ 1\ 0] = B_1 \cup B_2 = B$, then B $M_s^{t_2} = (B_1 \cup B_2)(M_1 \cup M_2^t)$. Like wise several bivectors was used to study the problem.

(By $M_s^{t_2}$ we mean the special transpose of the bimatrix, i.e., only the second matrix has been transposed and the other matrix remains the same i.e., only the second component of the matrix has been transposed and the first component remains as it is. $M_s^{t_1} = M_1^t \cup M_2$ we mean only the first component of the bimatrix alone is transposed the second component remains the same.

By $M^t = M_1^t \cup M_2^t$, both the component of the matrix M is being transposed).

$$= [1\ 0\ 1\ 0\ 0] \begin{bmatrix} 0 & 0 & 1 & 0 & 0 \\ 0 & 0 & 1 & 1 & 1 \\ 1 & 0 & 0 & 0 & 0 \\ 0 & 0 & 0 & 0 & 1 \\ 0 & 1 & 0 & 0 & 0 \end{bmatrix} \cup$$

$$[1\ 1\ 1\ 1\ 1\ 1\ 0] \begin{bmatrix} 0 & 1 & 1 & 0 & 1 \\ 0 & 0 & 1 & 0 & 1 \\ 1 & 1 & 1 & 0 & 1 \\ 1 & 1 & 1 & 1 & 1 \\ 1 & 1 & 1 & 1 & 1 \\ 0 & 0 & 1 & 1 & 0 \\ 0 & 0 & 0 & 1 & 0 \end{bmatrix}$$

$= [1\ 0\ 1\ 0\ 0] \cup [3\ 4\ 6\ 3\ 5]$.

Now we have to update and threshold only the second component of the bivector as the first component is in the acceptable form as it has only 0 and 1's and the on state of the $C_1$ is 1 with which we started. Infact we observe the value of the first component is equal to $B_1$ (only).

Thus we can say the FCM component of the FCRM bimodel is a fixed point. Now we update and threshold the second component and get

$$\begin{aligned} C &= C_1 \cup C_2 \\ &= (1\ 0\ 1\ 0\ 0) \cup (1\ 1\ 1\ 1\ 1) \\ &= B_1 \cup C_2. \end{aligned}$$

Now we do not work with first component of the binode. We find

$$CM = C_1 M_1 \cup C_2 M_2$$

$$= B_1 \cup (1\ 1\ 1\ 1\ 1) \begin{bmatrix} 0 & 0 & 1 & 1 & 1 & 0 & 0 \\ 1 & 0 & 1 & 1 & 1 & 0 & 0 \\ 1 & 1 & 1 & 1 & 1 & 1 & 0 \\ 0 & 0 & 0 & 1 & 1 & 1 & 1 \\ 1 & 1 & 1 & 1 & 1 & 0 & 0 \end{bmatrix}$$

$$= B_1 \cup (3\ 2\ 4\ 5\ 5\ 2\ 1)$$

We have to only threshold the second component which gives

$$D = D_1 \cup D_2$$
$$= B_1 \cup (1\ 1\ 1\ 1\ 1\ 1\ 1).$$

It has become essential to mention at least at this juncture just like in the FCM and FRM models described in [44], whose dynamical system can function only with state vectors which has to be zero or ones, i.e., the state vector to be in on or off state. For other numbers has no meaning in the state vectors. Likewise we see the FCRM bimodel can only recognize bistate bivectors.

$$X = X_1 \cup X_2$$

$$= (x_1^1,...,x_n^1) \cup (x_1^2,...,x_p^2) \ (\text{or} \ (x_1^2,...,x_m^2))$$

where $x_i^1 \in \{0, 1\}$ and $x_k^2$, $x_j^2 \in \{0, 1\}$, $1 \leq i \leq n$, $1 \leq k \leq p$, $1 \leq j \leq m$.

Only when the bivectors are in the on or off state alone can be dynamical bisystem recognize the bivectors. That is why the process of updating and thresholding the bivector by which all positive integers are replaced by 1 and negative integers by 0.

The zeros remain, as it is if at the starting the coordinates was zero.

Now we find $DM_s^{t_2} = (D_1 \cup D_2) \; (M_1 \cup M_2^t)$

$= \; D_1 M_1 \cup D_2 M_2^t$,

$$= B_1 \cup (1\;1\;1\;1\;1\;1\;1) \begin{bmatrix} 0 & 1 & 1 & 0 & 1 \\ 0 & 0 & 1 & 0 & 1 \\ 1 & 1 & 1 & 0 & 1 \\ 1 & 1 & 1 & 1 & 1 \\ 1 & 1 & 1 & 1 & 1 \\ 0 & 0 & 1 & 1 & 0 \\ 0 & 0 & 0 & 1 & 0 \end{bmatrix}$$

$= B_1 \cup (3\;4\;6\;4\;5)$.

After updating and thresholding the resultant bivector we get the hidden bipattern to be a pair of fixed bipoint given by $\{(1\;0\;1\;0\;0) \cup (1\;1\;1\;1\;1)\}$ and $\{(1\;0\;1\;0\;0) \cup (1\;1\;1\;1\;1\;1\;1)\}$.

The merits of this new FCRM bimodel is as follows:

1. FCRM bimodel allows experts to represent factual and evaluative concepts in an interactive frame.

2. Experts can quickly draw FCRM bimodel pictures or respond to the questionnaire.

3. Experts can consent or dissent to the local causal structure and perhaps the global equilibrium.

4. FCRM knowledge representation and inferencing structure reduces to simple bivector bimatrix operations, favours integrated circuit implementation.

5. Yet an FCRM equally encodes the experts' knowledge or ignorance, wisdom or prejudice.

It is important to mention here that we make use of the FCRM bimodel in a special way so that the comparison is possible at each stage.

We know the FCRM bimodel has two fuzzy models viz. the FCMs and FRMs as its components. Without loss of generality we can assume the first component of the FCRMs to be FCMs and always the second component to be the FRMs.

Further it is assumed in our book pertinently that in one model the nodes taken by the FCM to be identical with the nodes of the domain space of the FRM the second component of the FCRMs. Such a choice would certainly help for better comparison of the attributes on the dynamical bisystem. This is one of the special features of this new FCRMs which cannot be done using any other model. Thus this new model would cater to the needs of our analysis carried out in this book.

**Chapter Four**

# ANALYSIS OF THE SOCIO ECONOMIC PROBLEMS FACED BY THE PWDS USING SFCMS

This chapter has five sections. Section one gives a brief literature survey of fuzzy cognitive maps and their applications to unsupervised data. Further we give the notion of positive FCMs, which are FCMs that take only edge weights from the set {0, 1}. Here we define in this section one special type of FCMs, which are also positive FCMs and briefly how these models are used in analyzing social and psychological problems. Section two gives the expert opinion of the PWDs problems by 82 PWDS using special FCMs model. The corresponding expert opinion of 82 caretakers of the same PWDs is given using special FCMs is given in section three of this chapter. Section four gives the expert opinion of 12 NGO's using special FCM on the PWDs problem. In section five the combined special FCMs are used to analyse the problems of the PWDs. Further conclusions and suggestions are given in this section.

## 4.1 Application of FCMs models to Social and Psychological Problems

We proceed onto describe briefly the use of FCMs in studying the social problems. Cognitive maps or fuzzy cognitive maps are techniques that attempt to depict and to analyse the cognitive process of human thinking and human behaviour on specific domains by creating models. These models are represented as signed directed graphs of concepts and by the various causal relationships that exists between the concepts. They are mainly used for decision making and prediction. Fuzzy cognitive maps are used to analyse any real time system and study their behaviour with changes in the state of the system.

Axelrod proposed cognitive maps in 1976 as a formal tool for decision making. He had used the matrix representation of the directed graph to represent and study the social scientific knowledge. Cognitive maps were used by Maruyama to study the public health scenario of a society, where the edges of the map were assigned a positive or a negative value. Based on the Cognitive map structure in 1986, Bart Kasko proposed fuzzy cognitive maps.

In 1992 he had applied fuzzy cognitive maps to neural networks and studied the hidden patterns in combined and adaptive knowledge net works. In 1996 Craiger et al studied modeling organizational behavior with fuzzy cognitive maps. They had studied fuzzy cognitive maps as a computer simulation methodology for representing and predicting the behavior of models of arbitrary complexity. In 1996, Tsadiras and Margaritis studied the application of using Certain Neurons in fuzzy cognitive maps. In 1997 they had defined Certain Neurons Fuzzy Cognitive Maps and introduced time and decay mechanism.

Fuzzy Cognitive Maps (FCMs) had been applied to transportation problem to predict the route utility and to analyse the passenger preference in 2000 by Vasantha and V. Indra.

Vasantha and S. Florentein have used FCMs to analyse the social problems of the migrant labourers living with  HIV/AIDs in 2004.  They have used FCMs to study the views of Periyar on untouchability (in 2005) and study the views of the public about Vedic mathematics (in 2006).  Fuzzy Cognitive Maps have been used to study behaviour of stock market, drug trafficking etc.

Here we proceed onto define the notion of positive FCM model and special FCM model

**DEFINITION 4.1.2.1:**  *If in the FCMs the edge weights take values only from the set {0, 1} then we call the FCM to be a positive FCM. If E is the connection matrix associated with the positive FCM, E is called the positive connection matrix of the positive FCM.*

We see the class of positive FCMs from a subclass of FCMs.

**Definition 4.1.2.2:** *Let $E_1$, ..., $E_n$ be a collection of positive connection matrices associated with positive FCMs on the same set of t-nodes i.e., each $E_i$ is a $t \times t$ positive matrix, $1 \leq i \leq n$.*

*Let $E = \sum_{i=1}^{n} E_i$ be the combined positive FCM. If each entry in $E = (e_{ij})$ is divided by n and if  $\dfrac{e_{ij}}{n}$ greater than or equal to 0.5 replace it by 1  and if $\dfrac{e_{ij}}{n}$ is less than 0.5 replace it by 0.*

*Then the transformed matrix $E_s = \dfrac{E}{n}$ is called the positive special connection matrix or in short special connection matrix of the set of FCMs $E_1$, ..., $E_n$. The FCM associated with this* *positive special connection matrix is defined to be the special FCM.*

## 4.2 Expert opinion of the PWDs using Special FCM to analyse the problem of PWDs

We have applied the tool of FCMs to study and analyse the economic problems faced by PWDs in rural areas of Tamil Nadu. Since we have problems in defining the very concept of disability one can understand how it is impossible to describe fully their problems using statistics. Further the shame and stigma associated with disability in India exacerbates the problem.

From the survey and the interviews taken from the rural PWDs and the caretakers were the vital problems suffered by them.

1. No proper healthcare: They said their economic conditions forced them to neglect the PWDs from giving them proper medical care and or taking them for routine check up or even taking them to local doctor at time of emergency. For to take them out for medical aid one has to spend also on the conveyance which is an economic burden on the caretakers. Thus healthcare is one of the problems faced by the PWDs.

2. Poor nutrition: Due to poor economic conditions or reasons best known to the caretakers it was very difficult for the caretakers as well as family members of the PWDs to give "good" nutritious food to the PWDs. Only the male breadwinner alone was in a position to get somewhat nutritious food. It was also observed that not only the PWDs but also the family members suffered from lack of nutrition.

3. Improper clothing: When the very basic needs such as food and medical aid were under strain one cannot think of comfortable clothing for the disabled to comfort his / her disability. Thus clothes depending on the weather conditions cannot be provided by the caretakers to the dependent PWDs.

Baring this as the PWDs are assumed not only as a burden but a curse on the family and a social stigma they were not given new cloths even during the festivals. Infact they were made to wear the old dresses available in the family of the caretakers. Also they were not given frequent changes of the dress, which we observed to be dirty and torn. Further this poor hygiene made them unhealthy and gave a distressed out look.

4. No proper shelter: The PWDs from the rural poor were mostly segregated from the family in the verandah or in some common places such as temples or in the back yard under a thatch for reasons best known to the caretakers. They were alone in most cases and found to be lonely, sad and depressed.

5. No recreation: The sole recreation in villages as of today is watching TV (Television) or hearing music. As they (PWDs) happen to hold stigma they are never taken out for any public function, even for festivals to temples / churches or for recreation to cinemas, parks, picnic, tours or visits to relatives. They are made to be at home as watchman when the entire family leaves the house for a function or festival or recreation.

Games / painting or indoor games etc can certainly be provided to them if at all they have better economy. For in rich society PWDs are given nice recreations.

6. No proper school education: When they are denied the basic needs the question of sending them to school under these circumstances can only be construed as a waste of time for the caretakers. That is why there were several PWDs who had not even entered the school. They would by no means be interested in giving a special training depending on the disability. Thus the skills, talents and the gifts in these economically poor disabled never come to light. Yet another major reason is both the caretakers and the PWDs felt it to be a disgrace as the disability in a family was considered to be a curse and so they suffered a social stigma.

7.   No employment / self-employment:  The poor PWDs can be technically trained in certain jobs fetching skills depending on their deformity so that they can sustain themselves and lead a dignified life.  But unfortunately this is absent in most of the cases.  They are treated as if they are some condemned non-living objects.  The caretakers and the family fail to make use fo them in a positive way.

8.   No information about the SSHGs:  The information about Special Self-Help Groups (SSHG) and the structure and functional mechanism of SSHGs are unknown to PWDs.  It never reaches them the reasons being that the rural PWDs are ignorant as uneducated and segregated.  The caretakers never want the public or anyone to know about the presence of a PWD in their family as they view it as a social stigma for the marriages in these families will turn out to be  a problem.  Further the PWD in their family is viewed as a curse so they tend to isolate them consequent to which the information about SSHG or benefits announced by government or their rights as PWDs or free medical aid etc is never known to them.  So the question of utilizing these never occurs.  Further even if the information is found in dailies or advertisements in TVs PWDs are ignorant of it as they do not have any access to them.

9.   Welfare  measures  of  government  never  reach  the  rural PWDs: Welfare measures never reach  the rural poor PWDs for the reasons said in above point.

10. Poor economy:  The economic condition as per the survey done by the NGO LAMP NET is that only 2.1% of the families with a PWD as its member have an income over Rs. 3000 per month.  75.9% of them get an income less than Rs. 1000 per month.

11. Marriage remains a question mark:  When the very basic needs are not given it has become a near impossibility for the caretakers to even think of marriage for the PWDs under their custody.  This is evident from the survey data that in the reproductive age, they remain  unmarried where as in the unproductive age the percentage of married PWDs is considerably larger than the percentage of married in the reproductive age.

Another important reason is 'social stigma'. It is difficult to get them married off due to this major social hurdle. Poverty is also another reasons for the PWD's remaining unmarried.

Mainly these were the attributes given by the experts as problem of the PWDs. The group of experts were drawn from the PWD's caretakers / guardians / parents of PWDs and NGOs.

We took the experts opinions of 100 PWDs and their caretakers. However only 82 of the responses survived as PWDs and their caretakers. The remaining 18 had to be discarded as either the caretaker or the PWDs failed to give proper opinion. When we say proper opinion it mathematically implies that using the information given by them a dynamical system cannot be constructed. Now we have carefully worded the attributes so that none of the edge $e_{ij}$ of the attributes can be negative. Now one could think of using combined FCMs if the number of experts was a small number. So the chance of using combined FCMs was ruled out due to time consuming calculations. To overcome all these types of shortcomings we constructed or modified the combined FCM model and we have termed this new FCM model as the special FCM model. We found the combined FCM of the 82 experts by adding the 82 connection matrices. Then we divided each and every entry of the combined connection matrix by 82. If the entry was greater than or equal to 0.5 we replaced it by 1. All values less than 0.5 was replaced by 0. Now this combined FCM, which we define as special FCM is just like a FCM and functions exactly like the ordinary FCM.

**Justification for using special FCMs:** In the first place the problems faced by the PWDs due to poverty cannot be analyzed using the survey data as it is incapable of captivating the intricate and the inexpressible feelings associated with them. Since the data is an unsupervised one we can use only fuzzy models. As the attributes cannot be divided as disjoint classes only FCM alone can help. Further this model alone can give the hidden pattern of any attribute. Finally as the issue involves lots of sensitive problems only fuzzy cognitive model alone can captivate the feelings and the inter relation between the attributes / problems.

We give one integrated model associated with each of the groups of experts. This is made possible by special FCMs. We have taken three groups viz. 82 PWDs, 82 caretakers and 12 NGOs who solely work for the PWDs. Also it is evident from the attributes given by them. The economic status of the PWD's drastically affects not only their day to day life but also the longstanding attributes like unemployment, marriage etc. The problems of PWDs is also due to the fact that they are a very small percentage and hence it will not affect the vote bank of the political leaders, so like SHG for women the SSHG is not given importance and developed by the politicians. No one cares for the rural poor; more so for the rural poor with disability.

We give here the opinion of the 82 PWDs, all of them were aged above 15 years. Some of them were also educated and employed / self-employed. They were suffering from all types of disabilities. First they said for education of their children in a common school may not help them in any way. Instead it would make them only more depressed. All of them opted for special education depending on their disability. They wanted an education, which would sustain them economically to lead an independent life. They felt regular education has nothing to do with their employment.

They also accepted that it was very difficult for them to go to school. So they did not want to blame the caretakers or parents. They were interested to learn to read and write. More than 59% of them did not know to read or write. Further they blamed the caretakers saying that they did not get uniform treatment. It is their (caretakers) mindset which in the first place stops them from spending on the PWDs. Might be they may think they cannot expect any returns for their investments on education of the PWDs. In the first place they are denied good food and shelter leave alone the school education. This is well substantiated from our data that 59% of the PWDs have not entered the school premises.

The illiteracy among the women PWDs is much higher than that of the male PWD's. Also it is to be noted from the data that out of the total school dropouts 18% of the school dropouts are from the PWDs. The educational status is very frustrating among the PWDs. Among the reasons attributed by them for their neglect is that the caretakers don't find time to spend on them for they are always in a hurry to earn (as most of them are daily or weekly wagers) as over 75 percent of them earn less than Rs. 1000 a month. Abstaining from work would directly deny the family even this paltry income. Thus at the outset we feel the economic condition has a very great direct impact on the rural poor PWDs. The expert now using the 11 attributes described in section one gives the connection matrix associated with it.

We have taken these eleven attributes as they feel the economic status of the PWD's or their caretakers directly depend on their health care, nutrition, clothing, shelter, recreation, education, employment and so on. So we have taken these attributes to find the interrelation between them. The study pertains to the rural PWDs, as the economic condition of them is questionably poor.

Expert opinion of 82 PWDs who were literates and could follow our questions and answered well is taken as the special FCM model. As it would not be easy to discuss all the 82 connection matrices and the resultants we first found the 82-11 × 11 matrices associated with the 82 experts. We added all the 82 connection matrices and divided each of the term of the sum of the matrix by 82. If the value was greater than or equal to 0.5 we replaced it by 1; if it was less than 0.5, we replaced it by zero.

Thus the resultant matrix was also a connection matrix only with zeros and ones. We have avoided –1 to be an entry in all the connection matrices. This is appropriately done by using the proper adjectives to describe the attributes. So by this we have avoided the occurrence of –1 or to be more exact; $e_{ij} \in \{0,1\}$ i.e., the edge weights $e_{ij}$ are either 0 and 1. We define the combined FCM of the n matrices after dividing each component by n and replacing the values which are greater than or equal to 0.5 by 1 and less than 0.5 by 0 as the special connection matrix and the corresponding combined FCM is defined in page 49 of this book.

We now proceed onto give the special connection matrix of the special FCM model.

Let M denote the $11 \times 11$ special connection matrix of the special FCM model. M is calculated and is given as;

$$M = \begin{bmatrix} 0 & 1 & 0 & 0 & 0 & 0 & 0 & 0 & 0 & 1 & 0 \\ 1 & 0 & 0 & 0 & 0 & 0 & 0 & 0 & 0 & 1 & 0 \\ 0 & 0 & 0 & 0 & 0 & 0 & 0 & 0 & 0 & 1 & 0 \\ 0 & 0 & 0 & 0 & 0 & 0 & 0 & 0 & 0 & 0 & 0 \\ 0 & 0 & 0 & 0 & 0 & 0 & 0 & 0 & 0 & 0 & 0 \\ 0 & 0 & 0 & 0 & 0 & 0 & 0 & 0 & 0 & 0 & 0 \\ 0 & 0 & 0 & 0 & 0 & 0 & 0 & 0 & 0 & 1 & 1 \\ 0 & 0 & 0 & 0 & 0 & 0 & 0 & 0 & 0 & 0 & 0 \\ 0 & 0 & 0 & 0 & 0 & 1 & 0 & 0 & 0 & 0 & 0 \\ 1 & 1 & 1 & 0 & 0 & 0 & 0 & 0 & 0 & 0 & 0 \\ 0 & 0 & 0 & 0 & 0 & 0 & 1 & 0 & 0 & 0 & 0 \end{bmatrix}$$

Now we find the effect of the on state of poor economy and find its effect on the dynamical system.

Suppose X = (0 0 0 0 0 0 0 0 0 1 0) i.e., only poor economy is kept in on state. XM = (1 1 1 0 0 0 0 0 0 0 0) after updating XM we get

$$X_1 = (1\ 1\ 1\ 0\ 0\ 0\ 0\ 0\ 0\ 1\ 0).$$

Now

$$X_1M \hookrightarrow (1\ 1\ 1\ 0\ 0\ 0\ 0\ 0\ 0\ 1\ 0)$$
$$= X_2 \text{ (say)}.$$

(where '$\hookrightarrow$' the resultant vector has been updated and thresholded).

$$X_2M \hookrightarrow (1\ 1\ 1\ 0\ 0\ 0\ 0\ 0\ 0\ 1\ 0)$$
$$= X_3 \text{ (say)}.$$

But as $X_2 = X_3$, we see the hidden pattern of the state vector X is a fixed point. Thus the hidden pattern of the state vector X gives one the information that poor economy implies; they cannot have proper health care, they suffer from poor nutrition and they have improper clothing.

However poor economy has nothing to do with no education or no employment according to the PWDs.

Next the expert's wishes to study the effect of the on state "No proper health care" on the dynamical system M.

Let

$$T = (1\ 0\ 0\ 0\ 0\ 0\ 0\ 0\ 0\ 0\ 0)$$
$$TM = (0\ 1\ 0\ 0\ 0\ 0\ 0\ 0\ 0\ 1\ 0)$$
$$= T_1 \text{ (say)}.$$

$$T_1M \hookrightarrow (1\ 1\ 1\ 0\ 0\ 0\ 0\ 0\ 0\ 1\ 0)$$
$$= T_2 \text{ (say)}.$$

$$T_2M \hookrightarrow (1\ 1\ 1\ 0\ 0\ 0\ 0\ 0\ 0\ 1\ 0)$$
$$= T_3\ (= T_2).$$

Thus the hidden pattern is a fixed point. It is important to observe that the poor health care is dependent on the poor economy.

Suppose we wish to study the effect of the node "No proper school education", i.e. let the state vector be S = (0 0 0 0 0 1 0 0 0 0 0). To find the effect of S on M

$$SM = \quad (0\ 0\ 0\ 0\ 0\ 0\ 0\ 0\ 0\ 0\ 0)$$

after updating we get

$$SM = \quad (0\ 0\ 0\ 0\ 0\ 1\ 0\ 0\ 0\ 0\ 0).$$

According to the PWD the school education is of no use to them. This is also evident when they spoke about education. Suppose the node marriage alone is in the on state and all other nodes are in the off state, to study the effect of it on the dynamical system M.

Let X = (0 0 0 0 0 0 0 0 0 0 1) be the state vector with the node "marriage remains a question mark" is in the on state. To find the effect of X on M we find XM = (0 0 0 0 0 1 0 0 0 0 0), to find the effect of Y on M.

$$\begin{aligned} YM \quad &\hookrightarrow \quad (0\ 0\ 0\ 0\ 0\ 0\ 1\ 0\ 0\ 1\ 1) \\ &= \quad Y_1 \text{ (say)}. \end{aligned}$$

Now we find

$$\begin{aligned} Y_1M \quad &\hookrightarrow \quad (1\ 1\ 1\ 0\ 0\ 0\ 1\ 0\ 0\ 1\ 1) \\ &= \quad Y_2 \text{ (say)}. \end{aligned}$$

Now the effect of $Y_2$ on the dynamical system M gives

$$\begin{aligned} Y_2M \quad &\hookrightarrow \quad (1\ 1\ 1\ 0\ 0\ 0\ 1\ 0\ 0\ 1\ 1) \\ &= \quad Y_3 \ ( = Y_2). \end{aligned}$$

Thus the hidden pattern is a fixed point.

It is however very surprising to see the impact of "marriage remains a question mark" makes on state no proper health care, poor nutrition, improper clothing and no employment / no self-employment. Thus we see the special FCM model in a way stabilizes only a few nodes and we can say only the nodes which has got an approval atleast from more than 40 alone sustain and stays as a non zero component. However the conclusions will be mentioned in the last chapter.

## 4.3 The Expert Opinion of the Caretakers of the PWDs using SFCM

In this section we give the opinion of 82 caretakers of the 82 PWDs mentioned in section 4.2 of chapter four of this book. We have used the same model as that of the PWD viz. the special FCM model. The same procedure mentioned in the section 4.2 is used in this section also.

The main observation of our interview was that the caretakers said they were forced to neglect their PWD due to poor economy, most of the daily wagers said if they do not go to work they will not get their wages.

So some of them kept the food near their Children With Disability (CWDs) and without attending for their other cares left for work. For they said they were unable to find time to cater to the needs of the CWDs. Only on Sundays or the days they did not go to work they could give some attention.

They said their neglect did not really mean they did not like their children. They liked them even if they are CWD. But their poor economy was the main cause of neglect.

Now we give the special connection matrix T associated with the caretakers.

$$T = \begin{bmatrix} 0 & 1 & 0 & 0 & 0 & 0 & 0 & 0 & 0 & 1 & 0 \\ 1 & 0 & 0 & 0 & 0 & 0 & 0 & 0 & 0 & 1 & 0 \\ 0 & 0 & 0 & 0 & 0 & 0 & 0 & 0 & 0 & 1 & 0 \\ 0 & 0 & 0 & 0 & 0 & 0 & 0 & 0 & 0 & 1 & 0 \\ 0 & 0 & 0 & 0 & 0 & 0 & 0 & 0 & 0 & 1 & 0 \\ 0 & 0 & 0 & 0 & 0 & 0 & 0 & 0 & 0 & 0 & 0 \\ 0 & 0 & 0 & 0 & 0 & 0 & 0 & 0 & 0 & 1 & 1 \\ 0 & 0 & 0 & 0 & 0 & 0 & 0 & 0 & 0 & 0 & 0 \\ 0 & 0 & 0 & 0 & 0 & 1 & 0 & 0 & 0 & 0 & 0 \\ 1 & 1 & 1 & 1 & 1 & 0 & 0 & 0 & 0 & 0 & 1 \\ 0 & 0 & 0 & 0 & 0 & 0 & 0 & 0 & 0 & 1 & 0 \end{bmatrix}.$$

Suppose the node "no proper health care" alone is in the on state and all other nodes are in the off state. To find the effect of the node on the dynamical system T.

Let

$$X \quad = \quad (1\ 0\ 0\ 0\ 0\ 0\ 0\ 0\ 0\ 0\ 0)$$

Now

$$XT \quad = \quad (0\ 1\ 0\ 0\ 0\ 0\ 0\ 0\ 0\ 1\ 0).$$

After updating we get
$$X_1 \quad = \quad (1\ 1\ 0\ 0\ 0\ 0\ 0\ 0\ 0\ 1\ 0).$$

Now

$$X_1 M \quad \hookrightarrow \quad (1\ 1\ 1\ 1\ 1\ 1\ 0\ 0\ 0\ 1\ 1)$$
$$= \quad X_2 \text{ (say)}.$$

$$X_2 T \quad \hookrightarrow \quad (1\ 1\ 1\ 1\ 1\ 1\ 0\ 0\ 0\ 1\ 1)$$
$$= \quad X_3\ (= X_2).$$

Thus the hidden pattern of the state vector X in T is a fixed point. It is seen from the dynamical system T, that no proper health care is directly dependent on poor nutrition and poor economy.

The three factors poor economy, poor nutrition and poor health care give way to improper clothing, no proper shelter, no recreation and the marriage remains a question mark. However they agree with the CWDs that their school education has nothing to do with economy for they are clear that their schooling is not gong to give them any form of good employment mainly because they are in the rural areas with poverty dominating them. Their mind refuses to think of the self-employment or employment of the PWDs for to be self employed also one need to have some initial capital they said. It is very clear the PWDs and their caretakers do not think formal school education, which is free up to 18 years can change their outlook.

Next we proceed onto study the on state of the node poor economy alone and all other nodes are in the off state. Let

$$X = (0\ 0\ 0\ 0\ 0\ 0\ 0\ 0\ 0\ 1\ 0).$$

To find the effect of the state vector X on the dynamical system T.

$$XT = (1\ 1\ 1\ 1\ 1\ 0\ 0\ 0\ 0\ 0\ 1).$$

after updating we get

$$X_1 = (1\ 1\ 1\ 1\ 1\ 0\ 0\ 0\ 0\ 1\ 1).$$

We find the effect of $X_1$ on T

$$X_1 M = (1\ 1\ 1\ 1\ 1\ 0\ 0\ 0\ 0\ 1\ 1)$$

$$= X_2 \text{ (say); but } X_2 = X_1.$$

Thus the hidden pattern of X is a fixed point. We see the poor economy of the caretakers is the cause of poor (or no) health care, poor nutrition, no good clothing, poor shelter, no recreation and marriage a question mark. Here also as they do not bother about school education poverty has no impact on school education.

We found that their busy schedule has forced them to neglect the very knowledge of SSGH or welfare measure of the government, which never reaches them. Further our study reveals both the PWDs and the caretakers are least interested about school education. They want only that they of vocational training which would make the PWD earn or eke out his / her living by some self-employment. The reason one can attribute to this is the social stigma they face in their society.

One of the basic conclusions we get is that they don't wish the PWDs to get employed under any one for they are afraid of ill treatment and stigma. We are unable to understand their ideology about this.

## 4.4    The expert opinion of the NGOs working for the PWDs using Special FCM model

We have obtained the opinion of 12 NGOs working in different organizations like ROSE, LAMP, SERVITES, VDS, MAKAM, WCDS, SHARE AND CARE etc. We have given them the 11 attributes described in the section one of this chapter. We have also used the new special FCM model described in section two of chapter four. The NGOs mainly complained that the caretakers were not very keen to get their help. The cause of their ignorance may be due to lack of time, poor economy and illiteracy of the caretakers and social stigma.

Now we give the $11 \times 11$ special connection matrix of the special FCM model given by the NGOs.

$$N = \begin{bmatrix} 0 & 1 & 0 & 0 & 0 & 0 & 0 & 1 & 1 & 1 & 0 \\ 0 & 0 & 0 & 0 & 0 & 0 & 0 & 1 & 0 & 1 & 0 \\ 0 & 0 & 0 & 0 & 0 & 0 & 0 & 0 & 1 & 1 & 0 \\ 0 & 0 & 0 & 0 & 0 & 0 & 0 & 0 & 0 & 0 & 0 \\ 0 & 0 & 0 & 0 & 0 & 0 & 0 & 0 & 0 & 1 & 0 \\ 0 & 0 & 0 & 0 & 0 & 0 & 0 & 1 & 1 & 0 & 0 \\ 0 & 0 & 0 & 0 & 0 & 0 & 0 & 1 & 1 & 1 & 1 \\ 0 & 0 & 0 & 0 & 0 & 0 & 1 & 1 & 0 & 0 & 0 \\ 0 & 0 & 0 & 0 & 0 & 0 & 1 & 0 & 0 & 1 & 0 \\ 0 & 0 & 0 & 0 & 0 & 0 & 0 & 0 & 0 & 0 & 1 \\ 0 & 0 & 0 & 0 & 0 & 0 & 1 & 0 & 0 & 0 & 0 \end{bmatrix}$$

Now we find the effect of the on state of the vector "No proper health care" of PWDs, and all other nodes are in the off state. Let

$$X \quad = \quad (1\,0\,0\,0\,0\,0\,0\,0\,0\,0\,0)$$

be the given state vector. To find the hidden pattern of X using N.

$$XN \quad = \quad (0\,1\,0\,0\,0\,0\,0\,1\,1\,1\,0).$$

After updating XN let the resultant be $X_1$.

$$X_1 \quad = \quad (1\,1\,0\,0\,0\,0\,0\,1\,1\,1\,0).$$

$$X_1 M \quad \hookrightarrow \quad (1\,1\,0\,0\,0\,0\,1\,1\,1\,1\,1)$$
$$= \quad X_2 \text{ (say)}.$$

$$X_2 T \quad \hookrightarrow \quad (1\,1\,0\,0\,0\,0\,1\,1\,1\,1\,1)$$
$$= \quad X_3 \text{ (say)}.$$

Thus $X_3 = X_2$. The hidden pattern is a fixed point; we see according to the NGOs the "no proper health care" of the PWDs is due to poor nutrition, no employment, no information about the SHG, welfare measures of government never reaches them, poor economy and marriage remains a question rank.

However the NGOs said that when they see some of the PWDs, they are not only depressed and insecure they showed several signs of mental problems which needed proper counseling.

Next they study the effect of the node 7 i.e., "no employment / no self-employment of the PWDs. The effect of it on them and the caretakers.

Let

$$T \quad = \quad (0\ 0\ 0\ 0\ 0\ 0\ 1\ 0\ 0\ 0\ 0)$$

be the given state vector. To find the effect of T on N.

$$TN \quad \hookrightarrow \quad (0\ 0\ 0\ 0\ 0\ 0\ 1\ 1\ 1\ 1\ 1)$$
$$= \quad T_1 \text{ (say)}.$$

$$T_1N \quad \hookrightarrow \quad (0\ 0\ 0\ 0\ 0\ 0\ 1\ 1\ 1\ 1\ 1)$$
$$= \quad T_2 \text{ (say)}.$$

We see $T_1 = T_2$, thus the hidden pattern is a fixed point. We see no employment / no self-employment is the cause of no information about SHG, welfare measures of government never reaches them, they suffer from poor economy and their marriage remains a question mark.

Suppose the experts wish to study the effect of the on state of the poor economy alone on the dynamical system N.

Let $\quad$ Y $\quad = \quad (0\ 0\ 0\ 0\ 0\ 0\ 0\ 0\ 0\ 1\ 0)$

be the given state vector.

Now

$$YN \quad \hookrightarrow \quad (0\ 0\ 0\ 0\ 0\ 0\ 0\ 0\ 0\ 1\ 1)$$
$$= \quad Y_1 \text{ (say).}$$

$$Y_1N \quad \hookrightarrow \quad (0\ 0\ 0\ 0\ 0\ 0\ 1\ 0\ 0\ 1\ 1)$$
$$= \quad Y_2 \text{ (say).}$$

$$Y_2N \quad \hookrightarrow \quad (0\ 0\ 0\ 0\ 0\ 0\ 1\ 1\ 1\ 1\ 1)$$
$$= \quad Y_3 \text{ (say).}$$

$$Y_3N \quad \hookrightarrow \quad (0\ 0\ 0\ 0\ 0\ 0\ 1\ 1\ 1\ 1\ 1)$$
$$= \quad Y_4 \text{ (say).}$$

Since $Y_3 = Y_4$, we get the hidden pattern of Y to be a fixed point. Thus poor economy is the cause of no employment / self - employment, no information about SSHG, welfare measures of government never reaches them and the marriage of the PWDs remains a question mark.

However we see the NGOs did not give an interlinking dynamical system with no shelter or improper clothing or no recreation or no school education. However several of the NGOs were only interested in giving the PWDs vocational training depending on their disability which would for certain lead to some self employment of them.

But when we verified the 12 matrices we saw only 5 of them said school education was important for the development of the PWDs.

Now we give in the final section of this chapter the combined opinion of the three sets of experts viz. the PWD, the caretakers of the PWDs and the NGO who work only for the PWDs.

## 4.5 Combined special FCMs to study the problems of PWDs and conclusions drawn from our study

In this section we first define the notion of the combined special FCMs. Special FCMs were introduced in section one of chapter four. Now we define combined special FCMs analogous to combined FCMs.

If $M_1$, $M_2$, …, $M_p$ are p connection special matrices of a special FCM with n rows and n columns where n denotes the total number of concepts under study. The combined special FCM is got by adding all the p matrices $M_1$, …, $M_p$ i.e.,

$M = M_1 + … + M_p$.

Here $p = 3$ for we have take the opinion from three groups of experts given in section 4.2, 4.3 and 4.4 which are given by the special connection matrices M, T and N respectively.

Let

$W = M + T + N$

$$= \begin{bmatrix} 0 & 1 & 0 & 0 & 0 & 0 & 0 & 0 & 0 & 1 & 0 \\ 1 & 0 & 0 & 0 & 0 & 0 & 0 & 0 & 0 & 1 & 0 \\ 0 & 0 & 0 & 0 & 0 & 0 & 0 & 0 & 0 & 1 & 0 \\ 0 & 0 & 0 & 0 & 0 & 0 & 0 & 0 & 0 & 0 & 0 \\ 0 & 0 & 0 & 0 & 0 & 0 & 0 & 0 & 0 & 0 & 0 \\ 0 & 0 & 0 & 0 & 0 & 0 & 0 & 0 & 0 & 0 & 0 \\ 0 & 0 & 0 & 0 & 0 & 0 & 0 & 0 & 0 & 1 & 1 \\ 0 & 0 & 0 & 0 & 0 & 0 & 0 & 0 & 0 & 0 & 0 \\ 0 & 0 & 0 & 0 & 0 & 1 & 0 & 0 & 0 & 0 & 0 \\ 1 & 1 & 1 & 0 & 0 & 0 & 0 & 0 & 0 & 0 & 0 \\ 0 & 0 & 0 & 0 & 0 & 0 & 1 & 0 & 0 & 0 & 0 \end{bmatrix} +$$

$$
\begin{bmatrix}
0 & 1 & 0 & 0 & 0 & 0 & 0 & 0 & 0 & 1 & 0 \\
1 & 0 & 0 & 0 & 0 & 0 & 0 & 0 & 0 & 1 & 0 \\
0 & 0 & 0 & 0 & 0 & 0 & 0 & 0 & 0 & 1 & 0 \\
0 & 0 & 0 & 0 & 0 & 0 & 0 & 0 & 0 & 1 & 0 \\
0 & 0 & 0 & 0 & 0 & 0 & 0 & 0 & 0 & 1 & 0 \\
0 & 0 & 0 & 0 & 0 & 0 & 0 & 0 & 0 & 0 & 0 \\
0 & 0 & 0 & 0 & 0 & 0 & 0 & 0 & 0 & 1 & 1 \\
0 & 0 & 0 & 0 & 0 & 0 & 0 & 0 & 0 & 0 & 0 \\
0 & 0 & 0 & 0 & 0 & 1 & 0 & 0 & 0 & 0 & 0 \\
1 & 1 & 1 & 1 & 1 & 0 & 0 & 0 & 0 & 0 & 1 \\
0 & 0 & 0 & 0 & 0 & 0 & 0 & 0 & 0 & 1 & 0
\end{bmatrix} +
$$

$$
\begin{bmatrix}
0 & 1 & 0 & 0 & 0 & 0 & 0 & 1 & 1 & 1 & 0 \\
0 & 0 & 0 & 0 & 0 & 0 & 0 & 1 & 0 & 1 & 0 \\
0 & 0 & 0 & 0 & 0 & 0 & 0 & 0 & 1 & 1 & 0 \\
0 & 0 & 0 & 0 & 0 & 0 & 0 & 0 & 0 & 0 & 0 \\
0 & 0 & 0 & 0 & 0 & 0 & 0 & 0 & 0 & 1 & 0 \\
0 & 0 & 0 & 0 & 0 & 0 & 0 & 1 & 1 & 0 & 0 \\
0 & 0 & 0 & 0 & 0 & 0 & 0 & 1 & 1 & 1 & 1 \\
0 & 0 & 0 & 0 & 0 & 0 & 1 & 1 & 0 & 0 & 0 \\
0 & 0 & 0 & 0 & 0 & 0 & 1 & 0 & 0 & 1 & 0 \\
0 & 0 & 0 & 0 & 0 & 0 & 0 & 0 & 0 & 0 & 1 \\
0 & 0 & 0 & 0 & 0 & 0 & 1 & 0 & 0 & 0 & 0
\end{bmatrix}
$$

$$= \begin{bmatrix} 0 & 3 & 0 & 0 & 0 & 0 & 0 & 1 & 1 & 3 & 0 \\ 2 & 0 & 0 & 0 & 0 & 0 & 0 & 1 & 0 & 3 & 0 \\ 0 & 0 & 0 & 0 & 0 & 0 & 0 & 0 & 1 & 3 & 0 \\ 0 & 0 & 0 & 0 & 0 & 0 & 0 & 0 & 0 & 1 & 0 \\ 0 & 0 & 0 & 0 & 0 & 0 & 0 & 0 & 0 & 2 & 0 \\ 0 & 0 & 0 & 0 & 0 & 0 & 0 & 1 & 1 & 1 & 0 \\ 0 & 0 & 0 & 0 & 0 & 0 & 0 & 1 & 1 & 2 & 2 \\ 0 & 0 & 0 & 0 & 0 & 0 & 1 & 1 & 0 & 0 & 0 \\ 0 & 0 & 0 & 0 & 2 & 0 & 1 & 0 & 0 & 1 & 0 \\ 2 & 2 & 2 & 1 & 1 & 0 & 0 & 0 & 0 & 0 & 2 \\ 0 & 0 & 0 & 0 & 0 & 0 & 2 & 0 & 0 & 1 & 0 \end{bmatrix}.$$

Now W is the combined special connection matrix of the combined special FCM.

Now we study the effect of the state vector in which only the node "no proper health care" is in the on state. Let

$X = (1\ 0\ 0\ 0\ 0\ 0\ 0\ 0\ 0\ 0\ 0)$ be the state vector. The effect of X on the dynamical system M is given by

$$XM \hookrightarrow (0\ 3\ 0\ 0\ 0\ 0\ 0\ 1\ 1\ 3\ 0)$$

after updating and thresholding XM we get

$$X_1 \quad = \quad (1\ 1\ 0\ 0\ 0\ 0\ 0\ 1\ 1\ 1\ 0).$$

Now

$$X_1 M \quad = \quad (4\ 5\ 2\ 1\ 3\ 0\ 2\ 2\ 1\ 7\ 2)$$

after thresholding $X_1 M$. Let

$$X_2 \quad = \quad (1\ 1\ 1\ 1\ 1\ 0\ 1\ 1\ 1\ 1\ 1).$$

Now we effect of $X_2$ on M is given by

$$X_2M \quad \hookrightarrow \quad (4\ 5\ 2\ 1\ 3\ 0\ 4\ 4\ 3\ 16\ 4).$$

It is very surprising to see that the PWDs or the caretakers or the NGOs working for PWDs are not interested to send the PWDs to school that is the reason why only the $6^{th}$ note which corresponds to no proper school education remains in the off state. Further the $10^{th}$ node gets the maximum value of 16 i.e., all the 3 groups of experts feel 'poor economy' of the rural poor PWDs is the root cause of all other problems.

However all agree the PWDs are poorly nourished for it takes the next highest value namely 5. The first node, the seventh, eighth, ninth and the eleventh node take the value 4. It means that all the experts give equal importance to these nodes that is why they all get a value 4 in the hidden pattern. The nodes five and nine viz. no recreation and ignorance of welfare measures of the government, get the lesser value of 3. They are least bothered about shelter, the fourth node that gets least score. Not much importance is given for clothing. The only explanation one can derive from the less value of the nodes five and nine is that they have less chances for recreation and are ignorant of the welfare measures given by the government for they are used well by urban people and by the rich and the middle class. So their ignorance about it has not only made them poor but also unemployed.

Thus no proper health care is attributed to poor economy (poverty) as that node takes the value 16.

Next we study the effect of the on state of the state vector "no employment / self-employment". Let $Y = (0\ 0\ 0\ 0\ 0\ 0\ 1\ 0\ 0\ 0\ 0)$ be the given state vector. Effect of Y on M is given by

$$YM \quad = \quad (0\ 0\ 0\ 0\ 0\ 0\ 1\ 1\ 1\ 2\ 2)$$

after updating

Let

$$YM = Y_1$$
$$= (0\ 0\ 0\ 0\ 0\ 0\ 1\ 1\ 1\ 1\ 1)$$

after thresholding,

$$Y_1M \hookrightarrow (2\ 2\ 2\ 1\ 3\ 0\ 4\ 1\ 1\ 4\ 4)$$
$$= Y_2$$
$$= (1\ 1\ 1\ 1\ 1\ 0\ 1\ 1\ 1\ 1\ 1)$$

after thresholding,

$$Y_2M \hookrightarrow (4\ 5\ 2\ 1\ 3\ 0\ 4\ 4\ 4\ 16\ 4)$$
$$= Y_3 \text{ (say)}.$$

We see the hidden pattern of Y is the same as that X thus according to the combined opinion the nodes no proper health care and no employment / self-employment hold the same status.

Next we study the on state of the node poor economy (poverty) and all other nodes being in the off state. Let

$$S = (0\ 0\ 0\ 0\ 0\ 0\ 0\ 0\ 0\ 1\ 0)$$

be the given state vector. The effect of S on the dynamical system M given by

$$SM \hookrightarrow (2\ 2\ 2\ 1\ 1\ 0\ 0\ 0\ 0\ 1\ 2).$$

After thresholding SM we get

$$S_1 = (1\ 1\ 1\ 1\ 1\ 0\ 0\ 0\ 0\ 1\ 1).$$

Now we study the effect of $S_1$ on M,

$$S_1M \hookrightarrow (4\ 5\ 2\ 1\ 1\ 0\ 2\ 2\ 2\ 2\ 13\ 2)$$

$S_1M$ after thresholding and updating gives $S_2$ where

$$S_2 \quad = \quad (1\ 1\ 1\ 1\ 1\ 0\ 1\ 1\ 1\ 1).$$

We see the resultant of $S_2$ on M is given by

$$S_2M \quad \hookrightarrow \quad (4\ 5\ 2\ 1\ 3\ 0\ 4\ 4\ 3\ 16\ 4).$$

This resultant / hidden pattern is the same as those when the nodes 1 or 7 were in the on state.

We now give the conclusions of this chapter.

**Conclusions**

In this chapter we have introduced a new or a modified form of combined FCMs which we choose to call as special FCM model. We are justified in using this model as:

(1) Our attributes do not give the value –1 in this connection matrix.

(2) This is better than combined FCMs for when we use combined FCM all the entries would lie in the interval 0 to 82 so while working on any state vector it would only complicate the model and take away both time and money. So we defined the new special FCM model to have only 0's and 1's to be entries of the connection matrix there by making the working simple and not as one which is time consuming. Now our main observation from this study is that poor economy (poverty) is the major cause of all the problems faced by the caretakers and the PWDs.

Secondly they feel school education for the PWDs will have no real effect on them (which will be analysed in chapter seven), which is evident from the off state of the node 6, in the resultant vector.

Another factor is that they have not used the government support in any form and they are ignorant of it as seen from the hidden pattern.

Finally the hidden pattern in all the cases resulted in a fixed point there by confirming that nothing has changed their poor economy (poverty) and the related problems. That is no cyclic change in any of the factors analyzed in this problem. This implies that the problems are constant and need to be addressed effectively.

PWDs (their) main interest first lies in employment more so in self-employment; hence steps should be taken to give them good school education followed by job fetching skill training and see to it they are self-employed or employed.

Chapter Five

# THE STUDY OF INTER RELATED PROBLEMS FACED BY THE PWDS AND THE CARETAKERS

This chapter has five sections. Section one is introductory in nature; gives a brief literature survey to the problem and describes the various nodes / attributes given by an heterogeneous group of experts. In section two, the analysis of the problems of the PWDs and their related problems using FRM is carried out by a public person from Melmalayanur block. The same problem is analysed by a PWD using FRMs in section three. Section four gives the expert opinion of caretakers of a PWD using FRM. The same problem is analysed using a set of attributes of his choice, by an NGO working solely for PWDs using the FRM model.

## 5.1 Introduction

In this chapter we study the inter relations between the problems faced by the PWDs and the caretakers. We from our study made in chapter four found that most of the PWDs were depressed and they felt that they (PWDs) did not get proper or good treatment from their caretakers. Most of them felt and said that the other members of the family were treated better in every aspect; food, clothing, recreation and so on. So in this chapter we venture to study the interrelation of the PWDs and the caretakers using FRMs.

Now we give a brief literature survey of the FRMs. In the first place FRMs are a special type of FCMs when the related nodes or attributes can be divided into two disjoint classes. Just like FCMs, FRMs are also powerful dynamical systems which are sued to analyze specially the unsupervised data. Unlike FCMs the relational matrix associated with FRM are in general rectangular (however we can also have square matrices associated with FRMs and these square matrices may not have the diagonal elements to be zero as in case of FCMs). The hidden patter in case of FRMs, occur as a pair of fixed points or as a pair of limit cycle or a limit cycle and a fixed point or a fixed point and a limit cycle. Thus FRMs can give the effect of a node in the domain space on the nodes of the range space and vice versa.

FRMs were constructed in the year 2000 by W.B. Vasantha Kandasamy and Y. Sultana [44,47,48,57]. In 2001 the FRMs were used by them to analyze the employee employer relationship [44]. For more about FRMs please refer [44,47, 48,57]. We use this FRMs to study the problems of PWDs and the factors which contribute to the problem.

At the outset we are justified in using FRMs as we cannot get any numerical data relating to the sufferings or ill-treatment or discriminative treatment the PWDs undergo. Since all the problems related with them are very sensitive, abstract and cannot be measured in numbers we are forced to use this model.

We have taken views from 150 experts who form and heterogeneous group. The group consists of PWDs, caretakers, NGOs associated with PWDs, socio scientists, politicians, educationists and the public from Melmalayanur and Kurinji Padi blocks of Tamil Nadu. Out of the 150 opinions sought by us 7 of them turned to be invalid for analysis as they had failed to follow the instructions given by us.

We have displaced only views from four experts however for our study and conclusions we have made use of all the experts opinion. We just recall that the working of the FRMs are described in [44,47,48,57]. We now proceed onto describe the domain and range attributes of the PWDs and their problems.

The attributes related with the PWDs are taken as the nodes of the domain space of the FRMs. We describe the attributes associated with the PWDs in a line or two.

1. From the survey we found most of the PWDs looked and they have also acknowledged that they were depressed due to the disability and the treatment they get because of the disability.

2. The PWDs suffer from inferiority complex. Most of the PWDs suffered from inferiority complex due to the criticism they face due to their disability.

3. They suffer from mental agony due to isolation and ill treatment. They suffer mental tensions / trauma.

4. Self image - many did not have self image. In other words they lack self esteem.

5. Happy and contended

6. Uninterested in life.

7. Wished for self-employment and wanted to earn.

8. Dependency on others made them frustrated.

9. Ill-treated by relatives and at times by strangers.

10. Longed for recreation and attention from the caretakers.

11. Longing for good food and good cloths.

12. They had more time left to mourn and ponder over the disability and their bleak future and felt insecure.

We will label the 12 attributes associated with disability by $D_1, D_2, \ldots, D_{12}$ where $D_i$ corresponds to the $i^{th}$ attribute; $1 \leq i \leq 12$. The attributes related with the range space are as follows:

$R_1$ - Poverty is the main cause for neglect of PWDs as they (caretakers) have no money to spend on the basic needs like medical aid, food etc of the PWDs.

$R_2$ - Indifferent to PWDs, as they are burden to the caretakers.

$R_3$ - Ashamed to having a PWD, as the society looks low upon them and think or speak of the family as a cursed one by God.

$R_4$ - Sympathetic to the PWDs.

$R_5$ - They are caring for PWDs.

$R_6$ - Show hatred towards PWDs.

$R_7$ - Fatalism - so no proper care is taken.

$R_8$ - Caretakers are not interested in marrying them (PWDs) off as they think marriage may become an additional burden to them.

$R_9$ - The PWDs are an economic burden on them.

$R_{10}$ - The PWDs are isolated from others for reasons best known to them.

$R_{11}$ - Ignorance of caretakers about government welfare measures.

$R_{12}$ - Do not plan to train the PWDs to be self-sufficient.

$R_{13}$ - The caretakers fail to give the PWDs vocational training, which may lead to their self - employment.

$R_{14}$ - Not interested to even sending them to school.

$R_{15}$ - Negative attribute of caretakers pains the PWDs.

## 5.2 Expert Opinion from a Public Person Living in Melmalayanur Block

The following attributes are given by an expert who is public person from Melmalayanur Block. The problems of PWDs are taken as the nodes of the domain space and the attributes with the close caretakers are taken as the nodes of the range space.

The attributes associated with the PWDs are given below. They are in certain cases described in line or two.

$D_1$ - Depressed. Majority of the PWDs looked and said they were depressed because of their disability and general treatment they get from the caretakers and society in general.

$D_2$ - Suffer from inferiority complex.

$D_3$ - Mental stress/agony. They often were isolated and sometimes kept in a small hut outside the house which made them feel sad as well as gave time to think about their disability with no proper work. So they were often under stress and mental tension.

$D_4$ - Self image

$D_5$ - Happy and contended.

$D_6$ - Uninterested in life.

$D_7$ - Dependent on others for every thing.

$D_8$ - Lack of mobility

$D_9$ - Ill-treated by close relatives.

The attributes $D_1$, $D_2$, …, $D_9$ are taken as the nodes of the domain space of the FRM. We give attributes associated with the range space:

$R_1$ - Poor so cannot find money to spend on basic requirements. The PWDs go to work for their livelihood.

$R_2$ - Ashamed - relatives were ashamed of the PWDs.

$R_3$ - Indifferent-they were treated indifferently by their caretakers.

$R_4$ - PWDs are a burden to them. So they neglected them totally.

$R_5$ - Fatalism - they said it was fate that they have a PWD as their child / relative.

$R_6$ - Sympathetic.

$R_7$ - Caring.

$R_8$ - Show hatred towards the PWDs.

$R_9$ - The caretakers were not interested in marrying them off.

$R_{10}$ - The PWDs are an economic burden on them.

$R_{11}$ - They were isolated from others for reasons best known to the caretakers.

Thus $R_1, R_2, \ldots, R_{11}$ are taken as the nodes of the range space of the FRM.

The directed graph related to the view of the expert is given in following figure:

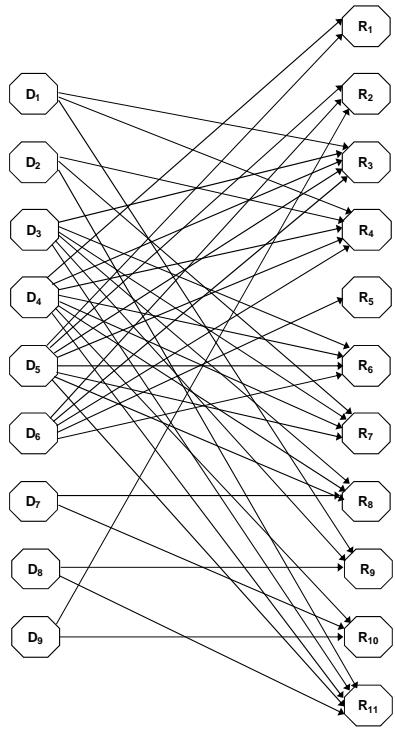

Let the relation matrix associated with the directed graph given above be given by T. T is a $9 \times 11$ matrix with entries from the set $\{0, -1, 1\}$.

$$T = \begin{array}{c c} & \begin{array}{c c c c c c c c c c c} R_1 & R_2 & R_3 & R_4 & R_5 & R_6 & R_7 & R_8 & R_9 & R_{10} & R_{11} \end{array} \\ \begin{array}{c} D_1 \\ D_2 \\ D_3 \\ D_4 \\ D_5 \\ D_6 \\ D_7 \\ D_8 \\ D_9 \end{array} & \left[ \begin{array}{c c c c c c c c c c c} 0 & 0 & 1 & 1 & 0 & 0 & 0 & 0 & 1 & 0 & 0 \\ 0 & 0 & 0 & 1 & 0 & 0 & -1 & 0 & 0 & 0 & 1 \\ 0 & 0 & 1 & 0 & 0 & -1 & -1 & 1 & 1 & 0 & 1 \\ -1 & 0 & -1 & -1 & 0 & 1 & 1 & -1 & 0 & -1 & -1 \\ -1 & -1 & -1 & -1 & 0 & 1 & 1 & -1 & 0 & 0 & -1 \\ 0 & 1 & 1 & 1 & 1 & -1 & 0 & 0 & 0 & 0 & 0 \\ 0 & 0 & 0 & 0 & 0 & 0 & 0 & 1 & 0 & 1 & 0 \\ 0 & 0 & 0 & 1 & 0 & 0 & 0 & 0 & 1 & 0 & 1 \\ 0 & 1 & 0 & 0 & 0 & 0 & 0 & 0 & 0 & 1 & 0 \end{array} \right] \end{array}$$

Now we study the effect of the state vectors on the dynamical system T. Suppose the expert wishes to study the 'on state' of the node $D_1$ and all other nodes are in the 'off state'. Let the state vector be $X = (1\ 0\ 0\ 0\ 0\ 0\ 0\ 0\ 0)$. The effect of X on the dynamical system T is given by

$$\begin{aligned} XT \quad &= \quad (0\ 0\ 1\ 1\ 0\ 0\ 0\ 0\ 1\ 0\ 0) \\ &= \quad Y\ (\text{say}). \end{aligned}$$

$$\begin{aligned} YT^t \quad &= \quad (3\ 1\ 2\ -2\ -2\ 2\ 0\ 2\ 0). \\ &\hookrightarrow \quad (1\ 1\ 1\ 0\ 0\ 1\ 0\ 1\ 0) \\ &= \quad X_1\ (\text{say}). \end{aligned}$$

($\hookrightarrow$ denotes that the resultant state vector $YT^t$ is updated and thresholded i.e., all negative values and 0 are replaced by 0 and all positive values greater than or equal to one are replaced by 1. By updating we mean the coordinate which we started in the 'on state' should remain in the 'on state' till the end).

Now we find

$$\begin{aligned} X_1 T \quad &\hookrightarrow \quad (0\ 1\ 1\ 1\ 1\ 0\ 0\ 1\ 1\ 0\ 1) \\ &= \quad Y_1\ (\text{say}). \end{aligned}$$

$$Y_1 T^t \quad \hookrightarrow \quad (1\ 1\ 1\ 0\ 0\ 1\ 1\ 1\ 1)$$
$$= \quad X_2 \text{ (say)}.$$

$$X_2 T \quad \hookrightarrow \quad (0\ 1\ 1\ 1\ 1\ 0\ 0\ 1\ 1\ 1\ 1)$$
$$= \quad Y_2 \text{ (say)}.$$

$$Y_2 T^t \quad \hookrightarrow \quad (1\ 1\ 1\ 0\ 0\ 1\ 1\ 1\ 1)$$
$$= \quad X_3 \text{ (say)}.$$
$$= \quad X_2.$$

Thus the hidden pattern gives a fixed pair given by {(1 1 1 0 0 1 1 1 1), (0 1 1 1 1 0 0 1 1 1 1)}. Thus when the node depressed alone in the domain space is in the 'on state' we see this makes the nodes $D_2$, $D_3$, $D_6$, $D_7$, $D_8$, $D_9$ to come to 'on state' in the domain space and $R_2$, $R_3$, $R_4$, $R_5$, $R_8$, $R_9$, $R_{10}$ and $R_{11}$ come to 'on state' in the range space.

Thus we see except the nodes the PWD has self image and she / he is happy and contended all other nodes come to 'on state'. Thus this reveals if a PWD is depressed certainly he has no self image and he is not happy and contended. Further it also reveals from the state vector in the domain space poverty is not a cause of depression for $R_1$ is in the 'off state'. Also $R_6$ and $R_7$ alone do not come to 'on state' which clearly shows that the caretakers are not sympathetic and caring which is one of the reasons for the PWDs to be depressed. Thus we see all negative attributes come to 'on state' in both the spaces when the PWD is depressed.

Next the expert is interested in studying the effect of the 'on state' of the node in the range space, viz. $R_6$ i.e., the caretakers are sympathetic towards the PWDs.

Let $Y = (0\ 0\ 0\ 0\ 0\ 1\ 0\ 0\ 0\ 0\ 0)$ be the state vector of the range space. To study the effect of Y on the dynamical system $T^t$.

$$YT^t \quad = \quad (0\ 0\ 0\ 1\ 1\ 0\ 0\ 0\ 0)$$
$$= \quad X_1 \text{ (say)}.$$

$$X_1T \quad = \quad (0\ 0\ 0\ 0\ 0\ 1\ 1\ 0\ 0\ 0\ 0)$$
$$= \quad Y_1 \text{ (say)}.$$

$$Y_1T^t \quad \hookrightarrow \quad (0\ 0\ 0\ 1\ 1\ 0\ 0\ 0\ 0)$$
$$= \quad X_2 \text{ (say)}.$$

But $X_2 = X_1$. Thus we see the hidden pattern of the state vector is a fixed pair of points given by {$(0\ 0\ 0\ 0\ 0\ 1\ 1\ 0\ 0\ 0\ 0)$, $(0\ 0\ 0\ 1\ 1\ 0\ 0\ 0\ 0)$}. It is clear when the PWD is treated with sympathy it makes him feel their caretakers are caring. So $R_7$ come to 'on state'. On the other hand, we see she/he is happy and contended with a self image.

Next the expert wishes to find the hidden pattern of the 'on state' of the domain node $D_4$, i.e., self image of the PWD alone is in the 'on state'.

Let $P = (0\ 0\ 0\ 1\ 0\ 0\ 0\ 0\ 0)$ be the given state vector. The effect of $P$ on $T$ is given by

$$PT \quad \hookrightarrow \quad (0\ 0\ 0\ 0\ 0\ 1\ 1\ 0\ 0\ 0\ 0)$$
$$= \quad S_1 \text{ (say)}$$

$$S_1 T^t \quad \hookrightarrow \quad (0\ 0\ 0\ 1\ 1\ 0\ 0\ 0\ 0)$$
$$= \quad P_1 \text{ (say)}.$$

$$P_1 T \quad \hookrightarrow \quad (0\ 0\ 0\ 0\ 0\ 1\ 1\ 0\ 0\ 0\ 0)$$
$$= \quad S_2 \text{ (say)}.$$

But $S_2 = S_1$ resulting in a fixed pair. Thus the hidden pattern of $P$ is a fixed pair. We see self image of the PWD makes him happy and contended. He / she also feel that the caretakers are caring and sympathetic towards them.

Now the expert studies the effect of the state vector in the range space when the PWD is isolated from the other, i.e., when $R_{11}$ is in the 'on state'. Let X = (0 0 0 0 0 0 0 0 0 0 1) be the given state vector. Its effect on the dynamical system T is given by

$$XT^t \hookrightarrow (0\ 1\ 1\ 0\ 0\ 0\ 0\ 1\ 0)$$
$$= Y \text{ (say)}.$$

$$YT \hookrightarrow (0\ 0\ 1\ 1\ 0\ 0\ 0\ 1\ 1\ 0\ 1)$$
$$= X_1 \text{ (say)}.$$

The effect of $X_1$ on T is given by

$$X_1T^t \hookrightarrow (1\ 1\ 1\ 0\ 0\ 1\ 1\ 1\ 1)$$
$$= Y_1 \text{ (say)}.$$

$$Y_1T \hookrightarrow (0\ 1\ 1\ 1\ 1\ 0\ 0\ 1\ 1\ 1\ 1)$$
$$= X_2 \text{ (say)}.$$

$$X_2T^t \hookrightarrow (1\ 1\ 1\ 0\ 0\ 1\ 1\ 1\ 1)$$
$$= Y_2 \text{ (say)}.$$

We see $Y_2 = Y_1$. Thus the hidden patter of the state vector is a fixed pair given by {(0 1 1 1 1 0 0 1 1 1 1), (1 1 1 0 0 1 1 1 1)}. Thus when the PWD is isolated from others he / she suffers all negative attributes and it is not economic condition that matters. Isolation directly means they are not taken care of and the caretakers are not sympathetic towards them. When they are isolated they are not happy and contended and they do not have self-image. All this is evident from the hidden pattern in which $R_1$, $R_6$ and $R_7$ are 0 and $D_4$ and $D_5$ are 0, i.e., in the 'off state'.

We have worked with the several other state vectors and the conclusions are based on that, as well as from the survey we have taken. This is given in the chapter eight of this book.

## 5.3 Expert Opinion of a PWD using FRMs

The attributes taken by him are enumerated as follows.

$D_1$ - Suffers from inferiority complex.

$D_2$ - No self-image.

$D_3$ - Discriminated from other family members.

$D_4$ - No good food or clothing even on festivals.

$D_5$ - Frustrated because they have to depend on caretakers for everything.

$D_6$ - Wishes to be trained in some trade so that they can be self employed.

The range attributes taken by him are as follows:

$R_1$ - Poverty

$R_2$ - The caretakers are depressed and dejected for having the PWDs as their children / wards

$R_3$ - No mind to spend on PWDs

$R_4$ - Ill-treated as PWDs are an economic burden and a social disgrace.

$R_5$ - Never takes any steps to make them (PWDs) get any form of vocational training.

This particular expert who is a PWD did not want to consider the attributes given by the others. He wishes to work only with these. As he happened to be open and also understanding the role of playing as an expert we opted to get his opinion.

The directed graph given by him is as follows.

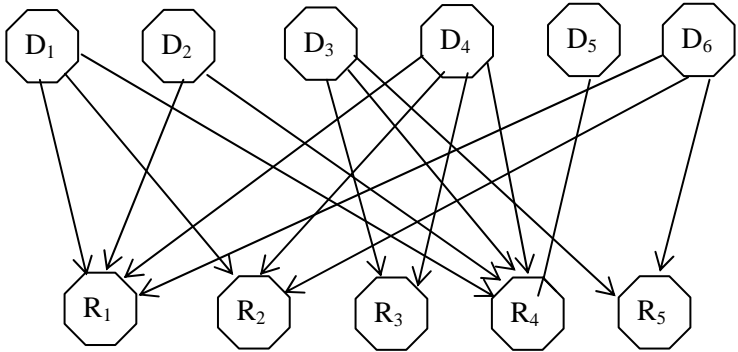

The related relational matrix of the FRM is given by M.

$$M = \begin{array}{c} \\ D_1 \\ D_2 \\ D_3 \\ D_4 \\ D_5 \\ D_6 \end{array} \begin{array}{c} R_1 \ \ R_2 \ \ \ R_3 \ R_4 \ R_5 \\ \begin{bmatrix} 1 & 1 & 0 & 1 & 0 \\ 1 & 0 & 0 & 1 & 0 \\ 0 & 0 & 1 & 1 & 1 \\ 1 & 1 & 1 & 1 & 0 \\ 0 & 0 & 0 & 1 & 0 \\ -1 & -1 & 0 & 0 & 1 \end{bmatrix} \end{array}$$

Now suppose the expert wishes to study the on state $D_1$ alone and all other nodes remain in the off state. Let $X = (1\ 0\ 0\ 0\ 0\ 0)$ be the given state vector. To find the effect of X on the dynamical system M.

$$\begin{aligned} XM \quad &\hookrightarrow \ (1\ 1\ 0\ 1\ 0) \\ &= \ Y \text{ (say)}. \end{aligned}$$

$$\begin{aligned} YM^t \quad &\hookrightarrow \ (1\ 1\ 1\ 1\ 1\ 0) \\ &= \ X_1 \text{ (say)}. \end{aligned}$$

$$X_1M \quad \hookrightarrow \quad (1\ 1\ 1\ 1\ 1)$$
$$= \quad Y_1 \ \text{(say)}.$$

$$Y_1M \quad \hookrightarrow \quad (1\ 1\ 1\ 1\ 1\ 0)$$
$$= \quad X_2 \ \text{(say)}.$$

We see $X_1 = X_2$. Thus the hidden pattern of the given state vector X is a fixed point given by the pair $\{(1\ 1\ 1\ 1\ 1\ 0), (1\ 1\ 1\ 1\ 1)\}$. Thus when the PWDs are suffering from the inferiority complex all nodes both in the domain space and the range space come to on state. Only the node $D_6$ is in the off state which makes one understand that if they are suffering from inferiority complex they are not employed or they don't earn and they are unaware of any job suiting their disability.

Now the expert wishes to study the on state of the node $R_5$ alone in the on state of the range space and all other nodes are in the off state. To find the effect of the state vector Y on the dynamically system $M^t$. Let

$$Y \quad = \quad (0\ 0\ 0\ 0\ 1)$$

$$YM^t \quad \hookrightarrow \quad (0\ 0\ 1\ 0\ 0\ 0)$$
$$= \quad X \ \text{(say)}.$$

$$XM \quad \hookrightarrow \quad (0\ 0\ 1\ 1\ 1)$$
$$= \quad Y_1 \ \text{(say)}.$$

$$Y_1M^t \quad \hookrightarrow \quad (1\ 1\ 1\ 1\ 1\ 1)$$
$$= \quad X_1 \ \text{(say)}.$$

$$X_1M \quad \hookrightarrow \quad (1\ 1\ 1\ 1\ 1)$$
$$= \quad Y_2 \ \text{(say)}.$$

$$Y_2M^t \quad \hookrightarrow \quad (1\ 1\ 1\ 1\ 1\ 1)$$
$$= \quad X_2 \ (= X_1).$$

Thus the hidden pattern reduces to a fixed pair given by [(1 1 1 1 1 1), (1 1 1 1 1)}. Therefore it is suggested that failure to give any form of vocational training to the PWDs by the caretakers leads to all negative attributes on the part of the PWDs and also they remain unemployed.

Thus the following conclusions are drawn.

(1)     All attributes from the domain space and range space come to on state when $R_5$ is in the on state.

(2)     It is suggested that the PWDs want to get employment more so self employment.

(3)     From our study and discussions it is learnt that the PWDs are very sensitive, when someone (say small boys or girls) laugh at them when they try to walk or move with their disability.  They become dejected, tensed and become uncomfortable.  It is suggested that counseling may be given to the PWDs to remain stable or unaffected by the gesture of the passers by or their own relatives.

(4)     PWDs crave for equal treatment and long for attention from the caretakers.

## 5.4   Expert Opinion of Caretakers of the PWDs using FRM

Here we give the expert opinion of the caretaker of the PWD.  The caretakers selected a set of attributes, which was much different from the PWDs and the public.  We first list the attributes chosen by the caretakers for the domain and range space.

The attributes given by the caretakers for the domain space:

$D_1$          -     Well taken care of by caretakers

$D_2$          -     Depressed / no self-image.

$D_3$     -     No training for future / marrying them off.

$D_4$     -     Dependent on every act.

$D_5$     -     Under nourished

$D_6$     -     Suffer from health problems.

$D_7$     -     Not property maintained physically by the caretakers.

The following are the nodes given by the caretakers for the range space.

$R_1$     -     Poverty - so cannot spend on PWD as much as they spend on the breadwinners.

$R_2$     -     Sufficient time cannot be spent on the PWDs as they are busy with their day-to-day earnings as most of them are daily wages.

$R_3$     -     No time or means to think of sending or even training the PWDs due to lack of time and money.

$R_4$     -     Social stigma as PWDs born in the family. There is a general belief that the family is cursed by God.

$R_5$     -     Fear of marrying PWDs as they may give birth to PWDs.

$R_6$     -     Not aware of any NGOs' working for PWDs or government aid given to PWDs or do not want to openly acknowledge that a PWD is living in their house.

Using these attributes the directed graph given by the expert is as follows:

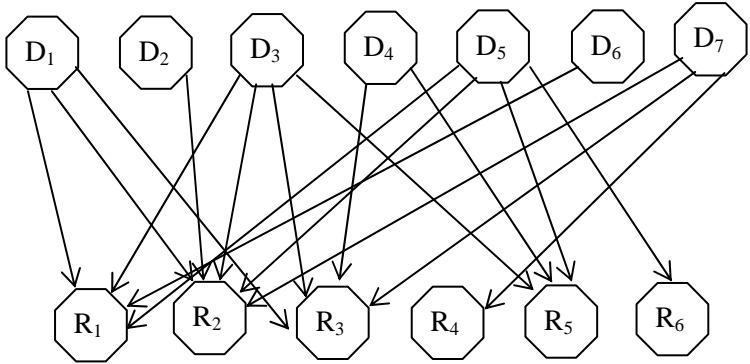

The relational matrix P obtained from the directed graph as follows.

$$P = \begin{array}{c} \\ D_1 \\ D_2 \\ D_3 \\ D_4 \\ D_5 \\ D_6 \\ D_7 \end{array} \begin{array}{c} R_1 \ \ R_2 \ \ R_3 \ R_4 \ R_5 \ R_6 \\ \begin{bmatrix} -1 & -1 & -1 & 0 & 0 & 0 \\ 0 & 1 & 0 & 0 & 0 & 0 \\ 1 & 1 & 1 & 0 & 1 & 0 \\ 0 & 0 & 1 & 0 & 1 & 0 \\ 1 & 1 & 0 & 0 & 1 & 1 \\ 1 & 0 & 0 & 0 & 0 & 0 \\ 0 & 1 & 1 & 1 & 0 & 0 \end{bmatrix} \end{array}$$

Now we see the effect of the state vector when the node $D_1$, i.e., the PWDs are well taken care of by the caretakers alone is in the on state and all other nodes are in the off state.

Let

$$X \quad = \quad (1\ 0\ 0\ 0\ 0\ 0\ 0)$$

$$XP \quad = \quad (-1\ -1\ -1\ 0\ 0\ 0);$$

After thresholding XP we get

$$Y \quad = \quad (0\ 0\ 0\ 0\ 0\ 0).$$

Now

$$YP^t \quad = \quad (0\ 0\ 0\ 0\ 0\ 0\ 0);$$

After updating $YP^t$ we get

$$X_1 \quad = \quad (1\ 0\ 0\ 0\ 0\ 0\ 0)$$
$$= \quad X.$$

Thus we see this state vector X gives a hidden pattern to be just $(X, (0\ 0\ 0\ 0\ 0\ 0))$. Thus this node of PWDs are well taken care of does not influence the other nodes of both the domain and range space. Now we study the on state of the node $D_2$ alone in the on state.

Let

$$T \quad = \quad (0\ 1\ 0\ 0\ 0\ 0\ 0)$$

be the given state vector.
To find the effect of T on P.

$$TP \quad = \quad (0\ 1\ 0\ 0\ 0\ 0)$$
$$= \quad S \text{ (say)}.$$

Now

$$SP^T \quad = \quad (-1\ 1\ 1\ 0\ 1\ 0\ 1);$$

after updating and thresholding $SP^T$, we get

$$T_1 \quad = \quad (0\ 1\ 1\ 0\ 1\ 0\ 1) \text{ (say)},$$
$$T_1P \quad \hookrightarrow \quad (1\ 1\ 1\ 1\ 1\ 1)$$
$$= \quad S_1 \text{ (say)}.$$

Now

$$S_1P \quad \hookrightarrow \quad (0\ 1\ 1\ 1\ 1\ 1\ 1)$$
$$= \quad T_2 \text{ (say)}$$

$$T_2 P^T \hookrightarrow (1\ 1\ 1\ 1\ 1\ 1)$$
$$= S_2 \text{ (say)}.$$

But $S_1 = S_2$ thus the hidden pattern of the state vector T is a fixed pair given by $\{(0\ 1\ 1\ 1\ 1\ 1\ 1),\ (1\ 1\ 1\ 1\ 1\ 1\ 1)\}$. We see when the PWDs are depressed all nodes from the range space and domain space come to on state except the node $D_1$ of the domain space; i.e., the PWD are well taken care of by the caretakers.

Now we study the on state of the node poverty from the range space. Let

$$Y = (1\ 0\ 0\ 0\ 0\ 0).$$

Now the effect of Y on the dynamical system P is given by

$$YP^T \hookrightarrow (0\ 0\ 1\ 0\ 1\ 1\ 0)$$
$$= X \text{ (say)}.$$

$$XP \hookrightarrow (1\ 1\ 1\ 0\ 1\ 1)$$
$$= Y_1 \text{ (say)}$$

$$Y_1 P^T \hookrightarrow (0\ 1\ 1\ 1\ 1\ 1\ 1)$$
$$= X_1 \text{ (say)}.$$

$$X_1 P \hookrightarrow (1\ 1\ 1\ 1\ 1\ 1)$$
$$= Y_2 \text{ (say)}$$

$$Y_2 P^t \hookrightarrow (0\ 1\ 1\ 1\ 1\ 1\ 1)$$
$$= X_2 \text{ (say)}.$$

Thus $X_1 = X_2$ leading to a fixed point. We see when poverty is in the on state all the nodes from both the domain and range space come to on state except the node $D_1$, i.e., the PWDs are well taken care of by the caretakers. Thus the on state of the node the PWDs are in the depressed state of the domain space or the on state of the node $R_1$, the poverty of the range space given the same hidden pattern. Thus we see poverty of the rural PWD really cripples both the PWDs and the caretakers. Let us now study the on state of the node $R_4$ in the range space; i.e., social stigma for having a PWD born their family. Let X = (0 0 0 1 0 0) be given state vector. To study the effect of X on the dynamical system P.

$$XP^t \hookrightarrow (0\ 0\ 0\ 0\ 0\ 0\ 1)$$
$$= \text{Y (say)}$$

$$YP \hookrightarrow (0\ 1\ 1\ 1\ 0\ 0)$$
$$= X_1 \text{ (say)}$$

$$X_1P^t \hookrightarrow (0\ 1\ 1\ 1\ 1\ 0\ 1)$$
$$= Y_1 \text{ (say)}$$

$$Y_1P \hookrightarrow (1\ 1\ 1\ 1\ 1\ 1)$$
$$= X_2 \text{ (say)}.$$

$$X_2P^t \hookrightarrow (0\ 1\ 1\ 1\ 1\ 1\ 1)$$
$$= Y_2 \text{ (say)}$$

$$Y_2P \hookrightarrow (1\ 1\ 1\ 1\ 1\ 1)$$
$$= X_3 \text{ (say)}.$$

Since $X_3 = X_2$ we see $X_3P^t = Y_3$ (=$Y_2$). Thus the hidden pattern of the state vector with $R_4$ in the on state in the range space is a fixed point given by the pair {(0 1 1 1 1 1 1), (1 1 1 1 1 1)}.

We see the FRM model given by the caretakers behaves in a unique way. For every hidden pattern is a fixed point. One of the nodes, that is the node $D_1$ in the domain space in on state makes all the nodes in the domain and range space to be in the off state. Another special feature of this system given by the caretakers is that whenever any node is taken in the on state either from the domain space or from the range space barring $D_1$ the hidden pattern is a fixed pair given by the pair {(0 1 1 1 1 1 1), (1 1 1 1 1 1)}.

This mainly shows the nodes given by the caretakers are so interdependent on each other so strongly that the on state of any one of them makes on the rest of the nodes from both the domain and range space.

In the following section we give the opinion of a NGO who solely works for the PWDs.

## 5.5 Expert Opinion of the NGO working for the PWDs using the FRM model

This NGO is an organization that has been working for the PWDs for more than a decade. The representative of the NGO has in the first place listed the following attributes associated with PWDs, which has been taken as the nodes of the domain space D.

$D_1$ - PWDs are aware of government special assistance to them.

$D_2$ - PWDs suffer chronic depression due to disability.

$D_3$ - PWDs are ill-treated by the caretakers as they feel PWDs are curse and social stigma to their family.

$D_4$ - PWDs are poorly maintained by the caretakers (Dress, Cleanliness etc).

$D_5$ - PWDs are not interested in school education, as it is painful to be humiliated by their own age group children or younger to them.

$D_6$ - PWDs aspire (wish) to be employed / self-employed or earn money; to be more exact want to be economically independent.

Now we given the nodes related with the caretakers and the society which forms the attributes of the range space.

$R_1$ - The caretakers do not take any effort to seek government support for the PWDs.

$R_2$ - Poverty of the caretakers is the root cause of all problems.

$R_3$ - Caretakers suffer humiliation in the society for having a PWD as their child.

$R_4$ - Caretakers do not have the mindset to invest on PWDs as they think are not going to repay them.

$R_5$ - Caretakers fail to consult NGO's working for the PWDs due to the reasons well known to them.

We now give the directed graph related with the FRM using the nodes $D_1$, ..., $D_6$ for the domain space and $R_1$, ..., $R_5$ as the attributes for the range space.

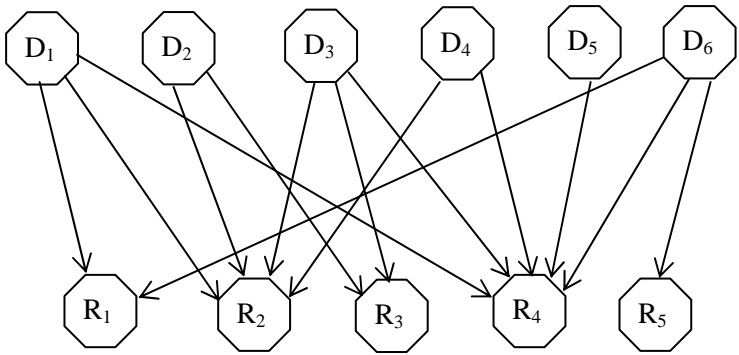

Let Z be the related connection matrix of the above directed graph given in figure.

$$M = \begin{array}{c} \\ D_1 \\ D_2 \\ D_3 \\ D_4 \\ D_5 \\ D_6 \end{array} \begin{array}{ccccc} R_1 & R_2 & R_3 & R_4 & R_5 \\ \left[ \begin{array}{ccccc} 1 & 1 & 1 & 1 & 0 \\ 0 & 1 & 1 & 0 & 0 \\ 0 & 1 & 1 & 1 & 0 \\ 0 & 1 & 0 & 1 & 0 \\ 0 & 0 & 0 & 1 & 0 \\ -1 & 1 & 0 & 0 & 1 \end{array} \right] \end{array}$$

Let us now consider the effect of the state vector in which the node $D_1$ alone in the on state and all other nodes remain in the off state. To find the effect of the state vector $X = (1\ 0\ 0\ 0\ 0\ 0)$ on the dynamical system Z.

$$\begin{array}{lll} XZ & \hookrightarrow & (1\ 1\ 1\ 1\ 0) \\ & = & Y \text{ (say)} \end{array}$$

$$\begin{array}{lll} YZ^t & \hookrightarrow & (1\ 1\ 1\ 1\ 1\ 0) \\ & = & X_1 \text{ (say)}. \end{array}$$

Now $\quad \begin{array}{lll} X_1M & \hookrightarrow & (1\ 1\ 1\ 1\ 0) \\ & = & Y_1 \text{ (say)}. \end{array}$

Since $Y_1 = Y$ we see the hidden pattern of the state vector X is a fixed pair given by $\{(1\ 1\ 1\ 1\ 0), (1\ 1\ 1\ 1\ 1\ 0)\}$ i.e., we see the all the nodes both in the domain space and the range space come to on state except the node $D_6$ and $R_5$ from the domain space the range space respectively. Thus we see when the PWDs are unaware of the government assistance they aspire to get self-employed remains in the off state in the domain space as well as the node; caretakers fail to consult the NGO's working for the PWDs are not sought for the reasons well known to the caretakers.

Thus we now study the effect of the node $R_3$ of the range space alone to be in the on state, i.e., the node caretakers suffer humiliation in the society for having a PWD as their child alone is in the on state and all other nodes are in the off state.

Let $Y = (0\ 0\ 1\ 0\ 0)$ to be given state vector. To find the effect of $Y$ on the dynamical system $Z$.

$$Y \quad = \quad (0\ 0\ 1\ 0\ 0) \quad \text{and}$$

$$Y\,Z^t \quad \hookrightarrow \quad (1\ 1\ 1\ 0\ 0\ 0)$$
$$= \quad X \text{ (say)}.$$

$$XZ \quad \hookrightarrow \quad (1\ 1\ 1\ 1\ 0)$$
$$= \quad Y_1 \text{ (say)}$$

$$Y_1\,Z^t \quad \hookrightarrow \quad (1\ 1\ 1\ 1\ 1\ 0)$$
$$= \quad X_1 \text{ (say)}$$

$$X_1Z \quad \hookrightarrow \quad (1\ 1\ 1\ 1\ 0)$$
$$= \quad Y_2 \text{ (say)}.$$

Thus $Y_1 = Y_2$, yields the hidden pattern to be fixed point. Only the nodes $D_6$ and $R_5$ alone are in the off state and all other nodes come to on state. Another significance is that we see the on state of the node in the domain space viz. PWDs are unaware of government special assistance to them and the node $R_3$ of the range space viz., caretakers suffer humiliation in the society for having a PWD as their child give the same hidden pattern.

Now we just study the on state of the node $D_4$ alone and all other nodes of the domain space remain in the off state. Let $T = (0\ 0\ 0\ 1\ 0\ 0)$ be the given state vector. The effect of $T$ on the dynamical system $Z$ is as follows:

$$TZ \quad \hookrightarrow \quad (0\ 1\ 0\ 1\ 0)$$
$$= \quad S$$

$$SZ^t \quad \hookrightarrow \quad (1\ 1\ 1\ 1\ 1\ 1)$$
$$= \quad T_1 \text{ (say).}$$

$$T_1Z \quad \hookrightarrow \quad (0\ 1\ 1\ 1\ 1)$$
$$= \quad S_1 \text{ (say).}$$

$$S_1Z^t \quad \hookrightarrow \quad (1\ 1\ 1\ 1\ 1\ 1)$$
$$= \quad T_2 \text{ (say).}$$

Thus $T_1 = T_2$, yields an hidden pattern which is a fixed point. However all the nodes come to on state except the node $R_1$, the caretakers don't take nay effect to seek government support for the PWDs.

Next we study the on state of $R_4$ alone from the range space. Let $W = (0\ 0\ 0\ 1\ 0)$ be the given state vector. Effect of W on Z is given by

$$WZ^t \quad \hookrightarrow \quad (1\ 0\ 1\ 1\ 1\ 0)$$
$$= \quad V \text{ (say).}$$

$$VZ \quad \hookrightarrow \quad (1\ 1\ 1\ 1\ 0)$$
$$= \quad W_1 \text{ (say).}$$

$$W_1Z^t \quad \hookrightarrow \quad (1\ 1\ 1\ 1\ 1\ 0)$$
$$= \quad V_1 \text{ (say)}$$

$$V_1Z \quad \hookrightarrow \quad (1\ 1\ 1\ 1\ 0)$$
$$= \quad W_2 \text{ (say).}$$

$$W_2Z^t \quad \hookrightarrow \quad (1\ 1\ 1\ 1\ 1\ 0)$$
$$= \quad V_2 \text{ (say).}$$

Thus again we get the hidden pattern to be fixed point.

We see the hidden pattern being a fixed point clearly describes that the problems related with the PWDs and the caretakers do not change with time or cyclic. Thus unless some reformative steps are taken to help the PWDs to get employed by the government and the NGOs not depending on the caretakers then alone the status of both the PWDs and the caretakers would change. More suggestions and conclusions are given in the final chapter of this book.

**Chapter Six**

# ANALYSIS OF PROBLEMS OF RURAL PWDS USING THE NEW FCRM BIMODEL

In this chapter we for the first time use the new FCRM bimodel to analyse the problems of poor rural PWDs. This new model is introduced in chapter three of this book. This model can function simultaneously as FCMs and FRMs. Thus when a node appears both in the FCM and FRM one can find the bihidden pattern and study them.

Further one can carry out the study of a PWD and caretaker simultaneously on a common node. Such a study is impossible with the existing model. Thus this new model is capable of such comparative study. Hence we are justified in using this new model for our study. For the description and working of the new FCRM bimodel one can refer chapter three of this book. Now we proceed onto see the simultaneous opinion of the caretakers and the PWDs.

## 6.1 Expert opinion of the caretaker and a PWD using the new FCRM bimodel

The nodes of the caretakers using special FCM model given in the chapter four of this thesis is taken as the FCM component of the new FCRM bimodel. The FRM model of a PWD given in chapter five of this book is taken as the FRM component of the new FCRM bimodel.

Let $B = B_1 \cup B_2$ be the bimatrix associated with the new FCRM bimodel. Here $B_1 = T$ given in page 60 of chapter four and $B_2 = M$ given in page 85 of chapter five.

$$B = B_1 \cup B_2$$

$$= \begin{bmatrix} 0 & 1 & 0 & 0 & 0 & 0 & 0 & 0 & 0 & 1 & 0 \\ 1 & 0 & 0 & 0 & 0 & 0 & 0 & 0 & 0 & 1 & 0 \\ 0 & 0 & 0 & 0 & 0 & 0 & 0 & 0 & 0 & 1 & 0 \\ 0 & 0 & 0 & 0 & 0 & 0 & 0 & 0 & 0 & 1 & 0 \\ 0 & 0 & 0 & 0 & 0 & 0 & 0 & 0 & 0 & 1 & 0 \\ 0 & 0 & 0 & 0 & 0 & 0 & 0 & 0 & 0 & 0 & 0 \\ 0 & 0 & 0 & 0 & 0 & 0 & 0 & 0 & 0 & 1 & 1 \\ 0 & 0 & 0 & 0 & 0 & 1 & 0 & 0 & 0 & 0 & 0 \\ 1 & 1 & 1 & 1 & 1 & 0 & 0 & 0 & 0 & 0 & 1 \\ 0 & 0 & 0 & 0 & 0 & 0 & 0 & 0 & 0 & 1 & 0 \end{bmatrix} \cup$$

$$\begin{bmatrix} 1 & 1 & 0 & 1 & 0 \\ 1 & 0 & 0 & 1 & 0 \\ 0 & 0 & 1 & 1 & 1 \\ 1 & 1 & 1 & 1 & 0 \\ 0 & 0 & 0 & 1 & 0 \\ -1 & -1 & 0 & 0 & 1 \end{bmatrix}$$

We see the node 10 'poor economy' of the FCM component be in the on state and the node $R_1$ (poverty) of the range space of the FRM component be in the on state.

Now we study the effect of the bistate vector $X = X_1 \cup X_2$ on the dynamical bisystem B.

$$
\begin{aligned}
X &= X_1 \cup X_2 \\
&= (0\ 0\ 0\ 0\ 0\ 0\ 0\ 0\ 0\ 1\ 0) \cup (1\ 0\ 0\ 0\ 0)
\end{aligned}
$$

be the state bivector, except the nodes 10 and $R_1$ all other nodes are in the off state.

$$
\begin{aligned}
XB &= (X_1 \cup X_2)\ (B_1 \cup B_2) \\
&= X_1\ B_1 \cup X_2\ B_2' \\
&\hookrightarrow (1\ 1\ 1\ 1\ 1\ 0\ 0\ 0\ 0\ 1\ 0) \cup (1\ 1\ 0\ 1\ 0\ 0) \\
&= Y_1 \cup Y_2 \text{ (say)}.
\end{aligned}
$$

$$
\begin{aligned}
(Y_1 \cup Y_2)B &= (Y_1 \cup Y_2)\ (B_1 \cup B_2) \\
&= Y_1\ B_1 \cup Y_2\ B_2 \\
&\hookrightarrow (1\ 1\ 1\ 1\ 1\ 0\ 0\ 0\ 0\ 1\ 0) \cup (1\ 1\ 1\ 1\ 0) \\
&= Z_1 \cup Z_2 \text{ (say)}.
\end{aligned}
$$

$$
\begin{aligned}
(Z_1 \cup Z_2)\ (B_1 \cup B_2^t) &= Z_1\ B_1 \cup Z_2\ B_2^t \\
&\hookrightarrow (1\ 1\ 1\ 1\ 1\ 0\ 0\ 0\ 0\ 1\ 0) \cup (1\ 1\ 1\ 1\ 1\ 0) \\
&= (P_1 \cup P_2) \text{ (say)}.
\end{aligned}
$$

$$
\begin{aligned}
(P_1 \cup P_2)\ (B_1 \cup B_2) &= P_1\ B_1 \cup P_2\ B_2 \\
&\hookrightarrow (1\ 1\ 1\ 1\ 1\ 0\ 0\ 0\ 0\ 1\ 0) \cup (1\ 1\ 1\ 1\ 1) \\
&= S_1 \cup S_2 \text{ (say)}.
\end{aligned}
$$

$$
\begin{aligned}
(S_1 \cup S_2)\ (B_1 \cup B_2^t) &= S_1\ B_1 \cup S_2\ B_2^t \\
&= (1\ 1\ 1\ 1\ 1\ 0\ 0\ 0\ 0\ 1\ 0) \cup (1\ 1\ 1\ 1\ 1).
\end{aligned}
$$

Thus the bihidden pattern of the state bivector is a fixed bipoint. Thus the final resultant is given by {(1 1 1 1 1 0 0 0 0 1 0) ∪ (1 1 1 1 1 0), (1 1 1 1 1 0 0 0 0 1 0) ∪ (1 1 1 1 1)}.

We see when the node poverty is in the on state the caretakers feel that PWDs are not given proper health care, they suffer from poor nutrition, improper clothing, no proper shelter, no recreation and the PWDs also suffer from inferiority complex with no self image, discriminated from other family members and no good food or clothing even on festivals. PWDs are frustrated because they have to depend on caretakers for every thing.

The caretakers are also depressed and dejected for having the PWDs as their children. Caretakers have no mind to spend on PWDs. They are ill treated as PWDs are a economic burden and a social disgrace. So caretakers never take any steps to make them get any form of vocations training.

Next we study the on state of the node "improper clothing and poor nutrition i.e., the nodes 3 and 2 in the on state in the FCM component of the FCRM and $D_4$ in the on state of the domain space from the FRM component of the FCR bimodel.

Let

$$X \quad = \quad (0\ 1\ 1\ 0\ 0\ 0\ 0\ 0\ 0\ 0\ 0\ 0) \cup (0\ 0\ 0\ 1\ 0\ 0)$$
$$\quad = \quad X_1 \cup X_2$$

be the given state bivector. To find the effect of X of the dynamical bisystem

$$B \quad = \quad B_1 \cup B_2$$

$$\begin{aligned}
XB \quad &= \quad (X_1 \cup X_2)\ (B_1 \cup B_2) \\
&= \quad X_1\ B_1 \cup X_2\ B_2 \\
&= \quad (0\ 1\ 1\ 0\ 0\ 0\ 0\ 0\ 0\ 0\ 0\ 0) \cup (0\ 0\ 0\ 1\ 0\ 0)\ B_2 \\
&\hookrightarrow \quad (1\ 1\ 1\ 0\ 0\ 0\ 0\ 0\ 0\ 1\ 0) \cup (1\ 1\ 1\ 1\ 1\ 0)
\end{aligned}$$

$$= \quad (Y_1 \cup Y_2) \text{ (say)}.$$

$$(Y_1 \cup Y_2) \, B_s^{t_2} = (Y_1 \cup Y_2) \, (B_1 \cup B_2^t)$$
$$= \quad Y_1 \, B_1 \cup Y_2 \, B_2^t$$
$$= \quad (1\ 1\ 1\ 1\ 1\ 0\ 0\ 0\ 0\ 1\ 0) \cup (1\ 1\ 1\ 1\ 1\ 0)$$
$$= \quad Z_1 \cup Z_2 \text{ (say)}.$$

Now we study the effect of $Z_1 \cup Z_2$ on B.

$$(Z_1 \cup Z_2) \, B = \quad Z_1 \, B_1 \cup Z_2 \, B_2$$
$$= \quad (1\ 1\ 1\ 1\ 1\ 0\ 0\ 0\ 0\ 1\ 0) \cup (1\ 1\ 1\ 1\ 1\ 0)$$
$$\hookrightarrow \quad (1\ 1\ 1\ 1\ 1\ 0\ 0\ 0\ 0\ 1\ 0) \cup (1\ 1\ 1\ 1\ 1).$$

Thus the hidden bipattern of the dynamical bisystem B is a fixed bipoint given by $\{(1\ 1\ 1\ 1\ 1\ 0\ 0\ 0\ 0\ 1\ 0) \cup (1\ 1\ 1\ 1\ 1\ 0), (1\ 1\ 1\ 1\ 1\ 0\ 0\ 0\ 0\ 1\ 0) \cup (1\ 1\ 1\ 1\ 1)\}$.

The main observation from this resultant is that the lack of good clothing for the PWDs are attributes to poverty so we see the resultant happens to be identical as that of when the state vector $(0\ 0\ 0\ 0\ 0\ 0\ 0\ 0\ 0\ 1\ 0) \cup (1\ 0\ 0\ 0\ 0)$ is taken. Both the hidden bipatterns are only fixed bipoints.

## 6.2  Special New FCRM bimodel of the PWDs and a PWD to study the problems of PWDs and the caretakers

In this bimodel for the FCM component we take the special FCM of the PWDs given by the connection matrix M given in chapter four and the FRM component of the bimodel is taken as M given in chapter five of this book. We shall denote M in chapter four by $N_1$ and M in chapter five by $N_2$.

Thus the CR bimatrix associated with the special connection matrix of the 82 PWDs and the FRM of the expert PWD is given by the following bimatrix.

$N = N_1 \cup N_2$

$$= \begin{bmatrix} 0 & 1 & 0 & 0 & 0 & 0 & 0 & 0 & 0 & 1 & 0 \\ 1 & 0 & 0 & 0 & 0 & 0 & 0 & 0 & 0 & 1 & 0 \\ 0 & 0 & 0 & 0 & 0 & 0 & 0 & 0 & 0 & 1 & 0 \\ 0 & 0 & 0 & 0 & 0 & 0 & 0 & 0 & 0 & 0 & 0 \\ 0 & 0 & 0 & 0 & 0 & 0 & 0 & 0 & 0 & 0 & 0 \\ 0 & 0 & 0 & 0 & 0 & 0 & 0 & 0 & 0 & 0 & 0 \\ 0 & 0 & 0 & 0 & 0 & 0 & 0 & 0 & 0 & 1 & 1 \\ 0 & 0 & 0 & 0 & 0 & 0 & 0 & 0 & 0 & 0 & 0 \\ 0 & 0 & 0 & 0 & 0 & 1 & 0 & 0 & 0 & 0 & 0 \\ 1 & 1 & 1 & 0 & 0 & 0 & 0 & 0 & 0 & 0 & 0 \\ 0 & 0 & 0 & 0 & 0 & 0 & 1 & 0 & 0 & 0 & 0 \end{bmatrix} \cup$$

$$\begin{bmatrix} 1 & 1 & 0 & 1 & 0 \\ 1 & 0 & 0 & 1 & 0 \\ 0 & 0 & 1 & 1 & 1 \\ 1 & 1 & 1 & 1 & 0 \\ 0 & 0 & 0 & 1 & 0 \\ -1 & -1 & 0 & 0 & 1 \end{bmatrix}$$

Suppose for instance one is interested in knowing the hidden bipattern when the on state of improper clothing, i.e., node 3 from the special FCM and node $D_4$ from the domain space of the FRM component is taken and all other nodes in the bivector is in the off state.

Let

$$\begin{aligned} X \quad &= \quad X_1 \cup X_2 \\ &= \quad (0\,0\,1\,0\,0\,0\,0\,0\,0\,0\,0) \cup (0\,0\,0\,1\,0\,0) \end{aligned}$$

be the corresponding state bivector. To find the effect of X on the dynamical bisystem $N = N_1 \cup N_2$.

$$
\begin{aligned}
XN \quad &= \quad (X_1 \cup X_2) \cup (N_1 \cup N_2) \\
&= \quad X_1\, N_1 \cup X_2\, N_2 \\
&\hookrightarrow \quad (0\ 0\ 1\ 0\ 0\ 0\ 0\ 0\ 0\ 1\ 0) \cup (1\ 1\ 1\ 1\ 0) \\
&= \quad Y_1 \cup Y_2 \text{ (say)}.
\end{aligned}
$$

Now

$$
\begin{aligned}
(Y_1 \cup Y_2)\,(N_1 \cup N_2)_s^{t_2} \quad &= \quad Y_1\, N_1 \cup Y_2\, N_2^t \\
&\hookrightarrow \quad (1\ 1\ 1\ 0\ 0\ 0\ 0\ 0\ 0\ 1\ 0) \cup (1\ 1\ 1\ 1\ 1\ 0) \\
&= \quad Z_1 \cup Z_2 \text{ (say)}.
\end{aligned}
$$

Now

$$
\begin{aligned}
(Z_1 \cup Z_2)\,(N_1 \cup N_2) &= Z_1\, N_1 \cup Z_2\, N_2 \\
&\hookrightarrow \quad (1\ 1\ 1\ 0\ 0\ 0\ 0\ 0\ 0\ 1\ 0) \cup (1\ 1\ 1\ 1\ 1) \\
&= \quad (T_1 \cup T_2) \text{ (say)}.
\end{aligned}
$$

$$
\begin{aligned}
(T_1 \cup T_2)\,(N_1 \cup N_2)_s^{t_2} &= T_1\, N_1 \cup T_2\, N_2^t \\
&\hookrightarrow \quad (1\ 1\ 1\ 0\ 0\ 0\ 0\ 0\ 0\ 1\ 0) \cup (1\ 1\ 1\ 1\ 1\ 0) \\
&= \quad W_1 \cup W_2 \text{ (say)} \\
&= \quad Z_1 \cup Z_2.
\end{aligned}
$$

Thus the hidden bipattern of the system is a fixed bipoint. Thus when the "node" the PWDs do not get good clothing (improper clothing) is in the on bistate we see the PWDs have "no proper health care" suffer from "poor nutrition", and suffer from "poor economy" which is indicated by the special FCM component.

The resultant given by the FRM component of the bimodel is that the PWDs suffer from inferiority complex, they have no self image, they are discriminated from other family members, frustrated because they have to depend on caretakers for every thing; they suffer because of poverty, caretakers are depressed and dejected for having PWDs as their children.

The caretakers have no mind to spend on PWDs, ill treated as PWDs are economic burden and a social disgrace and finally the caretakers do not take any steps to make them get any form of vocational training.

Thus the hidden bipattern is a fixed bipair given by

$$\{(1\ 1\ 1\ 0\ 0\ 0\ 0\ 0\ 0\ 1\ 0) \cup (1\ 1\ 1\ 1\ 1\ 0),$$

$$(1\ 1\ 1\ 0\ 0\ 0\ 0\ 0\ 0\ 1\ 0) \cup (1\ 1\ 1\ 1\ 1)\}.$$

All these effects cannot be got from a single model. FCRM bimodel alone can give the integrated effect of the bisystem on the on state of a binode.

Now the expert opinion of NGOs and the public about the problems of PWD is given in the following section.

## 6.3 Expert opinion of the NGOs and a public person to analyse the problems of the PWDs using a FCRM bimodel

The special FCM model given in chapter four of this book by the 12 NGOs is taken as the first component of the FCRM. The opinion and analysis of the native of Melmalayanur given in chapter five of this book is taken as the FRM component of the FCRM.

Let $C = C_1 \cup C_2$ be the bimatrix assoiciated with the FCRM bimodel. Here $C_1$ is taken as the matrix N given in chapter four of this thesis and $C_2$ is taken as matrix T given in chapter five of this book.

Thus

$$C = C_1 \cup C_2$$

$$= \begin{bmatrix} 0 & 1 & 0 & 0 & 0 & 0 & 0 & 1 & 1 & 1 & 0 \\ 0 & 0 & 0 & 0 & 0 & 0 & 0 & 1 & 0 & 1 & 0 \\ 0 & 0 & 0 & 0 & 0 & 0 & 0 & 0 & 1 & 1 & 0 \\ 0 & 0 & 0 & 0 & 0 & 0 & 0 & 0 & 0 & 0 & 0 \\ 0 & 0 & 0 & 0 & 0 & 0 & 0 & 0 & 0 & 1 & 0 \\ 0 & 0 & 0 & 0 & 0 & 0 & 0 & 1 & 1 & 0 & 0 \\ 0 & 0 & 0 & 0 & 0 & 0 & 0 & 1 & 1 & 1 & 1 \\ 0 & 0 & 0 & 0 & 0 & 0 & 1 & 0 & 0 & 0 & 0 \\ 0 & 0 & 0 & 0 & 0 & 0 & 0 & 0 & 0 & 1 & 0 \\ 0 & 0 & 0 & 0 & 0 & 0 & 0 & 0 & 0 & 0 & 1 \\ 0 & 0 & 0 & 0 & 0 & 0 & 1 & 0 & 0 & 0 & 0 \end{bmatrix} \cup$$

$$\begin{bmatrix} 0 & 0 & 1 & 1 & 0 & 0 & 0 & 0 & 1 & 0 & 0 \\ 0 & 0 & 0 & 1 & 0 & 0 & -1 & 0 & 0 & 0 & 1 \\ 0 & 0 & 1 & 0 & 0 & -1 & -1 & 1 & 1 & 0 & 1 \\ -1 & 0 & -1 & -1 & 0 & 1 & 1 & -1 & 0 & -1 & -1 \\ -1 & -1 & -1 & -1 & 0 & 1 & 1 & -1 & 0 & 0 & -1 \\ 0 & 1 & 1 & 1 & 1 & -1 & 0 & 0 & 0 & 0 & 0 \\ 0 & 0 & 0 & 0 & 0 & 0 & 0 & 1 & 0 & 1 & 0 \\ 0 & 0 & 0 & 1 & 0 & 0 & 0 & 0 & 1 & 0 & 1 \\ 0 & 1 & 0 & 0 & 0 & 0 & 0 & 0 & 0 & 1 & 0 \end{bmatrix}$$

is the CR bimatrix associated with the FCRM bimodel.

Now suppose one wants to study the effect of the state bivector

$$\begin{aligned} X \quad &= \quad X_1 \cup X_2 \\ &\hookrightarrow \quad (0\,0\,0\,0\,0\,0\,0\,0\,0\,1\,0) \cup \\ &\qquad (1\,0\,0\,0\,0\,0\,0\,0\,0\,0\,0) \end{aligned}$$

where $X_2 = (1\,0\,0\,0\,0\,0\,0\,0\,0\,0\,0) \in R$ the range space of the FRM component of the FCRM.

We now study the effect of $X = X_1 \cup X_2$ on C.

$$
\begin{aligned}
(X_1 \cup X_2)\ C_s^{t_2} \quad &= \quad (X_1 \cup X_2)\ (C_1 \cup\ C_2^t\ ) \\
&= \quad X_1\ C_1 \cup X_2\ C_2^t \\
&\hookrightarrow \quad (0\ 0\ 0\ 0\ 0\ 0\ 0\ 0\ 0\ 1\ 1) \cup \\
&\qquad (0\ 0\ 0\ 0\ 0\ 0\ 0\ 0\ 0) \\
&= \quad Y_1 \cup Y_2\ (\text{say}).
\end{aligned}
$$

Now
$$
\begin{aligned}
(Y_1 \cup Y_2)\ (C_1 \cup C_2) &= \quad Y_1\ C_1 \cup Y_2\ C_2 \\
&\hookrightarrow \quad (0\ 0\ 0\ 0\ 0\ 0\ 1\ 0\ 0\ 1\ 1) \cup \\
&\qquad (1\ 0\ 0\ 0\ 0\ 0\ 0\ 0\ 0\ 0\ 0) \\
&= \quad Z_1 \cup Z_2\ (\text{say}).
\end{aligned}
$$

$$
\begin{aligned}
(Z_1 \cup Z_2)\ (C_1 \cup C_2^t\ ) &\hookrightarrow (0\ 0\ 0\ 0\ 0\ 0\ 1\ 1\ 1\ 1\ 1) \cup \\
&\qquad (0\ 0\ 0\ 0\ 0\ 0\ 0\ 0\ 0) \\
&= \quad T_1 \cup T_2\ (\text{say}).
\end{aligned}
$$

$$
\begin{aligned}
(T_1 \cup T_2)\ (C_1 \cup C_2) &\hookrightarrow (0\ 0\ 0\ 0\ 0\ 0\ 1\ 1\ 1\ 1\ 1) \cup \\
&\qquad (1\ 0\ 0\ 0\ 0\ 0\ 0\ 0\ 0\ 0\ 0) \\
&= \quad P_1 \cup P_2\ (\text{say}).
\end{aligned}
$$

We see the hidden bipattern of the system is a fixed point. According to the NGO and the public person poverty is the root cause of unemployment, they are unaware of the SSHG group which can render them help, welfare measure of government never reach them as they are ignorant of the government scheme which help PWDs.

Marriage of the PWDs remain a question mark. From the public person it is clear that poverty does not have any impact on other nodes from both the domain space or the range space.

Thus the public person who is from their own panchayat and block feels and confidently says poverty has nothing to do with all the nodes described in the domain and the range space of the FRM. Thus we see something is a cause for these problems. From our discussions we see it is mainly the caretakers fear of social stigma, which is a cause of majority of the problems for the PWDs.

We now use the following state bivector

$$A \quad = \quad A_1 \cup A_2$$
$$\hookrightarrow \quad (0\ 0\ 0\ 0\ 0\ 1\ 0\ 0\ 0\ 0\ 0) \cup$$
$$(0\ 1\ 0\ 0\ 0\ 0\ 0\ 0\ 0\ 0\ 0)$$

where the nodes PWDs have no proper education is taken in the on state for the FCM component of the FCRM and the relatives are ashamed of having a PWD in their family is taken as the only node in the on state from the range space of the FRM.

To find the effect of the bistate vector on the dynamical bisystem C.

$$A\,C_s^{t_2} \quad = \quad (A_1 \cup A_2)\,(C_1 \cup C_2^t)$$
$$= \quad A_1\,C_1 \cup A_2\,C_2^t$$
$$\hookrightarrow \quad (0\ 0\ 0\ 0\ 0\ 1\ 0\ 1\ 1\ 0\ 0) \cup$$
$$(0\ 0\ 0\ 0\ 0\ 1\ 0\ 0\ 1)$$
$$= \quad B_1 \cup B_2 \text{ (say)}.$$

Now

$$(B_1 \cup B_2)\,(C_1 \cup C_2) = \quad B_1\,C_1 \cup B_2\,C_2$$
$$\hookrightarrow \quad (0\ 0\ 0\ 0\ 0\ 1\ 1\ 1\ 1\ 1\ 0) \cup$$
$$(0\ 1\ 1\ 1\ 1\ 0\ 0\ 1\ 1\ 0\ 1)$$
$$= \quad R_1 \cup R_2 \text{ (say)}.$$

Now
$$(R_1 \cup R_2)\,(C_1 \cup C_2^t) \quad \hookrightarrow \quad (0\ 0\ 0\ 0\ 0\ 1\ 1\ 1\ 1\ 1\ 1) \cup$$
$$(1\ 1\ 1\ 0\ 0\ 1\ 1\ 1\ 1)$$
$$= \quad S_1 \cup S_2 \text{ (say)}.$$

Now we see the effect of $S_1 \cup S_2$ on the dynamical bisystem $C = C_1 \cup C_2$

$$
\begin{aligned}
SC \quad &= \quad (S_1 \cup S_2)(C_1 \cup C_2) \\
&= \quad S_1\, C_1 \cup S_2\, C_2 \\
&\hookrightarrow \quad (0\ 0\ 0\ 0\ 0\ 1\ 1\ 1\ 1\ 1\ 1) \cup \\
&\qquad (0\ 1\ 1\ 1\ 1\ 0\ 0\ 1\ 1\ 1\ 1\ ) \\
&= \quad X_1 \cup X_2 \text{ (say).}
\end{aligned}
$$

$$
\begin{aligned}
X_1 \cup X_2\,(C1 \cup C_2^t) \quad &\hookrightarrow \quad (0\ 0\ 0\ 0\ 0\ 1\ 1\ 1\ 1\ 1\ 1) \cup \\
&\qquad (1\ 1\ 1\ 0\ 0\ 1\ 1\ 1\ 1) \\
&= \quad Y_1 \cup Y_2 \text{ (say).}
\end{aligned}
$$

Thus the bihidden pattern of the dynamical system is a fixed bipoint, given by the bipair $\{(0\ 0\ 0\ 0\ 0\ 1\ 1\ 1\ 1\ 1\ 1) \cup (0\ 1\ 1\ 1\ 1\ 0\ 0\ 1\ 1\ 1\ 1), (0\ 0\ 0\ 0\ 0\ 1\ 1\ 1\ 1\ 1\ 1) \cup (1\ 1\ 1\ 0\ 0\ 1\ 1\ 1\ 1)\}$.

Finally we give a FCRM bimodel in which the attributes of the FCM component of the FCRM bimodel coincides with the domain space of the FRM component of the FCRM bimodel.

## 6.4 Experts opinion using FCRM model to study the problems faced by the PWDs

Now two experts give their opinion. The first expert is a PWD's relative. He gives his opinion about the problems of a PWD and using FCMs he wants to analyse the problem. The attributes given by him for the analysis of the problems of the PWD is as follows.

$R_1$  -  PWDs suffer from acute inferiority complex.

$R_2$  -  They are unemployed / no means to earn

$R_3$  -  Ill treated by caretakers

$R_4$  -  Happy and contented

$R_5$ - PWDs suffer from mental depression.

$R_6$ - Denied good food, new clothing, proper shelter and recreation.

The other expert is a SSHG person working for PWDs. He wishes to work with FRM component of the FCRM. He has accepted the attributes of the PWD's relative viz. $R_1$, $R_2$, ..., $R_6$ as the domain space of the FRM. The attributes given for him for the range space is as follows:

$S_1$ - Fatalism and Karma.

$S_2$ - Poverty a cause of neglect

$S_3$ - Discriminated from others

$S_4$ - Relatives are ashamed of the PWDs

$S_5$ - PWDs are a curse upon the family

$S_6$ - PWDs are a burden

$S_7$ - Caretakers are sympathetic and carrying.

Now using these attributes FCRM bigraph of these experts is given by the following figures.

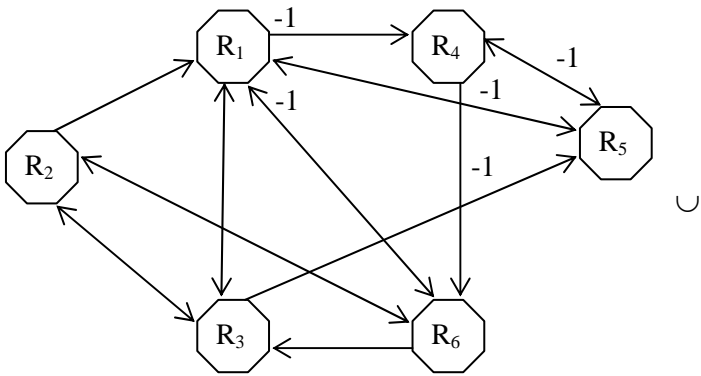

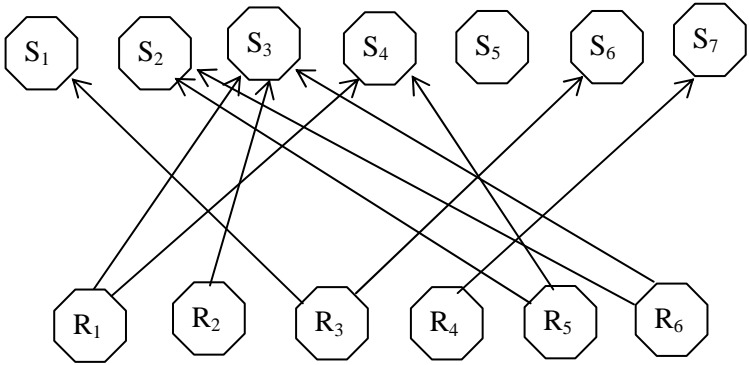

$$V = \begin{array}{c} \\ R_1 \\ R_2 \\ R_3 \\ R_4 \\ R_5 \\ R_6 \end{array} \begin{array}{cccccc} R_1 & R_2 & R_3 & R_4 & R_5 & R_6 \\ \left[ \begin{array}{cccccc} 0 & 0 & 0 & -1 & 1 & 1 \\ 1 & 0 & 1 & 0 & 0 & 1 \\ 1 & 1 & 0 & 0 & 0 & 0 \\ 0 & 0 & 0 & 0 & -1 & -1 \\ -1 & 0 & 1 & -1 & 0 & 0 \\ 1 & 1 & 1 & 0 & 0 & 0 \end{array} \right] \end{array} \cup$$

$$\begin{array}{c} \\ R_1 \\ R_2 \\ R_3 \\ R_4 \\ R_5 \\ R_6 \end{array} \begin{array}{ccccccc} S_1 & S_2 & S_3 & S_4 & S_5 & S_6 & S_7 \\ \left[ \begin{array}{ccccccc} 0 & 0 & 1 & 1 & 0 & 0 & 0 \\ 0 & 0 & 1 & 0 & 0 & 0 & 0 \\ 1 & 0 & 0 & 0 & 0 & 1 & 0 \\ 0 & 0 & 0 & 0 & 0 & 0 & 1 \\ 1 & 0 & 0 & 1 & 0 & 0 & 0 \\ 0 & 1 & 1 & 0 & 0 & 0 & 0 \end{array} \right] \end{array}$$

$= V_1 \cup V_2$ .

$V = V_1 \cup V_2$ is the related CR bimatrix associated with the bigraph with the bigraph given in figure.

Now if we consider the state bivector

$$X \qquad = \qquad (1\ 0\ 0\ 0\ 0\ 0) \cup (1\ 0\ 0\ 0\ 0\ 0\ 0)$$

where the nodes $R_1$, PWDs suffer from acute inferiority complex alone is in the on state from the domain space and $S_1$ - Fatalism and Karma alone is in the on state in the range space.

To find the effect of state bivector X on the dynamical bisystem V.

$$
\begin{aligned}
XV \qquad &= \quad (X_1 \cup X_2)\,(V_1 \cup V_2) \\
&= \quad X_1\,V_1 \cup X_2\,V_2 \\
&\hookrightarrow \quad (1\ 0\ 0\ 0\ 1\ 1) \cup (0\ 0\ 1\ 1\ 0\ 0\ 0) \\
&= \quad Y_1 \cup Y_2 \text{ (say).}
\end{aligned}
$$

$$
\begin{aligned}
Y\ V_s^{t_2} \qquad &= \quad (Y_1 \cup Y_2)\,(V_1 \cup V_2^t) \\
&\hookrightarrow \quad (1\ 1\ 1\ 0\ 0\ 1) \cup (1\ 1\ 0\ 0\ 1\ 1) \\
&= \quad Z_1 \cup Z_2 \\
&= \quad Z \text{ (say).}
\end{aligned}
$$

$$
\begin{aligned}
ZV \qquad &= \quad (Z_1 \cup Z_2)\,(V_1 \cup V_2) \\
&= \quad Z_1\,V_2 \cup Z_2\,V_2 \\
&\hookrightarrow \quad (1\ 1\ 1\ 0\ 1\ 1) \cup (1\ 1\ 1\ 1\ 0\ 0\ 0) \\
&= \quad T_1 \cup T_2 \\
&= \quad T \text{ (say).}
\end{aligned}
$$

$$
\begin{aligned}
T\ V_s^{t_2} \qquad &= \quad (T_1 \cup T_2)\,(V_1 \cup V_2^t) \\
&= \quad T_1\,V_1 \cup T_2\,V_2^t \\
&\hookrightarrow \quad (1\ 1\ 1\ 0\ 1\ 1) \cup (1\ 1\ 1\ 1\ 1\ 1) \\
&= \quad (P_1 \cup P_2) \\
&= \quad P \text{ (say).}
\end{aligned}
$$

$$
\begin{aligned}
\text{PV} \quad &= \quad (P_1 \cup P_2)\,(V_1 \cup V_2) \\
&= \quad P_1\,V_1 \cup P_2\,V_2 \\
&\hookrightarrow \quad (1\ 1\ 1\ 0\ 1\ 1) \cup (1\ 1\ 1\ 1\ 0\ 1\ 1) \\
&= \quad B_1 \cup B_2 \ (\text{say}).
\end{aligned}
$$

The hidden bipattern of the state vector is a fixed bipoint given by $\{(1\ 1\ 1\ 0\ 1\ 1) \cup (1\ 1\ 1\ 1\ 0\ 1\ 1),\ (1\ 1\ 1\ 0\ 1\ 1) \cup (1\ 1\ 1\ 1\ 1)\}$. Thus only the nodes $R_4$ - Happy and contented from the FCM component and $S_5$ from the FRM component of the FCRM remains in the off state. All other binodes from both the spaces come to on state.

Thus the conclusions given in the last chapter has carried several combinations of bistate vectors.

Chapter Seven

# FUZZY LINGUISTIC MODELS TO ANALYSE THE PROBLEMS OF MENTALLY CHALLENGED PEOPLE

Here in this chapter we analyse the problems related with mentally challenged people using fuzzy linguistic models. This chapter has three sections. Section one describe and develops the problems faced by the care takers / relatives of mentally challenged people. This section describes the problem of children who are mentally challenged and who suffer with immobility and situations where a care taker or a relative must always be present by his / her side. That is 'critical children'. Of course we have 'critical adults' who may not be born as a mentally challenged one but would have come to this state due to accident or some disease or old age. Section two recalls the basic notions related with fuzzy linguistic models and applies FLCM model and FLRMs model to this problem. Section three gives the conclusions and suggestions based on the analysis using these models and interviews.

## 7.1 Introduction

We have so far discussed and studied in this book mostly about physically challenged people and have given suggestion how to tackle the problems faced by them. In this chapter we discuss and describe the mentally challenged people more so children born with mental retardation. When we say mental retardation we mainly deal with those children and not those who have become so by some disease or accident or shock in their life. For the three categories mentioned later are people who have been understood, tended and taken care of so even after their encountering the state of mental disability they may be taken care depending on the circumstances they face in the society; for they are accepted by the society and not shunned but only pitied; so in this case they do not hold any social stigma for every one close to the family is well aware about their mental disability so the problem of marriage of their close relatives are not in any way affected.

Thus we in this chapter the authors only deal with mental disability by birth this includes also hyperactive children when not given proper medical attention who at times while attaining the adolescence become mentally disabled and hence suddenly loose their very cognitive senses which is very much visible from their gait and activities.

More problems are faced by the parents of mentally disabled children for they need special care and constant help.

We can broadly classify the mentally challenged children into the following: Born as a 'vegetable'. That is the children from the time of birth are incapable of any cognitive movement and they cannot recognize their parents or any relative. They eat (that is being fed by their relatives) and grow without any movement that is they cannot sit or stand or lift any article or hold any article. Their limbs movement are not under their control. Such are the worst type of mental disability which need constant attention and very special care. They grow in size but with no growth of their mind or cognitive senses.

This occurs due to the following reasons (their may be more reasons known to the medical expert) gathered by us which may be rational and some irrational for it is the information given by the close relatives caretakers of such type of children.

1.  At the time of delivery of the infant the oxygen is totally denied due to lack of medical expertise, which is irreversible and hence the infant completely looses the brain capacity or at time partially looses the brain capacity or any other complications at the time of delivery.

2.  Both parents / one of the parents happen to be alcoholic at the time of conception, this also may result in mentally retarded children.

3.  Marriage only from relatives which also results in very poor brain capacity and intelligence in general though not complete mental retardation.

4.  Genetic disorder which has been inherited from fore fathers and so on.

The above information is gathered by the experts who were collecting data from such families / non governmental organizations / hospitals.

However in this book we would by no means mention the names of any such place from which the information is gathered for we have promised confidentiality about the information obtained by us. Now for the problem, information was gathered from caretakers of the mentally retarded children, parents / brothers / sisters of these mentally retarded children, servants attached to their family for generations, NGOs who work for them, other organizations who work for them and doctors specialized in this field. Further some of the experts are from the special physical trainees of these children as well as special educational trainees of these children.

These experts had broadly classified the mental disability into the following types which is for our own study and it is not medically done by the experts in this field.

Type I.    This is the worst critical mental condition where the brain in the child does not function even a little. The child is in the 'vegetable' state lying all the time with inability to move even their limbs, in certain cases they move their limbs look at the caretaker and only make some sounds not able to speak. They shout for food or for any other inconvenience as the "natures call". They have to be fed by the caretakers and 24 hours of monitoring is needed for some suffer additional health problems.

Taking care of these children happens to be a big challenge for if they are poor they cannot take care of them for they have to go to work to earn their daily living. Further when the child is in this state from the time of birth it is very difficult to maintain the child in their home and give proper care.

Now in two ways the poor people tackle this situation. Either they throw the child in river or burry it after killing it or throw the child in the dustbin if they are city dwellers. Some just keep them in the house in the backyard and attend the child whenever they get time.

Infact only those charity organizations take care of these children. These children do not have a very long life they may live maximum upto 10 years; this information was given by some experts from the charity organizations and NGOs who work for them they die mainly due to lack of proper nutrition and medical care.

Some elders / adults reach this stage due to old age or disease or accidents where the brain gets damaged. Such are taken care of in most cases by their close relatives and they come under "worst critical condition".

Type II. This is the "bad critical mental condition" where the brain exists to 'some' extent; the child can recognize the caretakers, respond to some extent can be made to sit and at times they can crawl or move in a way they like. However they cannot stand or they cannot identify several things or their cognitive / associative part of the brain is dead or never functions. They cannot express any thing by words / speaking they also shout, the shout may vary so that the caretaker can understand the need of the child. Such children are a shade better than the worst critical situation. One does not know whether they respond to danger or any problem they can face!.

Type III. These children may be from birth or just when they enter teens become mentally disabled. They can walk with some support or their gait is very unstable so cannot walk alone. They cannot read or write or speak. They can recognize their caretakers. They are 'critical' in mental retardation. They survive and reach middle age. However they cannot be employed but can be trained to do their routine. They also do not talk only shout or produce noise. One is not able to understand why the mentally disabled lose their speech or cannot speak / talk properly. However constant care is to be given to them. For if neglected, they sitting in one place may even forget to walk or move their limbs for the movement of limbs in a cognitive way alone can make them to be self sufficient. Thus the caretakers say they should be given regular psyotherapy and exercise to keep them fit. They cannot be introduced or steam lined into any form of formal education. Some of the hyper active children after a stage turn to be this type. This was acknowledged by 9 caretakers that these children upto the age of 9 to 10 years were only hyperactive attended school regularly but suddenly became like this after some health problems and are unable to understand or recognize or read. Some experts feel as no medical aid was given to them for the hyperactivity so only they last their mind. Their gait and cognitive movements have damaged beyond repair. Some doctors say if they are taken special care from a very young age this could be averted not fully but atleast partially. This happening is mainly related with the genetics of their close relatives or some brain damage due to disease so at that stage noting can be done to revert them back to normalcy.

Type IV. These are hyperactive children (autism) who are slow learners and also difficulty in educating them their very gait and activities show they are different from normal. Such children can graduate and also can be employed and special type of attention ought to be given to them to attain this goal.

Type V. This includes adults who due to shock, deep depression or some other causes become mad or loose their mental balance; of which some are correctable and others non correctable but can be treated, at times they respond to treatment.

Now we will dealing be with these five types of mental disability and the problems faced by the caretakers using fuzzy linguistic models. When we say "problems" it includes social problems, economic problems and other problems which cannot be openly said by silently suffered.

## 7.2 Description of the Fuzzy Linguistic Relational Maps Model

Here we briefly describe the Fuzzy Linguistic Relational Maps (FLRMs) model.

Let $X = \{x_1, \ldots, x_n\}$ and $Y = \{y_1, \ldots, y_m\}$ ($m \neq n$) be two distinct and disjoint fuzzy linguistic sets with 0. We assume any two elements in $X$ may or may not be comparable and any two elements in $Y$ also may or may not be comparable. We call $X$ the domain space of linguistic fuzzy terms and $Y$ the range space of fuzzy linguistic terms. We denote by $S$ a set of linguistic terms / linguistic values which are comparable and contains 0. For any $x_i \in X$ we find its relation or its effect on $y_k \in Y$. If no relation exist we just mark '0'.

If a relation exists and say $s_j$ we just mark $x_i \xrightarrow{s_j} y_k$ where $s_j$ is a fuzzy linguistic term in $S$ and not a value in [0, 1]. We find relations of every element in $X$ on $Y$. This sort of relation can be described by a graph with vertices from $X$ and $Y$ and its edge values are from $S$. We define this graph as the fuzzy linguistic graph. Now let $M$ denote a $n \times m$ fuzzy linguistic matrix associated with the fuzzy linguistic graph. We call $M$ as the fuzzy linguistic dynamical system associated with the FLRM model.

Now if $T = (s_1, \ldots, s_n)$ where $s_i \in S$, $1 \leq i \leq n$ is called the fuzzy state linguistic vector associated with $X$. Similarly

V = {(s_1, ..., s_m)/s_j ∈ S, 1 ≤ j ≤ m} will be the fuzzy state linguistic vector associated with Y.

Now the effect of T on M is studied by the following way.

$$\text{If } M = \begin{bmatrix} m_{11} & ... & m_{1m} \\ m_{21} & ... & m_{2m} \\ \vdots & & \vdots \\ m_{n1} & ... & m_{nm} \end{bmatrix}$$

where $m_{ij} \in S$; $1 \leq i \leq n$ and $1 \leq j \leq m$ then T o M = min (min $\{s_i \, m_{ij}\}$) [(min or $\{max \{s_i \, m_{ij}\}\}$)] according to the experts need and the problem under investigation. Clearly (T o M) = V is a $1 \times m$ fuzzy linguistic vector with entries from S. We now find V o $M^t$ = $T_1$ and so on until we reach a fixed point or a limit cycle. We call this the fuzzy linguistic hidden pattern of the system M.

Thus in case of fuzzy linguistic relational map model we have a pair of state vectors to be the fuzzy linguistic hidden pattern. We will describe this by a simple example.

***Example 7.2.1:*** Let us study the performance of students and the teachers using FLRM model.

The performance of students are S = {0, best, good, better, poor, worst, very poor} and X = {good student, poor student, average student, best student, worst student} be the attributes related with the student and Y = {good teacher, average teacher, bad teacher, best teacher, devoted teacher, unconcerned teacher about students and kind teacher} be the attributes associated with a teacher taken as the range space of the fuzzy linguistic attributes.

X is taken as the domain fuzzy linguistic space. The fuzzy linguistic attributes associated with X is as follows which indicate the performance of students in studies.

$S_1$  -  best student
$S_2$  -  good student
$S_3$  -  average student
$S_4$  -  poor student
$S_5$  -  worst student

The fuzzy linguistic attributes associated with Y describing the teachers teaching ability is as follows:

$T_1$  -  best teacher
$T_2$  -  good teacher
$T_3$  -  devoted teacher
$T_4$  -  average teacher
$T_5$  -  kind teacher
$T_6$  -  bad teacher
$T_7$  -  teacher unconcerned about students.

The student performance graph given by an experts is given in the following:

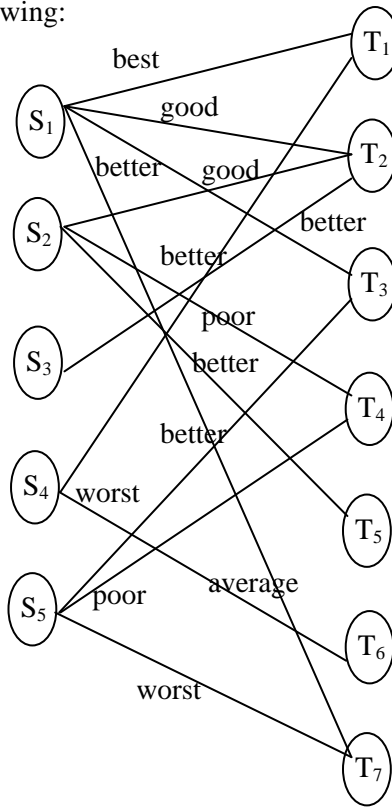

Likewise we can get the fuzzy linguistic graph using an expert.

The fuzzy linguistic matrix associated with the fuzzy linguistic graph is as follows:

$$
M = \begin{array}{c}
 \\
S_1 \\
S_2 \\
S_3 \\
S_4 \\
S_5
\end{array}
\begin{array}{ccccccc}
T_1 & T_2 & T_3 & T_4 & T_5 & T_6 & T_7 \\
\left[\begin{array}{ccccccc}
\text{best} & \text{good} & \text{better} & 0 & 0 & 0 & \text{average} \\
0 & \text{good} & \text{poor} & \text{better} & 0 & 0 & 0 \\
0 & \text{better} & 0 & 0 & 0 & 0 & 0 \\
\text{better} & 0 & 0 & 0 & 0 & \text{worst} & 0 \\
0 & 0 & \text{better} & \text{poor} & 0 & 0 & \text{worst}
\end{array}\right]
\end{array}
$$

It is important to mention here that we have not filled in all the rows we just give this only as an illustration for we want to show how the dynamical system functions.

Now every element (that is state vector) in the domain space is a 5 tuple the first tuple corresponds to the state of the best student influenced by the teachers attitude and depending on the teachers nature; likewise the second tuple corresponds to the state of the good students state depending on the teachers attitude and so on.

Suppose we take a state vector of performance of the student P = (average, 0, 0, 0, worst) related to the domain space. To find the effect of P on M is as follows.

$$
\begin{aligned}
P \circ M \ &= \max \{\min \{p_i, m_{ij}\}\} \\
&= (\text{average, average, average, worst, 0, 0, worst}) \\
&= Y \ (\text{say}).
\end{aligned}
$$

Now we find $Y \circ M^t$ using the same max min operation.

$Y \circ M^t =$ (average, average, average, average, worst) $= P_1$.
$P_1$ is updated and so on. We will arrive at a fixed point or a limit cycle as the number of elements in S is finite so we get the resultant in a finite number of steps.

It is pertinent to mention here the result depends on the operation if we want the best results we use max-max operation medium result we use max-min or min-max, the worst possible solution we use the min-min operation.

It is in the hands of the expert to choose any of the operations depending on the type of results he expects from the problem / study / analysis.

Analysis of the problems faced by the caretakers / parents of the mentally disabled children using FLRM model and FLCM model is described in the following. We have at the outset made clear all information regarding the mentally retarded children can be got only from the parents / caretakers / brothers or sisters / close relatives. So we have discussed with them and obtained information. First we wish to state for type I mentally retarted children there was not even a single kith or kin who spoke or gave information; for all the persons whom we interviewed for type I children was from the charity homes who were caretakers in the charity homes. They only gave complete information parentage and other information related with these children. We nearly interviewed 21 type I children from Chennai and around Chennai the first observation was that these children were left in road sides / dust bin and saved by the charity groups. The second was some of these children had very rich heritage and those rich parents funded the charity home and helped in the running of these homes to some extent. However it was surprising to find that there was no children belonging to middle class or lower middle class or poor or from rural area. These children were found only in Chennai. Some NGOs / social workers said in the rural areas if they found that their child was type I, they would kill the child mainly fearing the social stigma that is why we cannot keep them alive they said. In city these people should have thrown such children in the dustbin or left in train or in the station or in public places like bus stand. The only claim by the experts was that even to accept the type I mentally retarded children's parents, as their parents was desperate for these children do not know any thing, love of parents or the very concept of parents; in such situation the discussion over these type I children was painful and very sad for we cannot arrive at any conclusion and most of them did grow large in size only medical researchers should do a lot of research about these children and find some medical discovery so that such children were aborted before birth and to save the existing ones and as well as prevent the birth of such children. Of course taking care of them needs a lot of mental strength and a very special heart, when their own parents and very close relatives do not like to take care of them for the reasons best known to them. It is also important to note that it is very unfortunate to see these children crossing 10 years of age which is mainly attributed to diseases and improper nutrition they get from lying in one posture etc. Medical experts claim that there is no cure for mental disability.

Type II children were less critical they can respond sit or crawl and can be made to sit, walk with full support in some cases. These children are in general not left in the hands of caretakers in the charity organization. They are cared by close relatives in some of the cases.

In poor economic strata in rural areas these children are chained and kept either in the backyard of the house or in the entrance whenever the relatives go for work. In case of rich and middle class living in city these children are taken care of at home by caretakers. One of the observation is that if brother is mentally retarded then his sister is the one who cares the maximum for him and vice versa. This sort of concern needs a psychological study. However these children lived longer than the type I children.

Regarding type III children, they are mentally disabled but they can walk move about and be self sufficient with constant assistance. They cannot be given the usual schooling for they cannot talk but can be trained to do physical work and that too depends on training. Since they are not violent they can be trained not only to be self sufficient but also can be employed which involves only physical labour. But these children constantly need to be given physical training otherwise the care takers say they have chances of forgetting everything and coming to the state of dependency. They live longer than type I children.

Finally the type IV children are the hyper sensitive kids who are slow learners in adolescence and need some way of training. This kind of kids cannot be recognized from the usual at first sight but can be easily spotted by the experts. They can by constant training from the very early stages and with proper medication they can be made to function in a normal way.

It is pertinent to mention that the division of these mentally disabled children into four types done by the authors in consultation with the experts is only to make use of the new fuzzy linguistic relational model and fuzzy linguistic cognitive maps and not in any way these four types are divided medically or clinically. We have made this division to adopt to our FLRM model.

$P_1$ - Social Stigma: It is very important to note that any family with a mentally retarded child suffers a social stigma and are ostracized by their close relatives and the neighbours who live around them.

One is not in a position to understand the reason. This sort of ostracism we claim as social stigma. Further parents of usual children do not allow their children to be friendly or say even to see or play with these children.

$P_2$ - Economic Burden: In the first place these children are not going to be in any way help to the family members, but the family has to spend on them and in cases of training them or sending them to special schools which is an economic burden on the family.

If they are from the lower economic strata they certainly feel this more than the other economic stratas.

$P_3$ - Marriage of siblings near to impossibility: We are not able to understand / unravel the cause behind this. We got some vague reply that they would get such children if they opt to make alliances in such families. Some disposed it of by saying this would certainly have a genetic impact if not in this generation in the next generation.

So we want to avoid such problems said the persons whom we interviewed.

$P_4$ - Family in a state of despair: Most of the parents feel very desperate for they think they can take care of the child but after their death no one will take care of them as they cannot dress on their own or attend 'natural call' or eat on their own. Though they can see yet while they walk they just dash against every thing such is their state of mind for they lack cognitive action. So the parents / family of these children are in despair. They cannot express any problem or need.

$P_5$ - Constant attention is an impossibility / painful issue: The caretakers / parents feel that constant attention on these children who are so dependent is a painful issue, for they try to walk and dash against anything and fall down which makes more a misery so taking care of them at all times is a painful issue; they claimed.

$P_6$ - Partial / full time organization to take care of them: These parents / caretakers said if in the near by area they had a special school for them, they would put them in the day time in such schools; such schools were meagre. Some of the city dwellers cited that there around 250 mentally retarded children in the Saidapet area in Chennai and a school run privately to train them and so far 150 of them had been trained by this school claimed some experts and training them was very difficult they added. Further parents are very hesitant to send their children to this school fearing a social stigma.

$P_7$ - Training / making them organising them was very difficult.

They (caretakers / parents) claimed it was impossible for them to train / organize them to be self sufficient. For it needs specially trained persons to train them as the retardation is medically incurable.

$P_8$ - Time constraint in a busy schedule.

$P_9$ - Government can run full time or part time homes for them.
$P_{10}$ - More research in this direction to correct / treat / prevent metal retardation problem is to be carried out.
$P_{11}$ – Feelings of caretakers with mental retardation.

They felt even if they had the mind and heart to spend with them in the busy schedule they found it very difficult to spend time on them. Of course the parents liked to be with them but it was in vain.

Now we would be using the following fuzzy linguistic set L to study this problem. L = {0, low, medium, high, very low, just low, very high} is taken as the fuzzy linguistic set from which the edge values are taken. Here we see the term high to show the highest stress, low the lowest stress, medium stress and so on. L measures the stress of the caretakers now, also we will be using the same L to show high social stigma, high economic burden and so on. So these fuzzy linguistic terms are flexible enough to give results in a convenient and a sensitive way in which the relations exist.

These 11 attributes are taken and the Fuzzy Linguistic Cognitive Maps (FLCM) model is constructed [50].

The fuzzy linguistic graph given by the expert is as follows:

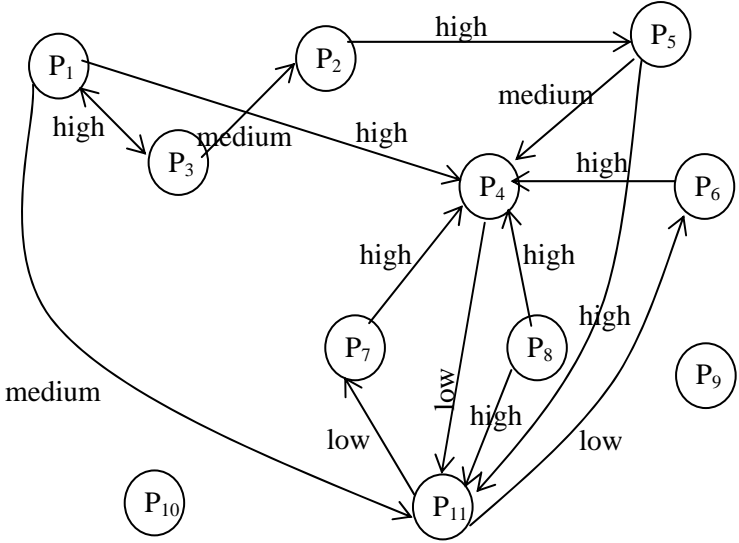

Using this fuzzy linguistic graph we have the following fuzzy linguistic matrix M which serves as the dynamical system for the fuzzy linguistic model. M =

|        | $P_1$ | $P_2$  | $P_3$ | $P_4$   | $P_5$ | $P_6$ | $P_7$ | $P_8$ | $P_9$ | $P_{10}$ | $P_{11}$ |
|--------|-------|--------|-------|---------|-------|-------|-------|-------|-------|----------|----------|
| $P_1$  | 0     | 0      | high  | high    | 0     | 0     | 0     | 0     | 0     | 0        | medium   |
| $P_2$  | 0     | 0      | 0     | 0       | high  | 0     | 0     | 0     | 0     | 0        | 0        |
| $P_3$  | high  | medium | 0     | 0       | 0     | 0     | 0     | 0     | 0     | 0        | 0        |
| $P_4$  | 0     | 0      | 0     | 0       | 0     | 0     | 0     | 0     | 0     | 0        | low      |
| $P_5$  | 0     | 0      | 0     | medium  | 0     | 0     | 0     | 0     | 0     | 0        | high     |
| $P_6$  | 0     | 0      | 0     | high    | 0     | 0     | 0     | 0     | 0     | 0        | 0        |
| $P_7$  | 0     | 0      | 0     | high    | 0     | 0     | 0     | 0     | 0     | 0        | 0        |
| $P_8$  | 0     | 0      | 0     | high    | 0     | 0     | 0     | 0     | 0     | 0        | 0        |
| $P_9$  | 0     | 0      | 0     | 0       | 0     | 0     | 0     | 0     | 0     | 0        | 0        |
| $P_{10}$ | 0   | 0      | 0     | 0       | 0     | 0     | 0     | 0     | 0     | 0        | 0        |
| $P_{11}$ | 0   | 0      | 0     | 0       | 0     | low   | low   | 0     | 0     | 0        | 0        |

Using this fuzzy linguistic matrix we can derive the effect of the on state of any fuzzy linguistic term on the dynamical system.

Let X = (high, 0, 0, 0, 0, 0, 0, 0, 0, 0, 0) be the fuzzy linguistic state vector in which only the state social stigma - $P_1$ is in the on state all other fuzzy linguistic terms are in the off state.

To find the effect of X on M

$$X \circ M = \max \min \{(X, M)\}$$
$$\hookrightarrow \text{(high, 0, high, high, 0, 0, 0, 0, 0, 0, medium)}$$
$$= X_1 \text{ (say)}.$$

($\hookrightarrow$ denotes the resultant has been updated).

Now $X_1 \circ M$ = max {min $(X_1, M)$}

↪ (high, medium, high, high, 0, low, low, 0, 0, 0, high)

= $X_2$ (say).

Consider $X_2 \circ M$ = max {min $(X_2, M)$}

↪ (high, medium, high, high, medium, low, low, 0, 0, 0, high)

= $X_3$.

We see $X_3 \circ M = X_3$ is a fixed fuzzy linguistic point. Hence if the fuzzy linguistic state vector X in which only the state "social stigma" $(P_1)$ is in the on state with the fuzzy linguistic attribute 'high' and all other vectors are in the 'off' state then the effect of X on M is as follows:

Thus X o M gives a hidden pattern which is a fixed point. We see if social stigma alone is in the on state then there is a medium suffering economically for the state $P_2$ is medium.

The marriage of the siblings is an impossibility is high i.e., the $P_3$ attribute in the state vector is high. Family being in the state of despair due to the social stigma is also high; that is $P_4$ is high. Constant attention $P_5$ is "medium" for we see some cases the mentally retarded children are just tied to a pole just like animals, so the constant attention remains medium. Partial / full time organization to take care of them is 'low' that is the '$P_6$' state attribute is low. Also training them or making them organized is very difficult is also 'low'; that is $P_7$ state vector is also 'low'. The fuzzy linguistic attributes $P_8$, $P_9$, $P_{10}$ is in the 0 state that is off state. Finally social stigma due to the mentally retarded children is 'medium'.

Thus we see the on state of the fuzzy attribute $P_1$ is high gives $P_2$, $P_3$, $P_4$, $P_5$, $P_6$, $P_7$ and $P_{11}$ in the on state (refer the fuzzy linguistic state vector $X_3$).

Suppose the expert considers the fuzzy linguistic state vector

S = (0, medium, 0, 0, 0, 0, 0, 0, 0, 0, 0)

in which only the fuzzy linguistic state vector $P_2$ is in the on state with 'medium' as its fuzzy linguistic term.

The effect of S on the fuzzy linguistic dynamical system M is as follows.

$$
\begin{aligned}
\text{S o M} \;\; &= \;\; \max\{\min(S, M)\} \\
&\hookrightarrow \;\; (0, \text{medium}, 0, 0, \text{medium}, 0, 0, 0, 0, 0, 0) \\
&= \;\; S_1 \text{ say.}
\end{aligned}
$$

$$
\begin{aligned}
S_1 \text{ o M} \;\; &\hookrightarrow \;\; (0, \text{medium}, 0, \text{medium}, \text{medium}, 0, 0, 0, 0, 0, \\
&\qquad \text{medium}) \\
&= \;\; S_2 \text{ (say).}
\end{aligned}
$$

$$
\begin{aligned}
S_2 \text{ o M} \;\; &\hookrightarrow \;\; (0, \text{medium}, 0, \text{medium}, \text{medium}, \text{low}, \text{low}, 0, \\
&\qquad 0, 0, \text{medium}) \\
&= \;\; S_3 \text{ (say).}
\end{aligned}
$$

$S_3$ o M gives back $S_3$ so the fuzzy linguistic hidden pattern is a fixed fuzzy linguistic point.

Thus the economic burden is medium, they suffer the related attributes $P_4$, $P_5$, $P_6$, $P_7$ and $P_{11}$ either in the medium state or in the low state.

Suppose Y = (0, 0, 0, ..., high) be the fuzzy linguistic state vector in which only the feelings of the caretakers is high and all other fuzzy linguistic states are in the off state that is the zero state to find the effect of Y on the fuzzy linguistic dynamical system M.

$$
\begin{aligned}
\text{Y o M} \;\; &= \;\; \max \min(Y, M) \\
&\hookrightarrow \;\; (0, 0, 0, 0, 0, \text{low}, \text{low}, 0, 0, 0, \text{high}) \\
&= \;\; Y_1 \text{ (say).}
\end{aligned}
$$

We see $Y_1$ is a fixed point of the dynamical system. Whatever be the problems faced by the caretakers it has only low impact on $P_6$ and $P_7$ is low.

The authors have worked with various states and the conclusions based on this is given in the end of this chapter.

Now we proceed onto use the FLRM [50] model to study the same problem relating to the economic situations of the mentally retarded children.

Let us take for this FLRM model take along the domain space the economic condition of the mentally retarded children / caretakers / parents as $S_1, S_2, S_3, S_4, S_5, S_6$ and $S_7$ where

$S_1$ - very rich

$S_2$ - rich

$S_3$ - upper middle class

$S_4$ - middle class

$S_5$ - lower middle

$S_6$ - poor

$S_7$ - very poor

and $P_1, P_2, \ldots, P_{10}$ described earlier in this chapter are taken as the range space of the fuzzy linguistic dynamical system.

We have taken the same L the fuzzy linguistic set to study the problem using FLRM-model.

We use an expert as a caretaker of the mentally retarded child and the fuzzy linguistic graph given by him is as follows:

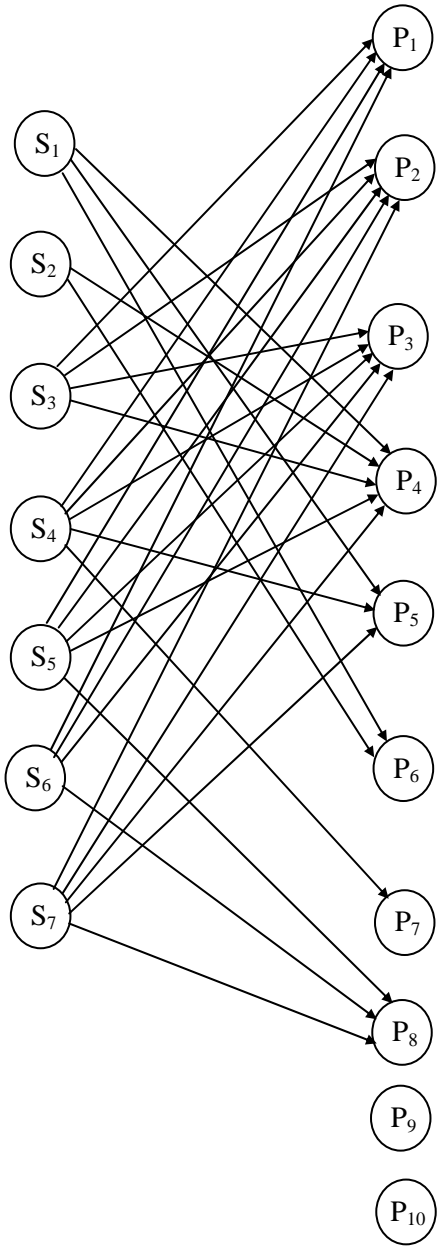

The fuzzy linguistic matrix P associated with the fuzzy linguistic graph is as follows.  P =

|       | $P_1$   | $P_2$ | $P_3$ | $P_4$  | $P_5$  | $P_6$ | $P_7$ | $P_8$   | $P_9$ | $P_{10}$ |
|-------|---------|-------|-------|--------|--------|-------|-------|---------|-------|----------|
| $S_1$ | 0       | 0     | 0     | medium | high   | high  | 0     | 0       | 0     | 0        |
| $S_2$ | 0       | 0     | 0     | high   | 0      | high  | 0     | 0       | 0     | 0        |
| $S_3$ | high    | low   | high  | high   | 0      | 0     | 0     | 0       | 0     | 0        |
| $S_4$ | v.high  | low   | high  | 0      | high   | 0     | high  | 0       | 0     | 0        |
| $S_5$ | high    | high  | high  | high   | 0      | 0     | 0     | high    | 0     | 0        |
| $S_6$ | medium  | high  | high  | 0      | 0      | 0     | 0     | v.high  | 0     | 0        |
| $S_7$ | 0       | high  | low   | high   | v.high | 0     | 0     | v.high  | 0     | 0        |

Now we study the effect of fuzzy linguistic state vectors on the dynamical system using max min operation.  Certainly in this case the fuzzy linguistic hidden pattern will be fixed that is a pair of state vector or a pair of state vectors for the limit cycle. Consider the state vector in which the attribute 'rich' alone is in the on state all other attributes in the domain space are in the off state.  The effect of X = (medium, 0, 0, 0, 0, 0, 0) is the given state vector.

We find max min ({X, S}), the hidden pattern gives medium effect on all the people who have mentally retarded children with no difference in the economic state.  Likewise all the nodes $P_1$, $P_2$, …, $P_8$ have medium effect and $P_9$ and $P_{10}$ are no effect as they do not exist or it is a state of zero state.

Also if social stigma alone is in high state in the range space and all other nodes are in the off state still the hidden pattern shows rich and poor certainly suffer the same amount of social stigma.  Only criteria the rich and very rich do not openly admit they suffer a social stigma but at the back of them the society certainly views them as family with social stigma or a curse which cannot be denied.  That is even if the family with a mentally retarded child do not think as a social stigma but the public / society certainly associates a social stigma and as a cursed family of god for their past karma. This cannot be wiped away for they feel "some thing" is wrong with such families. Only constantly addressing the public by effective media / media persons can eradicate this social stigma.

One has to understand any family can have mentally retarded children. Further once the social stigma is associated it is very difficult in case of middle class and poor to get married their siblings. The relatives also shun such families so marriage in their families at the reproductive age becomes an impossibility and even if they marry they are gripped with a constant fear that their offspring may be a mentally retarded one. Also observations prove the evidence of any government institutions solely for the mentally retarded children. Also one is not in a position to say whether rigorous research to stop the occurrence of such children by genetic modification has been substantially carried out.

## 7.3 Suggestions and Conclusions:

1.  Taking care of type I mental disability is a big challenge to the caretakers. Most of the cases we see these children are left under the care of charity organization. We from our data could see not even a single very poor family in rural Tamil Nadu ever tended or grew these children. Several of them innocently / ignorantly confided they have made mercy killing. Many lamented only after a period of three to ten months only they found that these children were totally disabled (that is both mentally and physically). So they did not care to tend them which also resulted in a natural death, for no proper hygiene or no good food or no proper care, made them an easy victim to all diseases and even at the time of sickness the caretakers had no mind to give them medical aid for they felt up bring such children was not only a social stigma but a economic burden on their meager earnings.

These children were not given bath or change of dress fully prayed by flees and bed sores. So the very poor / poor shunned to bring up these children for poverty, social stigma and the prosperity of the marriage alliances in such families was very difficult. However it is pertinent to make one more observation that such children are very rare in villages or rural area among the poor.

In case of middle class when they identified that the child is of type I mental and physical disability they silently disposed of the child. They never tended if they had no mind to kill it left in dustbins / public places so no one knew the identity of the child.

In case of rich when such children were born some admitted them in charity homes by regularly paying the organizations which took care of them. However the bringing up or birth or existence of such children were never known even to the close relatives of them. For in India the birth of such children is a social stigma and a curse of god. It is realized more as a social stigma. For the parents of these children are looked down by the public as sinners. So they do not claim but help the charity institutions to run so that they visit the child in privacy.

Some of the caretakers take care of these children with love and affection and some for the sake of money. Some of the caretakers said these children cannot be admitted in the hospital for the cost is high and second they will be easily infected by the hospital atmosphere and may be this information will become public.

Thus over all conclusions is that these children who do not know anything live/die at the mercy of their parents / caretakers. We suggest that doctors at the time of birth of such type of children can inform their parents about the pros and cons of the bringing up and government can run a centre purely for these children with constant monitoring of the doctors. Otherwise the study reveals that such children have no hopes for they are really a burden and social stigma to their parents.

Another suggestion is that the medical field should carry out special research on this sort of mental disability so that such children can be aborted in a very short period of confinement. This will do a lot of good to both parents and their children. Serious research in this direction must be made so that this sort of causality can be averted.

Since the mentally retarded children are not vote bankers, the government does not built or even take any steps to built schools for them. Infact they are the least taken care of citizens. It is not known whether the parents themselves have enrolled their name in the voter ID list lamented a care taker.

It is pertinent to record here they are not only neglected but never recognized as humans as they lack the main feature 'mind' or brain or the 'sixth sense' said the experts. So based on this the experts gave the following suggestions.

1. They said at the time of pregnancy the mother should be tested and if the infant has any of the symptoms or signs of down syndrome the child should be aborted.

2. The parents should be given awareness about why such children are born and steps to prevent it by leading a proper life.

3. Government should by means of advertisement, TV programmes and announcements, educate the people using medical experts why such children are born and how to avoid it.

4. Extreme care must be taken by the doctors while delivering the baby for if due to neglect the child has problems it is their moral responsibility.

5. Some of the experts said that just like HIV/AIDS eradication the education for the married couples should be given so that the chances of giving birth to a mentally retarded child is minimized, for the stigma associated with metal retardation is as good as or as bad as the HIV/AIDS.

Unless steps are taken to council the couples about giving birth to mentally retarded children or children with down syndrome or about hyper active children it is not easy to stop the birth of such children.

6. It is important and pertinent to keep in record that this research or study on mental retardation still remains dormant as experiment performed on rats / rabbits / monkeys cannot hold for studying the functioning of the human brains. So only special type of observation and research should be carried out on samples of pregnant women from the date of pregnancy and results based on this must be studied and analysed. After delivery we can study the child's mental ability. This will certainly throw light on the study of the mental retardation.

7. Research should be carried out about down syndrome children and their parents so that some break through can be made and some preventive methods can be carried out.

8. Some study and research can be done by the medical experts so that by spotting the gene present in these children and their parents and means to correct it may be developed by some alternative research or genetic modification.

9. Research must be done on the parents of the mentally retarded children and means and methods should be discovered so that this disability is not carried out to the future generation. If this is achieved certainly without any hesitation the social stigma can be rooted and people would come freely to have marriage proposals from these families.

10. Special training must be given by the government to experts to take care of this children which includes giving psychological counseling to their parents and training these children to be self sufficient.

11. These families which has mentally disabled children and if the economic condition of these families are poor or lower middle class certainly government can give them some money or allowances to their parents to support these children. This will not only lessen the social stigma but also will help the parents to take care of them and they can spend this money on nutritious food, medical care and proper clothing.

12. Several experts claimed it was a gross injustice done by the government to these children as they are not taking any effort to build day care centres for these children. This will lessen the burden of the parents on one hand and give employment to several people who have the service mind to such children. They are neglected by the government as they are not the vote bankers claimed the experts.

13. Another claim by the experts was that most of the mentally retarded children were obese and this also we are not in a position to give any reason only medical experts can substantiate this claim.

14. Also it was very surprising to note that type I, type II and type III mentally retarded children did not know the value / denomination of currency notes. They did not know any basic counting even with fingers. This shows their mind's ability or level.

15. Type I, Type II and Type III mentally retarded children do not have any human sense expect the human appearance claimed the experts.

That is why they do not have voting rights or to be more precise the experts claimed they are unaware of the fact about whether these mentally retarded people had voting rights.

16. Finally it is suggested that rigorous research must be carried out by the medical experts to stop birth of such children. For in case of Type I, Type II and Type III it is impossible to rehabilitate them in any way so it would be good if we can prevent the birth of these children.

# CONCLUSIONS AND SUGGESTIONS

In this chapter we give the overall conclusions and suggestions arrived from the fuzzy analysis and from the collected data through social survey, interviews and discussions with the PWDs caretakers and the other stakeholders. Some of the conclusions and suggestions obtained from the fuzzy models and the fuzzy bimodels are given in the respective chapters.

The main conclusions derived from the fuzzy analysis and the suggestions based on it and the conclusions and suggestions based on the collected data described in chapter two and discussion and interviews with them are also given as supportive evidence of the mathematical analysis.

Special FCM models of both the caretakers and the PWDs clearly showed that they suffer from acute poverty and they attributed their poverty to several problems like poor nutrition, improper clothes, poor health care and marriage a question. However we say they (PWDs and caretakers) never showed any inclinations towards the importance of school education.

However their marriage remained a question mark because most of them in their interviews said that the presence of a PWD in their family was a social stigma which prevented even the marriage of their brothers / sisters who were not suffering from any disability. However from the FRM model it was very clear that the family in which a PWD was born always suffered a social stigma. Infact they are considered to be cursed by god. They attributed the PWDs birth in their family to be due to the supernatural power. So they ill-treated and discriminated the PWDs. In fact they ignored the PWDs by denying the specialized care they require. However throughout the discussions interviews and from the data collected we were very surprised to see the PWDs and their caretakers gave least or to be precise no importance to school education. This is very much evident from the fact.

The node $D_6$ was in the on state in the special FCM given by M in page 58 of chapter four We saw the resultant hidden pattern of the state vector X = (0 0 0 0 0 1 0 0 0 0 0) was X itself. It did not influence any other node. So the hidden pattern of X was invariant on the dynamical system, which is very strange in FCM models. The same was the resultant when the special FCM model of the caretakers was used. Thus we can say school education was given "zero" importance by both the PWDs and the caretakers, which is evident from the two special FCM models described in chapter four. So it is clear that they did not care for school education. However when we analyse the result using the special FCM given by the NGOs we found it different. Using the dynamical system N given by the NGOs we find the resultant of the state vector X = (0 0 0 0 0 1 0 0 0 0 0), that is the node $D_6$ alone in the on state and all other nodes are in the off state.

The hidden pattern of X on N is given by a fixed point P = (0 0 0 0 0 1 1 1 1 1 1), clearly showing due to lack of proper school education, they are not self employed or employed, they do not know any information about SSHG or welfare measures of government never reaches the PWDs as they are ignorant of it, they suffer from poor economy and their marriage remains a question mark. Thus we see from our fuzzy analysis the lack of school education in a PWD in a root cause of poverty, unemployment, no marriage, ignorant of SSHG and government rehabilitation measures to them. So in the first place the PWDs and the caretakers must be given proper counseling about the need of school education. For the on state of the node school education in both the special FCMs models of the 82 PWDs and 82 caretakers did not influence the other node. The hidden pattern was the same state vector, which was sent. Thus it is confirmed that they are least bothered about school or formal education.

However they fear the social stigma when the teachers and other students talk or laugh at the disability of their children. In the first place they fear that the society will view them as a family cursed by god and that is why a disabled child is born in their family. This unscientific superstitious belief is strong in people's mind due to the unshakeable faith they have in the theory of 'karma' of Hindu religion that dictates that one's present life is the resultant of the way in which he/she lived in the previous life.

1. It is observed from the document of the World Bank that there is a considerable fall in the PWDs working age employed in the early 1990 and in the early 2000 among the literate PWDs. Only in case of diploma / certificate course qualified employed PWDs in this period satisfactory. In all other cases there is a considerable fall say atleast up to a minimum of 10% in this time period from early 1990 to early 2000. This may be one of the causes that can attributed to the younger generation remaining unmarried. From analysis of the data we found they (PWDs) crave for an employment. They are not interested in school education. They prefer diploma degree / certificate course in certain job fetching skills to earn. This is also evident that 0% employment in diploma / certificate holders in early 1900 but in early 2000 60% were employed from the diploma / certificate course holders. This is also evident from May 2007, document of the World Bank. Certainly study shows that among the employment of literates is reduced in these 10 years only confirms all rules to reserve post for them or any other considerations have not literally reached them. If this is the case with every PWD one can imagine the position of the poor rural

PWDs. All these confirm they suffer from an acute social stigma so do not opt for employment in private or quasi government they prefer to be self-employed.

2.   Further it is observed that from the special FCM and FRM analysis that neither the PWDs nor the caretakers are interested in school education. This is clearly evident from the tables given in chapter two where 59% of the PWDs have never even entered the school premises. This is all in keeping with the fact that in these 10 years the literate employed PWDs were considerably reduced. However from our discussions with PWD it is clear that they do not get any unemployment allowances. We are not wrong if we say none of the government aid has reached these poor rural PWDs.

3.   The disinterest in the caretakers as well as PWDs in knowing about and applying for benefits through the government's rehabilitation measures and jobs among the rural poor is due to:

> i.    Lengthy and complex procedure
> ii.   Distance (problem of mobility not trained on the cause and causalities of disabilities and the rights of PWDs)
> iii.  Non co operative officials
> iv.   Lack of transport
> v.    Lack of information
> vi.   Poor communication of rules
> vii.  Physical barriers.

All these problems can be easily over come if alone government takes proper census of the PWDs in each state using the unemployed youth after sufficiently training them and render their services so that both services / benefits reach the PWDs.

This has two benefits

> i.    Some youth who are unemployed get employment
> ii.   The poor PWDs who are from rural area can get benefit.

4. It is further suggested just like the government has a health care center in every village to benefit the rural poor, the government can run a PWD healthcare center after identifying the reasonable number of PWDs in and around that place with the following aids:

  i. Proper transport to help the PWD to visit them and avail services periodically.

  ii. Doctor to take care of the health of the PWDs.

  iii. The PWDs should be advised to seek the support of the NGOs and the government for employment / self-employment. In case of very several disabilities this unit can mobilize some monetary help every month to maintain themselves in case of rural poor.

  iv. As they are least interested in formal school education these units can train them in some skills depending on their disability and facilitates appropriate employment only after they complete the school education.

  v. In case of rural poor PWDs., the government can distribute new cloths for PWDs during festival occasion and also can distribute sweets and free food.

5. So one of the powerful ways is to make the public aware of the fact disability is not due to karma but only due to some factors like

  i. Malnutrition

  ii. Marrying within the relatives (consanguinity)

  iii. Early or late marriage.

  iv. Lack of hospitals to administer institutional delivery.

  v. Neglect of periodical health check-up and follow up which has resulted in disability etc.

  vi. Not proper medical attention at time of delivery of the child.

  vii. Polio attack

  viii. Accidents

ix.  Some disease which can cause disability like brain fever, fits when they are children.

So those families should be given moral support to over come their mental trauma by the society and the society should not shun them.

6. From our study using these special FCM model, FRM model and the special FCRM bimodel it is clear that the PWDs and the caretakers shun the usual schools attended by others. This was also evident from the fall in the educated PWDs employed in the 10 years from early 1990 to early 2000.

This is a major issue and should not be taken lightly for education is the main sign of development. They (PWDs) have the right to get education regardless of their physical, intellectual, social, emotional status as every child has a right to education.

It is shocking to know the PWDs are ill treated both by the students (who are not PWDs) and the parents for they fear that their child may be affected by constantly looking at the PWDs. Teachers fear that they may give birth to PWDs as they constantly see them. Further the PWDs caretakers want to hide the existence of a PWD in their home. That is why they do not wish to send him / her to school or public places.

So in the first place the public, students, parents and teachers must be given education to create an attitudinal change in them. PWDs and their caretakers must also be given counseling to make their disabled children to attend regular school. The recent Community Based Rehabilitation (CBR) approach calls for such a collective time bound effort in realizing the rights of the rural PWDs.

Both the children with disability (CWDs) and their caretakers must be given appropriate incentives so that the children are not stopped from attending school. Only if the children with disability even after say a time span of ten years do not show any sign of formal learning in their classes then alone can they are given special training in self-employment depending on their disability.

Parents often blame teachers and teachers blame parents and by this blame game the children with disability (CWDs) becomes a school drop out or never enters school. So it is advised that first government insists the parents of children with disability to send their children to preschools and failing to do it, the parents should be punished legally. Unless parents work together collaboratively with government and teachers it is impossible for the CWDs to enjoy their rights or lead a life on par with other children. It is impossible to make CWDs live with dignity until the parents get involved in the learning process. Parents have no right to feel or say they have no money or no time or feel ashamed of the child and consequently neglect its education. Government should also take special steps to support the rural poor caretakers of these disabled children so that they can cooperate with teachers and government and educate their children. The parents should not play only the passive role of just bringing the children to schools. They should take special interest in their education than the interest they take in case of their normal children. At this juncture it is pertinent to say social workers must be trained and employed by the government to council the caretakers and insist on the education of their disabled children. Thus the education of disabled child can be achieved only by collaborative groups – which include apart from teachers and parents educational psychologists, social welfare workers, medical doctors, nurses, physiotherapists and so on. These professionals have diverse roles-ranging from quality monitoring, whole school development, staff student support to meet the needs of the disabled children. Only in case of children with disability who suffer from critical / serious problems and cannot be given usual education can be given learner support and training in special education.

It is suggested some six months training to teachers may be given to equip the teacher with the theory and practice in special education and a specialization in one of the areas.

These teachers who are willing not only to undergo special training but also willing to work for the disabled children may be paid an extra allowance every month. This will be an incentive for them to work for the special children. However the progress of these special children should be regularly and constantly monitored. For that alone can give the expected results. We are however not discussing in this thesis the gender issues.

With the development made in science and technology it is the duty of the nation to provide the rural poor disabled children freely get the learning materials, assisting devices (like hearing aids, wheelchairs, calipers, white canes etc.). Even if these are not done properly it is a disgrace for us to call ourselves a nation with spiritualism and so on.

7.  Further it is the duty of everyone to respect the disabled or treat them in par with others. When we say treat them on par with others, we mean they should not be laughed at or mocked at or avoided or discriminated for fear of giving birth to disabled children. If karma is the cause of disability then how will anyone who sees a disabled give birth to a disabled?

The Hindu dharma speaks of karma barbarically describing the PWDs and their caretakers as persons who have done bad deeds in previous lives and as a consequent curse they suffer in this life. If this stigma is to be wiped out the only solution the experts feel is wiping out of the concept of karma in the Hindu religion which is not only barbaric in its out look but unscientific in its practice. It is a pity that most of the rural poor are so badly trained by the Hindu religion they feel it is a curse and torture the PWDs or ill treated them and discriminate them. So if the PWDs are to be treated humanly this deep rooted concept of karma or curse should be wiped out from the minds of the poor caretakers. The rich upper caste man however does not feel his child who is disabled as a curse. His richness helps his disabled child to be happy and enjoy all sorts of equality. The poor rural must in the first place be sensitized, educated and made to realize that disability is not the result of supernatural intervention. Only special study and research can change the mindset of the poor rural Hindus. That is why we say eradication of the concept of karma in the Hindu religion is the only solution.

8. It is unfortunate to state that in India only CWDs from economically poor or the middle class who mainly live in rural areas are denied every form of support. They are not only ill treated and discriminated but treated worse than animals. The main reasons they give is poverty and social stigma. Government with the help from professionals should take proper steps to make these PWDs to live with self-esteem and feel that they are infact an integral part of the society. When this is achieved the PWD will not suffer any more from frustration, depression or inferiority complex or discrimination. After all advancement if it is real must be from barbarism to humanism and not be inhuman especially with the PWDs.

For this government should come up with time bound programs that are properly backed by policy and budgets. Most initiatives in developing countries particularly on disability are left to charity that has destroyed even very brilliant programs since charity has no legal commitment and it can be with drawn any time depending on religion, caste, place, language, colour and gender and the popularity / profit they derive from it.

Though PWD Act of 1995 declared free education for disabled children upto the age of 18, we see the rural poor caretakers and PWDs are unaware of it. So in the first place government should give programs on TVs and radios and in public places where the rural poor meet and make the rural poor caretakers and PWDs to become aware of their rights. Social development workers can also be employed to carry out this task.

It is suggested that higher education for rural poor PWDs can be absolutely made free with free accommodation in hostels. Unless this is made with some (monetary help) stipend given monthly to them there will only be a decline in the already small percentage of PWDs who get higher education.

9. The other causes for denying school education by the caretakers are:

    i.   Caretakers fear the open social stigma they suffer for having a PWD child.

    ii.   The PWDs are laughed at by their classmates and school mates which invariably hurts them to the core and stops them from perusing their studies.

    iii.   The usual comments by teachers, public, students are that even people with proper functioning don't get opportunities so they look down upon the PWDs. This is evident from table 2.4 given in chapter two that majority of school dropouts are PWDs. Thus it would be appropriate to offer them short term training for 45 to 60 days after their school education and in case they do not wish to pursue higher education to be self sufficient and introduce them to some type of self employment with this opportunity they (caretakers) will not be stopping their PWDs from education.

This sort of training can be given to PWDs by the unemployed graduates who are willing and have a mind set to work for PWD they (unemployed graduates) can first be trained for 3 month with specialization in one of the disabilities. They can also be given power to visit schools government and private offices to see whether the stipulated percentage of PWDs are given employment and train the students after they complete their school final. Government must take steps with the help of these trained graduates. This will not only channelize the energy of these youth but also the nation can improve the man power. Under the guidance of these unemployed trained youth the PWDs can be given self employment depending on the type of the disability, they suffer from. From the data we collected majority of them wanted to be self employed even without school education.

This will not only make the PWDs happy but also they can be contributing to the development of the nation. Once they are self-sufficient economically their marriage will not remain a question. Certainly they can get married at an appropriate age and lead a normal non discriminated life with dignity. Since school education is of no importance to the PWDs and caretakers, we have to take special steps to bring them under the compulsory education.

India being second largest in population in the world having one-thirds of the worlds poor people. Twenty percent of the worlds out of school children live in India and these children have little access to medical care or economic security. Stark inequality also describes the situation for disabled people whose illiteracy has been estimated at a staggering 55%. Unless the disabled are integrated with the general community as equal partners the normal growth in PWDs to face life with courage and confidence will remain a question mark. The first major intervention should be enrolment and retention of all children with disability in the main stream of education system. .

10. The hidden pattern from the FRM of the public given in chapter five of this thesis clearly showed that poverty of the caretakers is not the sole cause for neglecting PWDs or discriminating PWDs. For, if one has the heart, what is shared with others can also be extended and shared with the PWDs.

We have tried to avoid offering judgments on the relative effectiveness of the different works but have read them with sincerity. Thus as far as India is concerned, the disability by poor rural Indians are attributed to Karma theory. Karma is a Sanskrit word means action. According to this theory, privileges and deficits of the current life are to be attributes to the actions of the past lives.

11. If one has any kind of impairment, it is attributes to the actions of the past lives. Likewise if one enjoys a privileged life it is attributed to the good deeds of his / her previous life.

This fatalism and linkage of disability to karma is found in several researches of disability in India (14 Gabel). Thus for a vast majority of people living on the Indian subcontinent disability thus is irrevocable, since the cause is believed supernatural. While the disabled were objects of pity and sympathy, prevention was considered unthinkable and rehabilitation not possible. Families of disabled persons also resigned themselves to their fate and suffered in silence.

Thus in our thesis it is suggested from the fuzzy models education can alone wipe out this mind set up. Also it would be proper if the concept of Karma is wiped out from the Hindu religion. When Hindu religion is throwing away several inappropriate concepts in present time why not the concept of Karma?

Secondly it is essential to state that India's urban centers now realize that disability is not the result of super natural intervention because the urban people are better educated than the rural poor. It is only in remote rural villages illogical unscientific belief systems are still cherished and practiced.

Thus it is understood that the rural makers have to be sensitized on the cause and causalities of disability. Then only they will realize that the PWDs have the right to special care and they too can demand all the basic rights like any other normal individuals to live a life of dignity and peace.

# INDEX

**A**

Abicyclic, 36-7
Adjacency bimatrix, 35-6

**B**

Biconcepts, 35
Bimatrix, 35
Binodes, 35
Bivectors, 35

**D**

Directed bicycle, 36-7
Domain fuzzy linguistic space, 121-3
Dynamical bisystem, 37

**E**

Edge biweight, 36
Equilibrium bistate, 37-8

## F

Feedback, 37
Fuzzy Cognitive Maps (FCMs), 11
Fuzzy Cognitive Relational Maps (FCRMs) bimodel, 11,29-33,
99
Fuzzy Linguistic Cognitive Maps (FLCM), 12
Fuzzy linguistic dynamical system, 120-3
Fuzzy Linguistic graphs, 12, 120-2
Fuzzy Linguistic hidden pattern, 12, 120-3
Fuzzy Linguistic matrices, 12, 120-3
Fuzzy Linguistic Relational Maps (FLRM), 120-2
Fuzzy Linguistic terms, 12, 120-2
Fuzzy Relational Maps (FRMs), 11
Fuzzy state linguistic vector, 120-3

## H

Hidden bipattern, 37

## I

Instantaneous state bivector, 36

## L

Limit bicycle, 38-8

## N

Non Government Organizations (NGOs), 13

## P

People with Disability (PWD), 13
Positive Special connection matrix, 49-51

# ABOUT THE AUTHORS

**Dr.W.B.Vasantha Kandasamy** is an Associate Professor in the Department of Mathematics, Indian Institute of Technology Madras, Chennai. In the past decade she has guided 13 Ph.D. scholars in the different fields of non-associative algebras, algebraic coding theory, transportation theory, fuzzy groups, and applications of fuzzy theory of the problems faced in chemical industries and cement industries. She has to her credit 646 research papers. She has guided over 68 M.Sc. and M.Tech. projects. She has worked in collaboration projects with the Indian Space Research Organization and with the Tamil Nadu State AIDS Control Society. She is presently working on a research project funded by the Board of Research in Nuclear Sciences, Government of India. This is her 67[th] book.

On India's 60th Independence Day, Dr.Vasantha was conferred the Kalpana Chawla Award for Courage and Daring Enterprise by the State Government of Tamil Nadu in recognition of her sustained fight for social justice in the Indian Institute of Technology (IIT) Madras and for her contribution to mathematics. The award, instituted in the memory of Indian-American astronaut Kalpana Chawla who died aboard Space Shuttle Columbia, carried a cash prize of five lakh rupees (the highest prize-money for any Indian award) and a gold medal.

She can be contacted at vasanthakandasamy@gmail.com
Web Site: http://mat.iitm.ac.in/home/wbv/public_html/
or http://www.vasantha.in

---

**Dr. Florentin Smarandache** is a Professor of Mathematics at the University of New Mexico in USA. He published over 75 books and 200 articles and notes in mathematics, physics, philosophy, psychology, rebus, literature. In mathematics his research is in number theory, non-Euclidean geometry, synthetic geometry, algebraic structures, statistics, neutrosophic logic and set (generalizations of fuzzy logic and set respectively), neutrosophic probability (generalization of classical and imprecise probability). Also, small contributions to nuclear and particle physics, information fusion, neutrosophy (a generalization of dialectics), law of sensations and stimuli, etc. He got the 2010 Telesio-Galilei Academy of Science Gold Medal, Adjunct Professor (equivalent to Doctor Honoris Causa) of Beijing Jiaotong University in 2011, and 2011 Romanian Academy Award for Technical Science (the highest in the country). Dr. W. B. Vasantha Kandasamy and Dr. Florentin Smarandache got the 2011 New Mexico Book Award for Algebraic Structures. He can be contacted at smarand@unm.edu

**Dr. A.Praveen Prakash** has obtained his doctorate degree on Problems faced by PWDs from University of Madras in 2009. He is presently working as Dean in Hindustan College of Engineering, Chennai. He has worked with several NGOs on projects connected with disability. He has 34 years of teaching and research experience. He has published several research articles. He can be contacted at apraveenprakash@gmail.com